\documentclass[11pt]{amsart}
\usepackage{amsmath,amsthm, amscd, amssymb, amsfonts, mathrsfs,mathtools,pb-diagram,pb-xy}
\usepackage{amsxtra, color}
\usepackage{enumitem}
\usepackage[all]{xy}
\usepackage{amssymb}
\usepackage{hyperref}
\usepackage[T2A,T1]{fontenc}
\usepackage[ansinew]{inputenc}
\usepackage{graphicx,fancyhdr}

\usepackage[table]{xcolor}   
\usepackage{booktabs}

\usepackage{tikz}
\usetikzlibrary{matrix,arrows.meta,decorations.pathmorphing,calc,positioning,fit,backgrounds}
\definecolor{accent}{HTML}{1F4E79}

\usepackage{tikz-cd}

\usepackage{booktabs}
\allowdisplaybreaks

\newcommand{\supera}[2]{{\mathbf A}(#1|#2)}
\newcommand{\superb}[2]{{\mathbf B}(#1|#2)}

\newcommand{\superd}[2]{{\mathbf D}(#1|#2)}
\newcommand{\superda}[1]{{\mathbf D}(2,1;#1)}
\newcommand{\superf}{{\mathbf F}(4)}
\newcommand{\superg}{{\mathbf G}(3)}

\newcommand{\wtoba}{\widetilde{\mathcal{B}}}
\newcommand{\cP}{\mathcal{P}}

\newcommand{\bracket}[3]{\binom{K_{#1};#2}{#3}}
\newcommand{\bra}[2]{\binom{K_{#1}}{#2}}

\newcommand{\lu}{\mathcal{L}}
\newcommand{\luq}{\lu_{\bq}}
\newcommand{\luqm}{\lu_{\bq^{(-1)}}}

\newcommand{\I}{\mathbb{I}}
\newcommand{\N}{\mathbb{N}}
\newcommand{\bq}{\mathfrak{q}}

\newcommand{\h}{\mathfrak{h}}
\newcommand{\g}{\mathfrak{g}}
\newcommand{\n}{\mathfrak{n}}

\newcommand{\bp}{\mathfrak{p}}
\newcommand{\bpf}{\mathbf p}

\newcommand{\eb}{\mathbf{e}}
\newcommand{\hb}{\mathbf{h}}
\newcommand{\fb}{\mathbf{f}}
\newcommand{\hbt}{\widetilde{\mathbf{h}}}

\newcommand{\mund}{\underline{m}}
\newcommand{\nund}{\underline{n}}
\newcommand{\Lt}{\mathtt{L}}

\numberwithin{equation}{section}
\theoremstyle{plain}

\newtheorem{theorem}{Theorem}[section]
\newtheorem{lemma}[theorem]{Lemma}
\newtheorem{definition-theorem}[theorem]{Definition-Theorem}
\newtheorem{proposition}[theorem]{Proposition}

\newtheorem{question}[theorem]{Question}
\newtheorem{corollary}[theorem]{Corollary}

\theoremstyle{definition}
\newtheorem{definition}[theorem]{Definition}
\newtheorem{example}[theorem]{Example}

\newtheorem{remark}[theorem]{Remark}
\newtheorem*{remark*}{Remark}
\newtheorem{remarks}[theorem]{Remarks}
\newtheorem{conjecture}[theorem]{Conjecture}

\newcommand\bth[1] { \begin{theorem}\label{t#1} }
\newcommand\ble[1] { \begin{lemma}\label{l#1} }

\newcommand\bpr[1] { \begin{proposition}\label{p#1} }
\newcommand\bqu[1] { \begin{question}\label{q#1} }
\newcommand\bco[1] { \begin{corollary}\label{c#1} }
\newcommand\bde[1] { \begin{definition}\label{d#1}\rm }
\newcommand\bex[1] { \begin{example}\label{e#1}\rm }
\newcommand\bre[1] { \begin{remark}\label{r#1}\rm }
\newcommand\bres[1] { \begin{remarks}\label{r#1}\rm }
\newcommand\bcj[1] { \begin{conjecture}\label{j#1}\rm }
\renewcommand {\eth} { \end{theorem} }
\newcommand{\ele} { \end{lemma} }

\newcommand{\epr} { \end{proposition} }
\newcommand{\equ} {\end{question} }
\newcommand{\eco} { \end{corollary} }
\newcommand{\ede} { \end{definition} }
\newcommand{\eex} { \end{example} }
\newcommand{\ere} { \end{remark} }
\newcommand{\eres} { \end{remarks} }
\newcommand{\ecj} { \end{conjecture} }
\newcommand{\enota} { \end{notation} }

\newcommand\pf{\begin{proof}}
\newcommand\epf{\end{proof}}

\newcommand{\Cset} {\mathbb C}
\newcommand{\Zset} {\mathbb Z}

\newcommand{\Z}{\mathbb Z}

\newcommand{\Qset} {\mathbb Q}

\newcommand{\ku}{\mathbb C}

\newcommand{\toba}{{\mathcal B}}

\newcommand{\fO}{\mathfrak O}

\newcommand{\II} {{\mathcal{I}}}

\newcommand{\JJ} {{\mathcal{J}}}

\newcommand{\Ss}{\mathcal{S}}

\newcommand\cW{\mathcal{W}}

\newcommand{\cX}{\mathcal{X}}

\newcommand{\ee}{\eta}

\newcommand{\wfO}{\underline{\mathfrak O}}
\newcommand{\wbeta}{\underline{\beta}}
\newcommand{\walpha}{\underline{\alpha}}
\newcommand{\wgamma}{\underline{\gamma}}

\newcommand{\uDe}{\underline{\Delta}}

\newcommand{\vpi}{\varpi}

\newcommand{\id} {\operatorname{id}}

\newcommand\Acf {\mathcal{A}}

\newcommand{\bqf}{\mathbf{q}}
\newcommand{\qf}{\mathbf{q}}

\newcommand{\Ab}{\mathbf{A}}
\newcommand{\Bb}{\mathbf{B}}

\newcommand{\Db}{\mathbf{D}}

\newcommand{\Fb}{\mathbf{F}}
\newcommand{\Gb}{\mathbf{G}}
\newcommand{\Ub}{\mathbf{U}}

\newcommand{\Ib}{\mathbf{I}}
\newcommand{\uf}{\mathfrak{u}}

\newcommand{\opp} {\operatorname{op}} 
\newcommand{\lcm}  {\operatorname{lcm}}

\newcommand{\ord} {\operatorname{ord}} 
\newcommand{\ad}{\operatorname{ad}}

\renewcommand {\max} {\operatorname{max}}

\begin{document}
\title[Super quantum groups of divided powers \`a la Lusztig]
{Super quantum groups of divided powers \`a la Lusztig}

\address{\noindent FaMAF-CIEM (CONICET) \\ 
Universidad Nacional de C\'ordoba \\
Medina Allende s/n, Ciudad Universitaria \\
5000 C\'ordoba \\
Rep\'ublica Argentina
}
\author[Iv\'an Angiono]{Iv\'an Angiono}
\email{ivan.angiono@unc.edu.ar}

\author[Cristian Vay]{Cristian Vay}
\email{cristian.vay@unc.edu.ar}



\thanks{
The work of I. A. and C. V. was partially supported by CONICET PIP 11220200102916CO and 11220250100219CO, FONCYT PICT 2020-SERIEA-02847 and SECYT-UNC 33620230100517CB.}

\begin{abstract}
We construct super analogues of the quantum groups of divided powers introduced by Lusztig. We call them super analogues because the Lie-theoretical data governing their structure is comparable to that of contragredient Lie superalgebras. Indeed, we show that they admit PBW-type bases parametrized by root systems of super type, and that they also have Lusztig isomorphisms giving a representation of the corresponding Weyl groupoids. Moreover, we prove that our super objects are equipped with a Frobenius-type morphism whose kernel is the Drinfeld double of a finite-dimensional Nichols algebra of super type, a super analogue of Lusztig's small quantum group, and whose cokernel is the enveloping algebra of the semisimple Lie algebra determined by the Cartan roots of the Nichols algebra. We carry out our construction within the framework of Nichols algebras of diagonal type, specifically, those whose braiding is obtained by specializing at a root of unity. This allows us to consider a divided power integral form inside a large quantum group of De Concini--Kac--Procesi type and then, via specialization, obtain our desired object.
\end{abstract}

\maketitle

%

\section{Introduction}
\label{intro}

Since their appearance in the 1980s, quantum groups have spread across various fields of mathematics and theoretical physics. These Hopf algebras were conceived as deformations of the enveloping algebras $U(\g)$ of semisimple Lie algebras $\g$. They are defined over the field $\mathbb{C}(q)$ of rational functions in one variable and admit various forms over $\mathbb{C}[q^{\pm1}]$, or even over $\Z[q^{\pm1}]$, giving rise, upon specialization, to Hopf algebras over the complex numbers that exhibit different behaviors with regard to their representations. When $q$ is evaluated at a complex number $\xi$ which is not a root of unity, the algebras obtained behave in a quite similar way, resembling $U(\g)$. The story is quite different when $\xi$ is a root of unity. There, three objects stand out.

\begin{enumerate}[leftmargin=*,label=\rm{(\roman*)}]
\item\label{item:intro-dCKP} The De Concini--Kac--Procesi form $U_{\xi}(\g)$, introduced and studied in \cite{DK,DKP,DP}: a version of the Drinfeld--Jimbo algebra obtained by specializing the $\Z[q^{\pm1}]$-subalgebra generated by the $E_i$, $F_i$, $K_i^{\pm1}$. It is a finite module over a central  Hopf subalgebra, which endows it with a structure of Poisson order and makes its representation theory amenable to geometric methods.

\item\label{item:intro-small-qg} The small quantum group, or Frobenius--Lusztig kernel, $\mathfrak{u}_{\xi}(\g)$ \cite{L-JAMS,L-book}: the finite-dimensional quotient of the previous algebra by the augmentation ideal of that central Hopf subalgebra.

\item\label{item:intro-Lusztig} The divided power form $\Ub_{\xi}(\g)$ \cite{L-quantum,L-JAMS,L-book}: the specialization at $q = \xi$ of the $\Z[q^{\pm1}]$-subalgebra generated by the divided powers $E_i^{(n)}$, $F_i^{(n)}$ and the torus. It is neither finite-dimensional nor a deformation of $U(\g)$. Rather, it interpolates between the small quantum group and $U(\g)$, as it fits into an extension of Hopf algebras
\begin{equation}\label{eq:intro-lusztig-sequence}
	\begin{tikzcd}[column sep=large, arrow style=tikz, >=To]
		\mathfrak{u}_{\xi}(\g) \arrow[r, hook] & \Ub_{\xi}(\g) \arrow[r, two heads, "\operatorname{Fr}"] & U(\g),
	\end{tikzcd}
\end{equation}
where the surjection is given by Lusztig's quantum Frobenius map; see \cite{A-notes} for a rigorous proof. This extension is the algebraic counterpart of the Frobenius morphism of algebraic groups in positive characteristic, and lies at the root of the applications of quantum groups to modular representation theory.
\end{enumerate}

The analogues of the first two algebras for Lie superalgebras are already known. In the classical setting, these algebras are defined by explicit presentations in terms of generators and relations, but their generalization to the super setting is quite delicate, as the defining ideal requires higher-degree Serre-type relations whose form depends on the Dynkin diagram. Finding a complete and minimal set of such relations was a substantial task, carried out by Yamane \cite{Ya}. Moreover, the same superalgebra may admit non-equivalent Dynkin diagrams related by odd reflections rather than by the action of a Weyl group, so there is no distinguished presentation.

Nevertheless, these two objects can be treated more uniformly within the theory of Nichols algebras of diagonal type. Following the work of the first author \cite{Ang-dpn}, super quantum groups \`a la De Concini--Kac--Procesi can be constructed by starting with a distinguished pre-Nichols algebra of diagonal type and then applying the Drinfeld double construction. These quantum groups admit a central Hopf subalgebra, and quotienting by this subalgebra yields a small super quantum group. The resulting small super quantum groups are also Drinfeld doubles of the corresponding Nichols algebras.

\medbreak

The objective of this work is to construct an analogue of the quantum group of divided powers and of the exact sequence \eqref{eq:intro-lusztig-sequence} for simple Lie superalgebras, following the Nichols algebra approach used for the other two algebras. 

\medbreak

Nichols algebras of diagonal type are fundamental invariants in the classification program of pointed Hopf algebras with abelian group of group-like elements \cite{AS}. The finite-dimensional ones were classified by Heckenberger \cite{H-classif RS} and can be organized into the classes of Cartan type and super type, the latter consisting of those associated with contragredient Lie superalgebras of classical type, together with finitely many exceptional examples. The structure of these algebras is governed by root systems and Weyl groupoids \cite{H-Weyl gpd,HY-groupoid}, as is the case for contragredient Lie superalgebras. These structures underlie the distinction between the aforementioned classes. Moreover, the associated Drinfeld doubles are related by Lusztig isomorphisms \cite{H-Lusztig}, which are reminiscent of Weyl group isomorphisms and odd reflections in the theory of Lie superalgebras.

\begin{figure}[ht]
\centering
\begin{tikzpicture}[
>=To,
every node/.style={font=\small},
zig/.style={->, decorate, black!70,
decoration={snake, amplitude=0.35mm, segment length=2mm, post length=1.2mm}},
down/.style={->, thick, black!70},
monoL/.style={{Hooks[right]}->, black!70},
monoR/.style={{Hooks[left]}->, black!70},
epi/.style={-{To[length=1.2mm]To[length=1.2mm]}, black!70},
lab/.style={font=\scriptsize, text=black!65},
col/.style={font=\footnotesize\itshape, text=black!70, align=center, anchor=south},
grp/.style={draw=accent!55, fill=accent!6, rounded corners=4pt,
inner xsep=3mm, inner ysep=7mm},
big/.style={font=\large}
]

\matrix[matrix of math nodes, ampersand replacement=\&,
nodes={anchor=center, inner sep=3pt},
column sep=8mm, row sep=30mm] (m) {
|(qf)| \bqf\in\mathbb{C}[\nu^{\pm1}]^{n\times n}
\& |(Bf)| \toba_{\bqf}
\& |(Uf)| U_{\bqf}
\&[12mm] |(Uzf)| \Ub_{\bqf,\Acf} \\
|(q)| \bq\in\mathbb{C}^{n\times n}
\& |(B)| \toba_{\bq}
\& |[big](u)| \mathfrak{u}_{\bq}
\& |[big](Uz)| \textcolor{accent}{\Ub_{\bq}} \\
};

\path let \p1=($(Uz.west)-(u.east)$), \n1={veclen(\x1,\y1)} in
node[big, right=\n1 of Uz] (Ug) {$U(\g_{\bq})$};

\begin{scope}[on background layer]
\node[grp, fit=(u)(Uz)(Ug)] (grp) {};
\end{scope}
\node[col, text=accent, anchor=south east, inner sep=1pt]
at ([xshift=-2mm,yshift=0.3mm]Uzf.south |- grp.north) {Quantum group of};
\node[col, text=accent, anchor=north east, inner sep=1pt]
at ([xshift=-2mm,yshift=-0.3mm]Uzf.south |- grp.north) {divided powers};

\draw[zig]   (qf) -- (Bf);
\draw[zig]   (Bf) -- (Uf);
\draw[monoR] (Uzf) -- (Uf);
\draw[zig]   (q) -- (B);
\draw[zig]   (B) -- (u);
\draw[monoL, thick] (u) -- (Uz);
\draw[epi, thick] (Uz) -- node[above, lab] {$\operatorname{Fr}$} (Ug);

\draw[down] ([yshift=-1.5mm]qf.south |- Bf.south)
-- node[midway, lab, rotate=-90, anchor=south, inner sep=3pt] {$\nu\mapsto\xi$}
([yshift=3.5mm]qf.south |- B.north);
\draw[down] ([yshift=-1.5mm]Bf.south)
-- node[midway, lab, rotate=-90, anchor=south, inner sep=3pt] {$\nu\mapsto\xi$} ([yshift=3.5mm]Bf.south |- B.north);
\draw[down, accent] ([yshift=-1.5mm]Uzf.south |- Bf.south)
-- node[pos=0.35, right, lab] {$\otimes_{\Acf}\mathbb{C}\,/\,\Ib_{\bq}$}
([yshift=3.5mm]Uzf.south |- B.north);

\coordinate (hdr) at ([yshift=2mm]qf.north);
\node[col] at (qf.center  |- hdr) {braiding\strut\\matrix\strut};
\node[col] at (Bf.center  |- hdr) {Nichols\strut\\algebra\strut};
\node[col] at (Uf.center  |- hdr) {Drinfeld\strut\\double\strut};
\node[col] at (Uzf.center |- hdr) {Divided power\strut\\integral form\strut};

\end{tikzpicture}
\caption{Construction of the quantum group of divided powers of $\bq$.}
\label{fig:intro-diagram}
\end{figure}
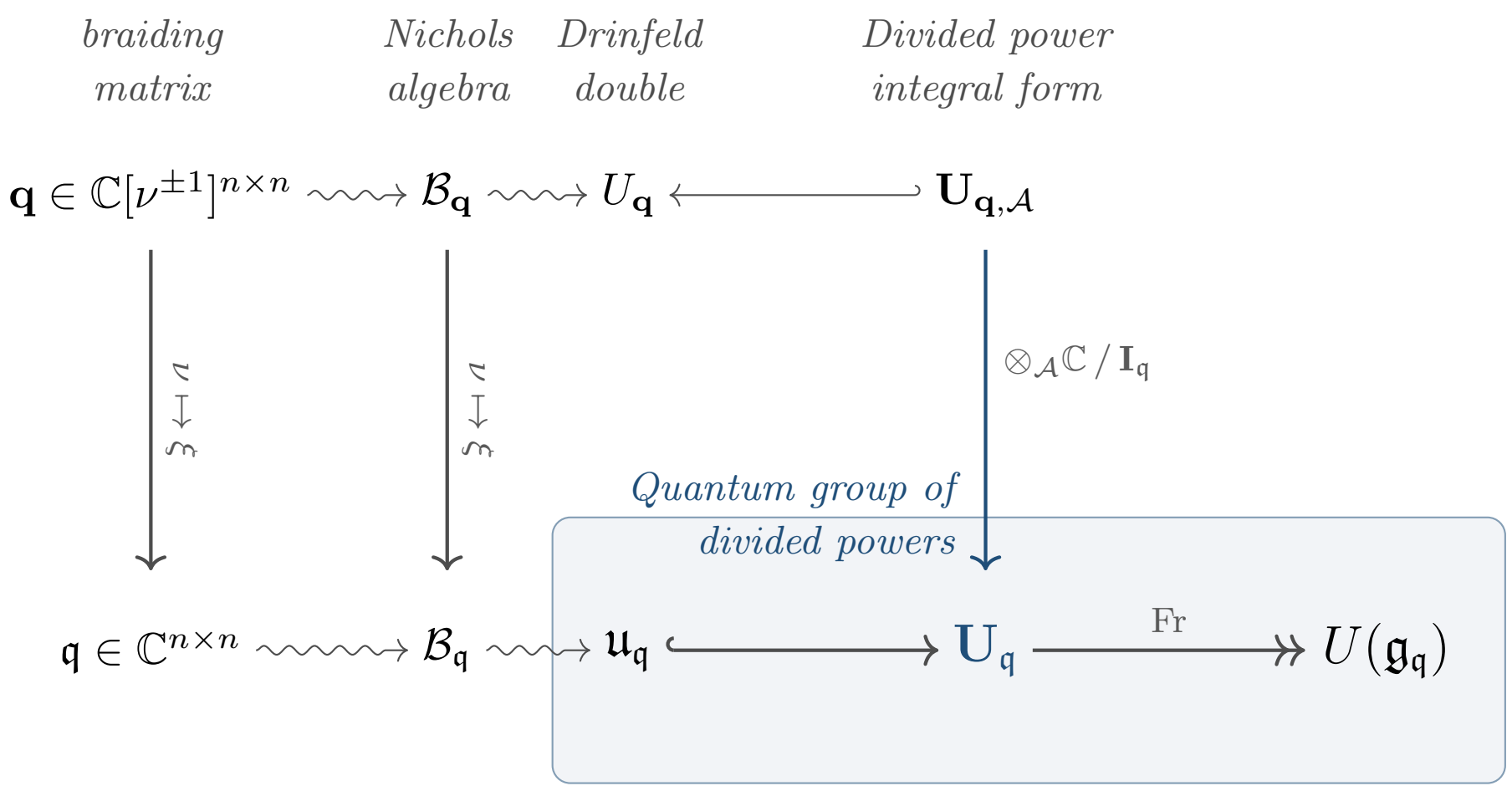

Here, we restrict our attention to Nichols algebras of diagonal type determined by braiding matrices belonging to one-parameter families in the sense of \cite{AAY}. These include those of Cartan type and super type, as well as the exceptional examples $\mathfrak{wk}(4)$ and $\mathfrak{br}(2)$ \cite{AA-diag-survey}. By definition, a braiding matrix $\bq\in\mathbb{C}^{n\times n}$ belonging to one of these families is obtained by specializing a braiding matrix $\bqf\in\mathbb{C}[\nu^{\pm1}]^{n\times n}$ at a root of unity $\xi$, that is, $\bqf(\xi)=\bq$. Notably, the corresponding Nichols algebras $\toba_{\bq}$ and $\toba_{\bqf}$ share the same root system. This framework allows us to deal with the specialization procedure needed for our construction. We thus achieve our objective for this class of Nichols algebras as we explain below, cf. Figure \ref{fig:intro-diagram}. 

\medbreak

We begin by considering the Drinfeld double $U_{\bqf}$ associated with $\toba_{\bqf}\#\ku\Z^{n}$. It is determined by a Hopf skew-pairing and admits a natural triangular decomposition $U_{\bqf}=U_{\bqf}^-U_{\bqf}^0U_{\bqf}^+$. 
Inside $U_{\bqf}$, we define the \emph{divided power integral form} $\Ub_{\bqf,\Acf}$ as the $\Acf$-subalgebra generated by $U_{\bqf}^0$ and the $\Acf$-duals $\Ub^{\pm}_{\bqf,\Acf}$, with respect to the Hopf skew-pairing, of the $\Acf$-subalgebras generated by $U^{\pm}_{\bqf}$. Here, $\Acf\subset\mathbb{C}(\nu)$ is a suitable $\mathbb{Z}$-subalgebra that allows us to consider Lusztig isomorphisms on the $\Acf$-form of $U_{\bqf}$. This definition is intrinsic and requires no presentation.

We show that $\Ub^{\pm}_{\bqf,\Acf}$ is generated as an algebra by quantum divided powers $E^{(n)}_{\beta}$ and $F^{(n)}_{\beta}$, where $\beta$ ranges over the positive roots of the root system of $\bqf$, and they provide PBW-type bases for these algebras. Moreover, we prove that $\Ub_{\bqf,\Acf}$ admits a triangular decomposition whose zero part is strictly larger than the group algebra of the torus: it involves \emph{quantum brackets} attached to the Cartan roots of $\bqf$, which are analogues of Lusztig's Gaussian binomials. We also show that Lusztig isomorphisms restrict to isomorphisms of $\Acf$-algebras
$\Ub_{\rho_i(\bqf),\Acf}\longrightarrow\Ub_{\bqf,\Acf}$,
using the Weyl groupoid of $\bqf$ in an essential way. In fact, these isomorphisms are instrumental in proving the triangular decomposition. 

\medbreak

At this point, we note that the specialization $U_{\bqf,\Acf}\otimes_{\Acf}\mathbb{C}$ is isomorphic to the so-called large quantum group $U_\bq$ \cite[Definition 4.2]{AAY}, whose positive part is the distinguished pre-Nichols algebra $\wtoba_\bq$; see \cite[Theorem 5.7]{AAY} and \cite{Ang-dpn}. In other words, $U_\bq$ is the quantum group \`a la De Concini--Kac--Procesi associated with $\bq$. The Cartan roots of $\bq$ enter the picture here as the roots $\beta$ for which the powers $E_{\beta}^{N_{\beta}}$ generate a central subalgebra $Z_\bq^+$, such that the Nichols algebra $\toba_{\bq}$ is the quotient of $\wtoba_{\bq}$ by $Z_\bq^+$ \cite{Ang-dpn}.

On the other hand, the graded dual of $\wtoba_\bq$ is the \emph{Lusztig algebra} $\lu_{\bq}$ of \cite{AAR-Belgium}, which fits into an extension of braided Hopf algebras by \cite[Theorem 3.3]{AAR-Belgium}:
$$
\begin{tikzcd}[column sep=large, arrow style=tikz, >=To]
	\toba_{\bq} \arrow[r, hook] & \lu_{\bq} \arrow[r, two heads] & U(\n^{+}_{\bq}).
\end{tikzcd}
$$

Here, $\n^{+}_{\bq}$ is the nilpotent part of the semisimple Lie algebra $\g_{\bq}$ constructed in \cite{AAR-MRL,AAR3}. Its root system is spanned by the Cartan roots of $\bq$, and it is explicitly determined in those references for each $\bq$. It turns out that $\lu_{\bq}\simeq\Ub^{-}_{\bqf,\Acf}\otimes_{\Acf}\mathbb{C}$. Thus, what is missing, and what we provide here, is the entire Hopf algebra that contains the small quantum group $\mathfrak{u}_{\bq}$ and admits a surjection onto the whole of $U(\g_{\bq})$.

\medbreak

Explicitly, our main construction is the \emph{quantum group of divided powers of} ${\bq}$:
\begin{align*}
\Ub_{\bq}= \Ub_{\bqf,\Acf}\otimes_{\Acf} \mathbb{C}/\Ib_{\bq},
\end{align*}
where $\Ib_{\bq}$ is the central Hopf ideal of $\Ub_{\bqf,\Acf}\otimes_{\Acf} \mathbb{C}$ generated by the elements $K_\beta^{N_\beta}-1$ for which $N_\beta\beta\in\varPi^{\bq}$, the set of simple roots of the root system of $\g_{\bq}$. We investigate the structure of $\Ub_{\bq}$ in detail and prove that it fits into an exact sequence of the desired form; see Figure~\ref{fig:intro-diagram}. In particular, when $\bq$ is of super type, this construction yields a super quantum group of divided powers \`a la Lusztig.

To establish these results, we first observe that $\Ub_{\bq}$ inherits the triangular decomposition. We then prove that it is generated by the generators of the small quantum group together with the divided powers $E^{(N_{\beta})}_{\beta}$ and $F^{(N_{\beta})}_{\beta}$ and the quantum brackets corresponding to roots $\beta$ such that $N_\beta\beta\in\varPi^{\bq}$; in particular, it is finitely generated. We also establish a PBW basis in terms of these generators. These results allow us to prove that $\uf_{\bq}$ is normal in $\Ub_\bq$. To identify the Lie algebra appearing in the quotient, we use the multiplication formulas satisfied by these generators.

Furthermore, we show that $\Ub_{\bq}$ is pointed and, under mild assumptions, describe its group of group-like elements by generators and relations. Finally, we prove that the Lusztig isomorphisms descend to isomorphisms $\Ub_{\rho_i(\bq)}\to \Ub_{\bq}$.

\medbreak

The question of whether one can construct quantum groups of divided powers for every finite-dimensional Nichols algebra of diagonal type has attracted the attention of several specialists. In this paper we give a partial answer. We were encouraged to work on this problem by Nicol\'as Andruskiewitsch, with whom we have discussed it extensively; some of these discussions are reflected in \cite{AAV}. Indeed, the Lusztig algebras $\lu_{\bq}$ were first introduced in passing in our joint paper \cite{AAGTV} with him, Garc\'ia Iglesias and Torrecillas, where it was suggested that they might play an important role. Later, Andruskiewitsch, the first author and Rossi Bertone carried out a thorough study of them in \cite{AAR-MRL,AAR-Belgium,AAR3}. The results obtained there play a fundamental role in the construction presented here. 

It is natural to expect that the super quantum groups introduced here also have remarkable properties in its representation theory, analogous to those of Lusztig's quantum groups. We are exploring this in \cite{AHV}, together with T. Heidersdorf.

Finally, we should mention other works constructing super analogues of Lusztig's quantum groups. In \cite{ElK}, El Turkey and Kujawa construct a $\Z[\nu^{\pm1}]$-form of divided powers for the quantum group $U_\nu(\mathfrak{gl}(m|n))$, as presented by Zhang \cite{zhang}, and provide a basis for it by exploiting a $\nu$-super Schur--Weyl duality associated with $\mathfrak{gl}(m|n)$. In a series of papers, Clark, Hill, Wang and collaborators generalized important results of Lusztig, see e.g. \cite{CHW}; however, the only new example of super type is the one associated with $\mathfrak{osp}(1|2n)$. 
A recent work \cite{CLZ} investigates the algebraic structure of the quantum queer superalgebra and its Lusztig integral form, determining explicitly quantum root vectors and a PBW-type basis; they also analyze the specialization at primitive roots of unity.
In \cite{AAR-MRL}, the Drinfeld double of $\lu_{\bq}\#\Z^n$ is called a quantum divided power algebra; however, unlike the construction presented here, it lacks a Cartan subalgebra.

\medbreak

The manuscript is organized as follows. Section \ref{sec:background} collects the background material. We recall the Drinfeld double of the bosonization
of a Nichols algebra of diagonal type and its triangular decomposition, the associated root system, Weyl groupoid and Lusztig isomorphisms, and the resulting PBW bases. We then introduce the ring of coefficients $\Acf$, and show that the PBW bases descend to $\Acf$-bases of $U_{\bqf,\Acf}$. We close the section with the Cartan roots of $\bq$, the semisimple Lie algebra $\g_{\bq}$ they determine, and the specialization at a root of unity, where the distinguished pre-Nichols algebra, the Lusztig algebra and the large quantum group appear.

Section \ref{sec:int-forms-divided-powers} is devoted to the divided power integral form $\Ub_{\bqf,\Acf}$, defined by duality with respect to the Hopf skew-pairing. We prove that its positive and negative parts are generated by the quantum divided powers attached to the positive roots, which moreover give PBW-type bases, and that $\Ub_{\bqf,\Acf}$ is a Hopf $\Acf$-subalgebra of $U_{\bqf}$. Its zero part is strictly larger than the group algebra of the torus, as it involves quantum brackets, which we study in Subsection \ref{subsec:quantum-brackets}. The main result of the section is the triangular decomposition of $\Ub_{\bqf,\Acf}$, proved in Subsection \ref{subsec:div-powers-triang-decomp} together with the fact that the Lusztig isomorphisms restrict to it.

In Section \ref{sec:specialization} we specialize at a root of unity and pass to the quotient by a central Hopf ideal, obtaining the quantum group of divided powers $\Ub_{\bq}$. Subsection \ref{subsec:specialization-first-prop} shows that it inherits the triangular decomposition and that it is generated by finitely many elements: those of the small quantum group together with the divided powers and the quantum brackets attached to the indecomposable Cartan roots. Subsection \ref{subsec:extension} contains the main result of the paper: the small quantum group is a normal Hopf subalgebra of $\Ub_{\bq}$ with quotient $U(\g_{\bq})$. In Subsection \ref{subsec:specialization-Lusztig-iso} we prove that the Lusztig isomorphisms descend to $\Ub_{\bq}$ compatibly with this extension, and describe the isomorphisms they induce between the enveloping algebras. We close in Subsection \ref{subsec:example-B11} with a fully explicit description of the quantum group of divided powers of super type $\superb{1}{1}$.

\subsection*{Acknowledgements}
We are deeply grateful to Nicolás Andruskiewitsch, our former PhD advisor and long-time mentor, for many years of collaboration and fruitful discussions on this and many other problems. We also thank M. Yakimov and T. Heidersdorf for stimulating discussions that inspired parts of this work. AI tools were used for proofreading and to improve the exposition; the mathematical content is entirely the authors' own. This work was partially done during the visit of I. A. to Guandong Technion Israel Institute of Technology during the first semester of 2026.

\section{Background}\label{sec:background}

In this section, we fix the notation regarding Nichols algebras of diagonal type and their Drinfeld doubles, and recall their main features that are instrumental for our construction. We refer to \cite{AA-diag-survey,AAY} for more details, and to \cite{Radford-book} for unexplained notation on Hopf algebras.

Our initial data are the matrices $\bq$ and $\bqf$ determining the Nichols algebras. Throughout this work $\theta\in\N$, $\I=\{1, \dots, \theta\}$ and $\nu$ is an indeterminate. Let
\begin{align}\label{eq:braiding-matrix}
\bq = (\bq_{ij} ) \in \big( \Cset^\times \big)^{\I \times \I}
\end{align}
be an indecomposable braiding matrix with finite-dimensional Nichols algebra $\toba_\bq$; these were classified in \cite{H-classif RS}. We further assume that $\bq$ satisfies \eqref{eq:condition-cartan-roots}, which will be explained later. Moreover, we assume that $\bq$ belongs to a one-parameter family \cite{AAY}. Among other things, this means that there exists a root of unity $\xi\in\Cset$ and
\begin{align}\label{eq:braiding-matrix-parameter}
\bqf = ( \qf_{ij}) \in   \big( \Cset[\nu^{\pm 1}]^\times \big)^{\I \times \I}
\end{align}
such that $\bq=\bqf(\xi)$ and the Nichols algebra $\toba_\bqf$ associated to $\bqf$ has the same arithmetic root system as $\bq$. For short we will say that {\it $\bqf$ specializes to $\bq$ at $\xi$}.

In what follows, the symbol $\dag$ will be used to refer to both $\bqf$ and $\bq$. We set:
\begin{itemize}\renewcommand{\labelitemi}{$\circ$}
\item $\Cset_\bq=\Cset$ and $\Cset_\bqf=\Cset(\nu)$; notice that $\dag\in\Cset_{\dag}^{\I \times \I}$. 
\item $V_\dag$ is the braided $\Cset_\dag$-vector space with basis $\{E_1, \dots, E_\theta\}$ and braiding 
\begin{align*}
c_{\dag}(E_i\otimes E_j)&=\dag_{ij}E_j\otimes E_i, & & i,j\in\I. 
\end{align*}
\item The Nichols algebra $\toba_\dag$ is the braided Hopf algebra quotient $T(V_\dag)/\JJ(V_\dag)$. A minimal set of generators for the ideal $\JJ(V_\bq)$ was given in \cite{A-presentation}, and from there one can deduce a corresponding set for $\JJ(V_{\bqf})$ as explained in \cite[Section 4.1]{AAY}. However, we do not need here a precise description of them. 
\item Similarly, $W_{\dag^{(-1)}}$ is the braided $\Cset_\dag$-vector space with basis $\{F_1, \dots, F_\theta\}$ and braiding matrix $\dag^{(-1)}=(\dag_{ij}^{-1})\in(\Cset_\dag^{\times})^{\I\times\I}$, i.e. $c_{\dag^{(-1)}}(F_i\otimes F_j)=\dag_{ij}^{-1}F_j\otimes F_i$, and $\toba_{\dag^{(-1)}}$ is the corresponding Nichols algebra.
\end{itemize}


\begin{remark}\label{rmk:N tilde}
We will identify the matrix $\dag$ with the  $\Z$-bilinear form $\Z^\I\times\Z^\I\longrightarrow\Cset_\dag^\times$ that it defines. Namely, let $\{\alpha_1, \dots,\alpha_\theta\}$ denote the canonical basis of $\Z^\I$. Then $\dag(\alpha_i,\alpha_j)=\dag_{ij}$, $i,j\in\I$. For $\alpha,\beta\in\Z^\I$, we define
\begin{align*}
\dag_{\alpha\beta}&=\dag(\alpha,\beta), &
\dag_\beta&=\dag_{\beta\beta} & &\mbox{and} & N^\dag_\beta&=\ord\dag_{\beta}\in\N\cup\{\infty\}.
\end{align*}
We also denote $\widetilde N_\beta =  N_{\beta}^\bqf$ and $ N_\beta =  N_{\beta}^\bq$, as this is the usual notation for these numbers.
\end{remark}

\subsection{Quantum groups \`a la Drinfeld}

In the next section, we will construct quantum groups of divided powers starting from a Drinfeld double as that in \cite[\S 4.2]{AAY}, see also \cite[\S 4]{H-Lusztig}, which we now recall. Let $D_\dag$ be the $\Cset_\dag$-algebra generated by the elements $E_i$, $K_i^{\pm1}$, $F_i$ and $L_i^{\pm1}$, $i\in\I$, subject to the relations
\begin{align}
K_iK_i^{-1}&=L_iL_i^{-1}=1, & K_iL_j &=L_j K_i, 
\\
K_iK_j&=K_jK_i, & L_iL_j &=L_j L_i, 
\\
\label{eqn:g-con-Ei} 
K_iE_j &= \dag_{ij}E_jK_i,& L_iE_j &= \dag_{ji} E_jL_i,
\\\label{eqn:g-con-Fi}
K_iF_j &= \dag_{ij}^{-1}F_jK_i,& L_iF_j &= \dag_{ji}^{-1} F_jL_i,
\\\label{eq:Ubqf-linking-EF}
E_iF_j - F_j E_i &= \delta_{ij} (K_i -L_i^{-1}), & i,j&\in\I,
\end{align}
In addition,  $E_1, \dots, E_\theta$ are subject to the relations of the Nichols algebra $\toba_\dag$, and $F_1, \dots, F_\theta$ to those of $\toba_{\dag^{(-1)}}$. 
Moreover, $D_\dag$ is a Hopf algebra with comultiplication determined by
\begin{align*}
\Delta(K_i) &=K_i\otimes K_i,&\quad\Delta(L_i)&=L_i\otimes L_i,\\
\Delta(E_i) &= K_i \otimes E_i + E_i \otimes 1,&\quad
\Delta(F_i) &=1 \otimes F_i + F_i \otimes L_i^{-1}.
\end{align*}

We set the following $\Cset_\dag$-subalgebras of $D_{\dag}$:
\begin{align*}
D_{\dag}^{+} &= \Cset_\dag\langle E_i: i\in \I \rangle, &
D_{\dag}^{+0} &= \Cset_\dag\langle K_i^{\pm 1}: i\in \I\rangle, &
D_{\dag}^{\geqslant} &= D_{\dag}^{+}D_{\dag}^{+0}, &
D_{\dag}^{0} &= D_{\dag}^{+0}D_{\dag}^{-0},
\\
D_{\dag}^{-} &= \Cset_\dag\langle F_i: i\in \I \rangle, &
D_{\dag}^{-0} &= \Cset_\dag\langle L_i^{\pm 1}: i\in \I\rangle, &
D_{\dag}^{\leqslant}&= D_{\dag}^{-}D_{\dag}^{-0}.
\end{align*}
The multiplication induces an isomorphism of vector spaces
\begin{align}\label{eq:triangular decomposition of Ubqf}
D_{\dag}^-\otimes D_{\dag}^0\otimes D_{\dag}^+\overset{\sim}{\longrightarrow} D_{\dag}.
\end{align}

We observe that $D_{\dag}^{\pm0}$ is a group algebra of a free abelian group of rank $\theta$, and $D_{\dag}^+$ (respectively, $D_{\dag}^-$) has a canonical structure of Hopf algebra in the category ${}_{D_{\dag}^{+0}}^{D_{\dag}^{+0}}\mathcal{YD}$ (resp. $\mathcal{YD}_{D_{\dag}^{-0}}^{D_{\dag}^{-0}}$) of left (resp. right) Yetter--Drinfeld modules over $D_{\dag}^{+0}$ (resp. $D_{\dag}^{-0}$). Moreover, we have isomorphisms
\begin{align*}
D_{\dag}^+ &\simeq\toba_\dag, & D_{\dag}^{\geqslant} &\simeq D_{\dag}^+\# D_{\dag}^{+0}, & 
D_{\dag}^- &\simeq\toba_{\dag^{(-1)}}, & D_{\dag}^{\leqslant}&\simeq D_{\dag}^{-0}\#D_{\dag}^-
\end{align*}
of (braided) Hopf algebras, as appropriate.

In other words, $D_{\dag}$ is the Drinfeld double of $D_{\dag}^{\geqslant}$. Indeed, there is a unique Hopf skew-pairing 
$\langle \cdot , \cdot \rangle: D_{\dag}^{\leqslant} \times D_{\dag}^{\geqslant} \to \Cset_{\dag}$ determined by
\begin{align}\label{eq:skew pairing Ubqf}
\langle L_i, K_j \rangle &= \dag_{ji}^{-1}, &
\langle F_i, E_j \rangle &= \delta_{ij}, &
\langle L_i, E_j \rangle &= \langle F_i, K_j \rangle = 0, & i,j &\in \I, 
\end{align}
and the associated Drinfeld double is isomorphic to $D_{\dag}$ \cite{H-Lusztig}. Moreover, the restriction $\langle \cdot , \cdot \rangle: D_{\dag}^{-} \times D_{\dag}^{+} \to \Cset_\dag$
is non-degenerate and is a 
Hopf skew-pairing of braided Hopf algebras. We notice that, for $x_{\pm}\in D_{\dag}^{\pm}$ and $g_{\pm} \in D_{\dag}^{\pm0}$, the following identity holds
\begin{align*}
\langle g_-x_-, x_+g_+ \rangle = \langle x_-, x_+ \rangle \langle g_-, g_+ \rangle.
\end{align*}

For each $\beta=n_1\alpha_1+\cdots+n_\theta\alpha_\theta\in\Z^\I$, we set $K_\beta=K_1^{n_1}\cdots K_\theta^{n_\theta}$ and $L_\beta=L_1^{n_1}\cdots L_\theta^{n_\theta}$. We notice that $D_{\dag}$ is a $\Z$-graded Hopf algebra 
by setting
\begin{align*}
\deg E_i&=\alpha_i=-\deg F_i, & \deg K_i&=\deg L_i=0, & \text{for all } &i\in\I. 
\end{align*}

\begin{remark}
Later, we will write $U_\bqf=D_\bqf$ following the notation in \cite[\S 4.2]{AAY}.
\end{remark}

\subsection{Root systems, PBW basis and Lusztig isomorphisms}\label{subsec:weyl-gpd}

By \cite{karchencko}, there exist $E_{\beta_1}, \dots, E_{\beta_\ell}\in\toba_\dag$ with $\deg E_{\beta_i}=\beta_i$ such that
\begin{align*}
\left\{E_{\beta_1}^{n_1}\cdots E_{\beta_\ell}^{n_\ell}\mid0\leq n_j<N^\dag_{\beta_j},\,1\leq j\leq\ell\right\}
\end{align*}
is a basis of $\toba_\dag$. By our assumption, the set of degrees is the same for $\bqf$ and $\bq$, and is called the set of positive roots. We set
\begin{align*}
 \varDelta^\dag_+=\{\beta_1, \dots, \beta_\ell\}\quad\mbox{and}\quad \varDelta^{\dag} = \varDelta_+^{\dag} \cup -\varDelta_+^{\dag}.
\end{align*}

The (generalized) Cartan matrix $(c_{ij}^{\dag})_{i,j\in \I}\in\Z^{\I\times\I}$ associated to $\dag$ is defined in the diagonal by $c_{ii}^{\dag} = 2$, and for $i\neq j$ by
\begin{align}\label{eq:defcij}
c_{ij}^{\dag}&:= -\min \left\{ n \in \N_0: (n+1)_{\dag_{ii}}
(1-\dag_{ii}^n \dag_{ij}\dag_{ji} )=0 \right\}.
\end{align}
This gives rise, for each $i\in\I$, to the reflection $s_i^{\dag}\in GL(\Z^\I)$ determined by 
\begin{align}\label{eq:siq-definition}
s_i^{\dag}(\alpha_j)&:=\alpha_j-c_{ij}^{\dag}\alpha_i,
\end{align}
for $j\in\I$, which, in turn, defines the matrix $\rho_i(\dag)$ by
\begin{align}\label{eq:rhoiq-definition}
(\rho_i (\dag))_{jk}&= \dag(s_i^{\dag}(\alpha_j),s_i^{\dag}(\alpha_k)) = \dag_{jk} \dag_{ik}^{-c_{ij}^{\dag}} \dag_{ji}^{-c_{ik}^{\dag}} \dag_{ii}^{c_{ij}^{\dag}c_{ik}^{\dag}},
\end{align}
for $j,k\in\I$. Set $\cX^\dag = \{\rho_{j_1} \dots \rho_{j_n}(\dag)\mid j_1, \dots, j_n \in \I, n \in \N \}$. 
It turns out that the Nichols algebra associated to $\rho_{j_1} \dots \rho_{j_n}(\bq)$ is also in the classification \cite{H-classif RS}, and $\rho_{j_1} \dots \rho_{j_n}(\bqf)$ specializes to $\rho_{j_1} \dots \rho_{j_n}(\bq)$. Moreover, $\rho_{j_1} \dots \rho_{j_n}(\varDelta^\dag)=\varDelta^{\rho_{j_1} \dots \rho_{j_n}(\dag)}$.

The family of reflections $s_{j_1}^{\rho_{j_1} \dots \rho_{j_n}(\dag)}:\rho_{j_2} \dots \rho_{j_n}(\dag)\rightarrow\rho_{j_1} \dots \rho_{j_n}(\dag)$ leads to the so-called Weyl groupoid $\cW_\dag$. It is known that there exists an element $\omega_0^{\dag}=s_{i_1}^{\dag} \cdots s_{i_\ell}:\rho_{i_\ell}\cdots \rho_{i_1}(\dag)\rightarrow\dag$ of maximal length in $\cW_\dag$, and every reduced composition of reflections can be extended to a reduced expression of $\omega_0^{\dag}$. Moreover, the positive roots $\varDelta_+^\dag$ can be constructed as
\begin{align} \label{eq:betak}
\beta_k &= s_{i_1}^{\dag}\cdots s_{i_{k-1}}(\alpha_{i_k}) \qquad \text{for all }k\in\{1, \dots, \ell\}.
\end{align}

The Lusztig isomorphisms $T_i^{\dag}:U_{\rho_i (\dag)} \to U_{\dag}$, $i\in\I$, were introduced in \cite[Proposition 6.8]{H-Lusztig}. These are algebra isomorphisms defined by
\begin{align}\label{eq:lusztig-isom-defn}
\begin{aligned}
T_i^{\dag}(K_j) &=  K_j K_i^{-c_{ij}^{\dag}};
&
T_i^{\dag}(E_i) &= \begin{cases} 
F_i L_i, & j=i, \\ (\ad_c E_i)^{-c_{ij}^{\dag}} E_j, & j\neq i,
\end{cases}
\\
T_i^{\dag}(L_j) &= L_j L_i^{-c_{ij}^{\dag}};
&
T_i^{\dag}(F_i) &= \begin{cases} 
K_i^{-1} E_i, & j=i, \\ 
(\lambda_{ij}^{\dag})^{-1} (\ad_c F_i)^{-c_{ij}^{\dag}} F_j, & j\neq i,
\end{cases}
\end{aligned}
\end{align}
by abuse of notation, we use the same symbols for the generators of $U_{\rho_i(\dag)}$ and of $U_{\dag}$; here
\begin{align}\label{eq:lambdaij-defn}
\lambda_{ij}^{\dag} &=(\dag_{ii}^{-c_{ij}^{\dag}} \dag_{ij}\dag_{ji})^{c_{ij}^{\dag}} (-c_{ij}^{\dag})_{\dag_{ii}}^!
\prod_{0\le s<-c_{ij}^{\dag}} (\dag_{ii}^s \dag_{ij}\dag_{ji}-1) \in \Cset_\dag^\times.
\end{align}
Notice that $\lambda_{ij}^{\dag} \neq 0$ by the very definition of $c_{ij}^\dag$ \eqref{eq:defcij}.

A set of PBW generators can be constructed using the Lusztig isomorphisms, cf. \cite[Theorem 4.5]{HY-shapov}. Namely, we fix a reduced expression $\omega_0^{\dag}=s_{i_1}^{\dag} \cdots s_{i_\ell}$ of the longest element in the Weyl groupoid and define
\begin{align}\label{eq:Ebeta-defn}
E_{\beta_k}=T_{i_1}^{\dag} \dots T_{i_{k-1}} (E_{i_k}) \in U_{\dag}^+\quad\mbox{and}\quad
F_{\beta_k}=T_{i_1}^{\dag} \dots T_{i_{k-1}} (F_{i_k}) \in U_{\dag}^-
\end{align}
for $k\in\{1, \dots, \ell\}$, where $E_{i_k}\in U^+_{\rho_{i_{k-1}\cdots \rho_{i_1}(\dag)}}$ and $F_{i_k}\in U^-_{\rho_{i_{k-1}\cdots \rho_{i_1}(\dag)}}$. Then
\begin{align}\label{eq:PBW-basis-U+-U-}
\begin{split}
\left\{  E_{\beta_1}^{n_1}E_{\beta_2}^{n_2} \dots E_{\beta_\ell}^{n_\ell} \mid 0\le n_j< N^\dag_{\beta_j}, 1\leq j \leq\ell \right\}\quad\mbox{and}\\
\left\{  F_{\beta_1}^{m_1}F_{\beta_2}^{m_2} \dots F_{\beta_\ell}^{m_\ell} \mid 0\le m_j< N^\dag_{\beta_j}, 1\leq j \leq\ell \right\}
\end{split}
\end{align}
are bases of $U_{\dag}^+$ and $U_{\dag}^-$, respectively. Moreover,
\begin{align}\label{eq:PBW-basis-U}
\left\{ E_{\beta_1}^{m_1} \dots E_{\beta_\ell}^{m_\ell} K_1^{a_1} \dots K_{\theta}^{a_{\theta}}  L_1^{b_1} \dots L_{\theta}^{b_{\theta}} F_{\beta_1}^{n_1} \dots F_{\beta_\ell}^{n_\ell}:
0\le m_j, \, n_j< N_{\beta_j}^\dag, \, a_i, \, b_i\in\Z \right\}
\end{align}
is a basis of $D_{\dag}$.

\subsection{Restricting scalars}\label{subsection:restricting}

We denote by $\Acf$ the $\mathbb{Z}$-subalgebra of $\Cset_\bqf=\Cset(\nu)$ defined according to the type of $\bqf$ as shown in Table \ref{table:Acf}. This algebra allows us to restrict the scalars while preserving all the Lusztig isomorphisms as we see below.

\begin{table}[ht]
\centering
\small
\renewcommand{\arraystretch}{1.6}
\arrayrulecolor{accent}
\setlength{\tabcolsep}{3pt}
\begin{tabular}{ll}
\toprule
\rowcolor{accent!8}
$\bqf$ & $\Acf$ \\
\midrule
Cartan type
& $\Z\left[\nu^{\pm1},(\nu^s-1)^{-1}\mid1\leq s<\max_{i\neq j}|c_{ij}^\bqf|\right]$ \\
Super $\Ab(k\mid\theta-k)$, $1\leq k\le \lfloor\frac{\theta+1}{2}\rfloor$
& $\Z\left[\nu^{\pm1},(\nu-1)^{-1}\right]$ \\
Super $\Bb(k\mid\theta-k)$, $1\leq k\le\theta-1$
& $\Z\left[\nu^{\pm1},(\nu-1)^{-1},(\nu^2-1)^{-1}\right]$ \\
Super $\Db(k\mid\theta-k)$, $1\leq k\le\frac{\theta}{2}$
& $\Z\left[\nu^{\pm1},(\nu-1)^{-1},(\nu^2-1)^{-1}\right]$ \\
Super $\Db(2,1;\alpha)$, $d_1,d_3\in\N$
& $\Z\left[\nu^{\pm d_i},(\nu^{d_i}-1)^{-1},(\nu^{d_1+d_3}-1)^{-1}\mid i=1,3\right]$ \\
Super $\Fb(4)$
& $\Z\left[\nu^{\pm1},(\nu^{s}-1)^{-1}\mid 1\leq s<4\right]$ \\
Super $\Gb(3)$
& $\Z\left[\nu^{\pm1},(\nu^{s}-1)^{-1},(\nu+1)^{-1}\mid 1\leq s<4\right]$ \\
$\texttt{wk}(4)$
& $\Z\left[\frac{1}{2},\nu^{\pm1},(\nu\pm1)^{-1}\right]$ \\
$\texttt{br}(2)$, $\zeta$ a $3$rd root of unity
& $\Z\left[\zeta,\nu^{\pm1},(\nu-1)^{-1},(\nu-\zeta)^{-1}\right]$ \\
\bottomrule
\end{tabular}
\caption{The $\mathbb{Z}$-subalgebra $\Acf$ is generated by $\bpf_{ii}^{\pm1}$
and $(\bpf_{ii}^s \bpf_{ij}\bpf_{ji}-1)^{-1}$ for all $i\neq j \in \I$,
$0\le s<-c_{ij}^{\bpf}$, and $\bpf$ Weyl equivalent to $\bqf$.}
\label{table:Acf}
\end{table}

\begin{lemma}\label{le:elements in Acf}
Let $n \in \N$, $j_1, \dots, j_n \in \I$, $\bpf = \rho_{j_1} \dots \rho_{j_n}(\bqf)$. Then
the elements 
\begin{align*}
\bpf_{\alpha\beta},\quad
\binom{m}{r}_{\bpf_\beta},\quad
(t)_{\bpf_{ii}}^{-1}\quad\mbox{and}\quad
(\lambda_{ij}^{\bpf})^{\pm1}
\end{align*}
belong to $\Acf$ for all $\alpha,\beta\in\Z^\I$, $i,j\in\I$, $0\leq r\leq m$ and $2\leq t\leq -c_{ij}^{\bpf}$. 
\end{lemma}

\begin{proof}
Our choice of the generators of $\Acf$ is restricted to the validity of this statement, after analyzing each one-parameter family. 
More precisely, as observed in the proof of \cite[Proposition A.3]{AAY}, we have that $\{\bpf_{ij}\mid i,j\in\I\} \subset \Z[\nu^{\pm d_1},\nu^{\pm d_3}]$ if $\bpf$ is of type $\Db(2,1;\alpha)$, $\{\bpf_{ij}\mid i,j\in\I\} \subset \Z[\zeta,\nu^{\pm1}]$ if $\bpf$ is of type $\texttt{br}(2)$, and $\{\bpf_{ij}\mid i,j\in\I\}\subset\Z[\nu^{\pm1}]$ in the remaining types. Therefore $\{\bpf_{\alpha\beta}\mid \alpha,\beta\in\Z^\I\}\subset\Acf$, which also implies that the $q$-binomial numbers belong to $\Acf$ because 
of the $q$-Pascal formula.

In turn, if $\bpf_{ii}$ is a root of unity, then $\bpf_{ii}\in\{-1,\zeta\}$ and hence $(t)_{\bpf_{ii}}^{-1}\in\Acf$. Otherwise, $\bpf_{ii}$ is a power of $\nu$ and $(\nu^t-1)^{-1}\in\Acf$ by definition. Then $(t)_{\bpf_{ii}}^{-1}\in\Acf$ as 
$\displaystyle(t)_{\bpf_{ii}}=\frac{1-\bpf_{ii}^t}{1-\bpf_{ii}}$. Finally, $(\bpf_{ii}^s \bpf_{ij}\bpf_{ji}-1)^{-1}$ can be constructed using the generators of $\Acf$ for all $0\le s<-c_{ij}^{\bpf}$, and therefore $(\lambda_{ij}^{\bpf})^{\pm1}$.
\end{proof}

Let us denote by $U_{\bqf,\Acf}^\star$ the $\Acf$-subalgebra of $U_\bqf=D_\bqf$ generated by the same elements as $U_\bqf^\star$, where $\star\in\{\pm,0,\pm0,\geqslant,\leqslant\}$. 
The triangular decomposition of $U_{\bqf}$ induces the corresponding one of $U_{\bqf,\Acf}$, i.e. the multiplication gives an $\Acf$-linear isomorphism
\begin{align}
\label{eq:triang-int}
U_{\bqf, \Acf}^+ \otimes_{\Acf} U_{\bqf, \Acf}^0 \otimes_{\Acf} U_{\bqf, \Acf}^- \overset{\sim}{\longrightarrow} U_{\bqf,\Acf}.
\end{align}

We can mimic the proofs in \cite[\S 5.2]{AAY} to get the following.

\begin{proposition}\label{prop:PBW-int}
\begin{enumerate}[leftmargin=*,label=\rm{(\alph*)}]
\item\label{item:Ti-lemma-restriction} Let $\bpf=\rho_{j_1} \dots \rho_{j_n}(\bqf)$ for some $j_1, \dots, j_n \in \I$ and $n \in \N$. For all $i\in\I$, $T_i^{\bpf}$ restricts to an $\Acf$-algebra isomorphism
\begin{align*}
T_{i}^{\bpf}:U_{\rho_i (\bpf), \Acf}\to U_{\bpf, \Acf}.
\end{align*}
\item\label{item:Ti-lemma-rootvectors} The PBW generators $E_{\beta_k}$ and $F_{\beta_k}$ of in \eqref{eq:Ebeta-defn} belong to $U_{\bqf, \Acf}$ for all $1\leq k\leq\ell$.
\item\label{item:PBW-int-form} The sets \eqref{eq:PBW-basis-U+-U-} and \eqref{eq:PBW-basis-U} are $\Acf$-bases of $U_{\bqf, \Acf}^\pm$ and $U_{\bqf, \Acf}$, respectively.
\end{enumerate}

\end{proposition}

\pf
The scalars $\lambda_{ij}^{\bpf}$ from \eqref{eq:lambdaij-defn} involved in the definition of $T_{i}^{\bpf}$ in \eqref{eq:lusztig-isom-defn} are invertible in $\Acf$ thanks to Lemma \ref{le:elements in Acf}. This implies \ref{item:Ti-lemma-restriction}, while \ref{item:Ti-lemma-rootvectors} follows from
\ref{item:Ti-lemma-restriction} and \eqref{eq:Ebeta-defn}.
For \ref{item:PBW-int-form}, we consider the case of $U_{\bqf, \Acf}^+$, the proof for the other cases are analogous. Let $Y$ be the set of
PBW monomials of  $U_{\bqf}^+$ from \eqref{eq:PBW-basis-U+-U-}.
By \ref{item:Ti-lemma-rootvectors}, $Y \subset U_{\bqf, \Acf}^+$.
The defining relations of $U_{\bqf}^+$ involve products of $E_i$ with coefficients in $\Acf$, hence we may prove recursively that, for $j>k$, $E_{\beta_j}E_{\beta_k} \in \Acf Y$, the $\Acf$-module generated by $Y$, where each monomial in the expansion has letters $E_{\beta_t}$, $j>t>k$;
see the proof of \cite[Theorem 4.8]{HY-shapov}. Thus $\Acf Y$ is a left ideal containing $1$, so $\Acf Y =  U_{\bqf, \Acf}^+$.
Since the elements in $Y$ are linearly independent over $\Cset(\nu)$, they are so over $\Acf$ and hence \ref{item:PBW-int-form} follows.
\epf

\begin{remark}
Compared with \cite{L-JAMS,L-book}, our ring $\Acf$ is bigger. We did so in order to ensure the invertibility of the algebra maps $T_{i}^{\bpf}$. On the other hand, we chose $\Acf$ subject to this condition, which will allow us to evaluate the parameter in any root of unity, with very few exceptions (only order $2$ or $3$ in most cases).
\end{remark}

\subsection{Cartan roots}\label{subsec:cartan-roots}
An index $i\in\I$ is called a Cartan vertex of $\dag$, if $\dag_{ij}\dag_{ji} = \dag_{ii}^{c_{ij}^{\dag}}$ for all $j\neq i$. Then, the set of Cartan roots of $\dag$ is
\begin{align*}
\fO^{\dag} &= \{s_{i_1}^{\dag} s_{i_2} \dots s_{i_k}(\alpha_i) \in \varDelta^{\dag}\mid
i\in \I  \text{ is a Cartan vertex of } \rho_{i_k} \dots \rho_{i_2}\rho_{i_1}(\dag) \}.
\end{align*}
Equivalently, $\beta\in\fO^\dag$ is a Cartan root if and only if
\begin{align}\label{eq:cartan-root}
\dag_{\beta\gamma}\dag_{\gamma\beta}&=\dag_{\beta\beta}^{a_{\beta\gamma}} \qquad \text{for all }\gamma\in\Z^\I,
\end{align}
by \cite[Remark 2.3]{AAR3}. We can check by inspection that $\fO^{\bq}=\fO^{\bqf}$. These roots are fundamental to our construction.

\begin{remark}
Fix $\beta\in\fO^\dag$, $\gamma\in\Z^{\I}$. The integer $a_{\beta\gamma}$ is unique modulo $\ord \dag_{\beta\beta}$. In case that $\ord \dag_{\beta\beta}<\infty$, we will choose $a_{\beta\beta}=2$, and for $\gamma\notin \Z\beta$ as the unique value such that $-a_{\beta\gamma}\in [0, \ord \dag_{\beta\beta})$.
\end{remark}

For $\beta\in\fO^\bq$, we introduce the element $\wbeta = N^{\bq}_{\beta} \beta$. By \cite[Theorem 3.7]{AAR3}, we know that the set $\wfO^{\bq}=\{\wbeta\mid \beta \in \fO^{\bq} \}$ is either empty or
a root system inside the real vector space generated by $\fO^{\bq}$, and the set $\varPi^{\bq}$ of all indecomposable elements\footnote{Those which cannot be represented as a non-trivial positive linear combination of elements of  $\wfO^{\bq}_+$} of $\wfO^{\bq}_+=\wfO^{\bq} \cap \N_{0}^{\theta}$ forms a basis of this root system. We set then accordingly $\varPi^{\bqf}\coloneqq \varPi^{\bq}$. The root system is precisely determined for each $\bq$ in \cite{AAR3}. We denote  by $\g_{\bq}$ the semisimple Lie algebra with this root system, and  fix a triangular decomposition $\g_{\bq}=\n_\bq^+ \oplus \h_{\bq} \oplus \n_\bq^-$.

\begin{remark}\label{rem:identity-cartan-entries}
Let $(a_{\wbeta\wgamma})_{\wbeta,\wgamma\in\varPi^\bq}$ be the Cartan matrix of $\g_{\bq}$. Then $N_{\beta}^\bq a_{\wbeta\wgamma}=N_\gamma^\bq a_{\beta\gamma}$ for all $\wbeta,\wgamma\in\varPi^\bq$.

This identity is immediate for $\beta=\alpha_i$, $i$ a Cartan vertex, and $\gamma=n\alpha_i+m\alpha_j$. The general case reduces to this one by the definition of Cartan root and \cite[Theorem 2.2]{AAR3}.
\end{remark}

\medbreak

We can now state our extra assumption on $\bq$, the same as in \cite{AAY}. Namely, throughout this work we will assume that
\begin{align}\label{eq:condition-cartan-roots}
q_{\alpha\beta}^{N_\beta}=1 \qquad \text{for all }\alpha \in \varDelta^{\bq}, \, \beta\in\fO^{\bq}.
\end{align}
This is equivalent to $q_{\beta\alpha}^{N_\beta}=1$ for all $\alpha \in \varDelta^{\bq}$ and $\beta\in\fO^{\bq}$; and to $q_{\alpha_i\beta}^{N_\beta}=1$ for all $i\in\I$ and $\wbeta\in\varPi^{\bq}$. Also, \eqref{eq:condition-cartan-roots} holds for $\bq$ if and only if it holds for $\rho_i (\bq)$, see \cite[Remark 4.4]{AAY}.

\subsection{Specializing at \texorpdfstring{$\bq$}{\bq}}
\label{subsec:dpn}

The distinguished pre-Nichols algebra $\wtoba_{\bq}=T(V_\bq)/\II(V_\bq)$ was introduced in \cite{Ang-dpn} where $\II(V_\bq)$ is a certain ideal contained in $\JJ(V_\bq)$. This algebra is also an $\N_0^\I$-graded braided Hopf algebra quotient of $T(V_\bq)$. The main difference, but not the only one, between the two ideals is that $E_\beta^{N_\beta^\bq}\in\JJ(V_\bq)\setminus\II(V_\bq)$ for all $\beta\in\fO^\bq$. 

Coming back to the previous subsections, recall there is a root of unity $\xi\in\Cset$ such that $\bqf(\xi) = \bq$. Let $\Bbbk$ be either $\Cset$ or $\Qset[\xi]$. We can consider $\Bbbk$ as an $\Acf$-module via the assignment $\nu \mapsto \xi$, cf. \cite[\S 5.3]{AAY}. Then the specialization $U_{\bqf,\Acf}\otimes_{\Acf}\Bbbk$ is isomorphic to the so-called large quantum group $U_\bq$ \cite[Definition 4.2]{AAY}, which is the Drinfeld double of the bosonization of the distinguished pre-Nichols algebra $\wtoba_\bq$, cf. \cite[Theorem 5.7]{AAY} and \cite{Ang-dpn}. Thus, $U_\bq$ admits a presentation by generators and relations analogous to that of $D_\bq$ but $\JJ(V_\bq)$ is replaced by $\II(V_\bq)$, and $\JJ(W_{\bq^{(-1)}})$ by $\II(W_{\bq^{(-1)}})$. Moreover, the PBW basis of $\toba_{\bqf^{(\pm1)}}$ given in \eqref{eq:PBW-basis-U+-U-} descends to a $\Bbbk$-basis of $\wtoba_{\bq^{(\pm1)}}$. Also, after specialization, the subalgebras $U_\bqf^\star$, for $\star\in\{\pm,0,\pm0,\geqslant,\leqslant\}$, and the Lusztig isomorphisms $T_i^\bqf$ descend to comparable algebraic data which we denote by $U_\bq^\star$ and $T^\bq_i$, respectively. 

Another relevant object in our construction is the graded dual of $\wtoba_{\bq}$ which is called the Lusztig algebra $\luq$ associated with $\bq$ in  \cite{AAR-MRL}. Thus, $\luq$ is also an $\N_0^\I$-graded braided Hopf algebra  equipped with a bilinear form $( \, , \, ) : \luq \times \wtoba_{\bq} \rightarrow \ku$ satisfying
\begin{align}
\label{eq:pair-LqBq}
( y, xx') &= ( y^{(2)},x ) ( y^{(1)},x' ),  &
( yy', x ) &= ( y,x^{(2)} ) ( y',x^{(1)} ), 
\end{align}
for all $x,x' \in \wtoba_{\bq}$, $y,y' \in \luq$. Given $n \in \N_0$ and $\beta_j\in\varDelta_+^\bq$, the element $\ee_{\beta_j}^{(n)}\in\luq$ was defined in \cite{AAR-MRL} in such a way that
\begin{align}\label{eq:fdual-definition}
( \ee_{\beta_j}^{(n)}, E_{\beta_1}^{m_1} \cdots E_{\beta_\ell}^{m_{\ell}} ) &= \begin{cases}
1, & m_j=n, \, m_k=0 \text{ for }k\neq j, \\ 0, & \text{otherwise}.
\end{cases}
\end{align}
Then $\luq$ is generated by $\{ \ee_{\alpha_i}\mid i\in\I \}  \cup  \{\ee_{\beta}^{(N_\beta)}\mid\beta\in \fO^{\bq}_+ \}$  and
\begin{align}\label{eq:PBW-basis-Lusztig-alg}
\left\{ \ee_{\beta_1}^{(n_1)} \cdots \ee_{\beta_\ell}^{(n_\ell)} \mid 0\le n_j <\widetilde{N}_{\beta_j},\,1\leq j\leq\ell \right\}
\end{align}
is a basis of $\luq$; cf. \cite[Corollary 4.9 and Lemma 4.6]{AAR-MRL}, our $\eta_{\beta_j}$ is denoted $y_{\beta_j}$ in loc.~cit.

Finally, we recall the subalgebra $Z_\bq$ of $U_\bq$ generated by $E_\beta^{N_\beta}$, $F_\beta^{N_\beta}$, $K_\beta^{\pm N_\beta}$, $L_\beta^{\pm N_\beta}$, $K_\gamma^{\pm \mathtt{N}}$ and $L_\gamma^{\pm \mathtt{N}}$, $\beta\in\fO_+^\bq$ and $\gamma\in\varDelta_+^\bq$, where $\mathtt{N}=\lcm\{N_\beta\mid\beta\in\varDelta_+^\bq\}$, cf. \cite[\S4.5]{AAY}. The subalgebras $Z_\bq^\star$, for $\star\in\{\pm,0,\pm0,\geqslant,\leqslant\}$, are defined in obvious way.

\begin{remark}\label{rem:Zq central}
The subalgebra $Z_\bq$ is central by the assumption \eqref{eq:condition-cartan-roots}, cf. \cite[\S4.1]{Ang-dpn}, and $(K_{\beta}^{N_{\beta}}, L_{\beta}^{N_{\beta}}| \wbeta \in \varPi^{\bq})=(K_{\beta}^{N_{\beta}}, L_{\beta}^{N_{\beta}}| \beta \in \fO_+^\bq)$ by \cite[(6.8)]{AAY}. Moreover, by \cite[Theorem 29]{Ang-dpn} and dually by 
\cite[Theorem 3.3]{AAR-Belgium}, there are extensions of braided Hopf algebras
\[
\begin{tikzcd}[column sep=large, arrow style=tikz, >=To]
	Z_{\bq}^+ \arrow[r, hook, "\iota"] & \wtoba_{\bq} \arrow[r, two heads, "\vpi"] & \toba_{\bq},
\end{tikzcd}
\qquad
\begin{tikzcd}[column sep=large, arrow style=tikz, >=To]
	\toba_{\bq} \arrow[r, hook, "\vpi^*"] & \luq \arrow[r, two heads, "\iota^*"] & U(\n_\bq^+).
\end{tikzcd}
\]
\end{remark}

\section{Divided power integral forms}\label{sec:int-forms-divided-powers}

Inspired by Lusztig's works \cite{L-JAMS,L-book}, we will construct a divided power integral form attached to each braiding matrix in the one-parameter family as an $\Acf$-subalgebra of $U_{\bqf,\Acf}$; this is a technical difference with what was done in \cite{AAY}, where some results stated in loc. cit. are for $\mathtt{A}=\Acf\otimes_{\Z}\Cset$ in place of $\Acf$. More precisely, our object of study in this section is the following algebra. Recall the Hopf skew-pairing from \eqref{eq:skew pairing Ubqf}.

\begin{definition}
The \emph{divided power integral form} is the $\Acf$-subalgebra $\Ub_{\bqf, \Acf}$ of $U_{\bqf, \Acf}$ generated by $U_{\bqf, \Acf}^0$ and $\Ub_{\bqf, \Acf}^{\pm}$, where
\begin{align}\label{eq:Uqres}
\Ub_{\bqf, \Acf}^{-}= \{ y\in U_{\bqf}^{-} | \langle y, U_{\bqf,\Acf}^+ \rangle \subset \Acf \}\quad\mbox{and}\quad
\Ub_{\bqf, \Acf}^{+} = \{ x\in U_{\bqf}^{+} | \langle U_{\bqf,\Acf}^-, x \rangle \subset \Acf \}.
\end{align}
We set $\Ub_{\bqf, \Acf}^{\leqslant}=\Ub_{\bqf, \Acf}^{-}U_{\bqf, \Acf}^{-0}$ and $\Ub_{\bqf, \Acf}^{\geqslant}=\Ub_{\bqf, \Acf}^{+}U_{\bqf, \Acf}^{+0}$. 
\end{definition}

The positive and negative quantum divided powers are the elements
\begin{align}\label{eq:defn-divided-power}
E_{\beta}^{(n)} &=\frac{E_{\beta}^n}{(n)_{\bqf_{\beta\beta}}^!}, &
F_{\beta}^{(n)} &=\frac{F_{\beta}^n}{(n)_{\bqf_{\beta\beta}}^!}, &
\beta\in\varDelta_+^\bqf, & \, \, 0\le n<N_\beta^\bqf.
\end{align}
We recall that $\widetilde{N}_{\beta}=N_\beta^\bqf$ by Remark \ref{rmk:N tilde}, and that $E_{\beta}^n=F_{\beta}^n=0$ if $\beta\notin\fO^{\bqf}_+$ and $n\geq \widetilde{N}_{\beta}$. We will see that \eqref{eq:defn-divided-power} provides a set of distinguished generators of $\Ub_{\bqf, \Acf}$.

We fix a reduced expression $\omega_0^{\bqf}=s_{i_1}^{\bqf} \cdots s_{i_\ell}$ of the longest element of the Weyl groupoid ending in $\bqf$, and set $\beta_k$, $E_{\beta_k}$ and $F_{\beta_k}$ as in \eqref{eq:betak} and \eqref{eq:Ebeta-defn}. Let 
\begin{align*}
\Lt=\left\{\nund=(n_i)\in\N_0^{\ell} \mid 0\le n_j< \widetilde{N}_{\beta_j},\,1\leq j\leq\ell\right\}. 
\end{align*}
For each $\nund\in\Lt$, we set
\begin{align*}
F^{\nund} &=F_{\beta_1}^{n_1} \dots F_{\beta_\ell}^{n_\ell}\in U_{\bqf, \Acf}^-, &
F^{(\nund)} &=F_{\beta_1}^{(n_1)} \dots F_{\beta_\ell}^{(n_\ell)} \in \Ub_{\bqf,\Acf}^{-},
\\
E^{\nund} &=E_{\beta_1}^{n_1} \dots E_{\beta_\ell}^{n_\ell}\in U_{\bqf, \Acf}^+, &
E^{(\nund)} &=E_{\beta_1}^{(n_1)} \dots E_{\beta_\ell}^{(n_\ell)} \in \Ub_{\bqf,\Acf}^{+}.
\end{align*}

Our first result is that quantum divided powers lead to a sort of PBW basis for $\Ub_{\bqf, \Acf}^{\pm}$. Notice that $\Ub_{\bqf, \Acf}^{\pm}$ is an $\Acf$-submodule of $U_{\bqf}^\pm$.

\begin{proposition}\label{prop:dual-A-bases}
The $\Acf$-submodules $\Ub_{\bqf, \Acf}^{-}$ and $\Ub_{\bqf, \Acf}^{+}$ are subalgebras in the corresponding categories of Yetter--Drinfeld modules. Moreover, $\Ub_{\bqf, \Acf}^{+}$ (resp. $\Ub_{\bqf, \Acf}^{-}$) is generated as an $\Acf$-algebra by the positive (resp. negative) quantum divided powers. Also,
\begin{align*}
\langle F_{\beta_1}^{(n_1)} \cdots F_{\beta_\ell}^{(n_\ell)} , E_{\beta_1}^{m_1} \cdots E_{\beta_\ell}^{m_\ell} \rangle =
\delta_{n_1, m_1} \cdots \delta_{m_\ell, n_\ell}=\langle F_{\beta_1}^{n_1} \cdots F_{\beta_\ell}^{n_{\ell}} , E_{\beta_1}^{(m_1)} \cdots E_{\beta_\ell}^{(m_\ell)} \rangle,
\end{align*}
for all $0 \le n_j, m_j< \widetilde{N}_{\beta_j}$, $1\leq j\leq \ell$, and hence the sets
\begin{align}\label{eq:PBW-basis-Ures+-Ures-}
\left\{  F^{(\nund)}\mid\nund\in\Lt \right\}
\quad\text{and}\quad
\left\{  E^{(\nund)}\mid\nund\in\Lt \right\}
\end{align}
are $\Acf$-basis of $\Ub_{\bqf,\Acf}^{-}$ and $\Ub_{\bqf,\Acf}^{+}$, respectively.
\end{proposition}

\begin{proof}
Since $U_{\bqf,\Acf}^{\pm}$ are braided Hopf $\Acf$-subalgebras of $U_{\bqf}^{\pm}$, we deduce that $\Ub_{\bqf,\Acf}^{\pm}$ are subalgebras from the properties of a Hopf skew-pairing. The other claims follows as in \cite[Proposition 5.5 and Corollary 5.6]{AAY}.
\end{proof}

\begin{corollary}\label{coro:convex-order-divided-powers}
Let $1\le j<k\le \ell$, $0\le M< \widetilde{N}_{\beta_k}$ and $0\le N< \widetilde{N}_{\beta_j}$.
Then there exist $\mathtt{a}_{n_{j+1},\ldots,n_{k-1}}\in\Acf$ such that
\begin{align*}
F_{\beta_k}^{(M)}F_{\beta_j}^{(N)} &= \bqf_{M\beta_k \, N\beta_j}F_{\beta_j}^{(N)}F_{\beta_k}^{(M)}
+\sum \mathtt{a}_{n_{j+1},\ldots,n_{k-1}} \, F_{\beta_{j+1}}^{(m_{j+1})}\cdots F_{\beta_{k-1}}^{(m_{k-1})},
\end{align*}
where the sum is over those $(n_{j+1},\ldots,n_{k-1})\in\N_0^{k-j-1}$ such that
$$m_{j+1}\beta_{j+1}+ \cdots + m_{k-1}\beta_{k-1}=M\beta_k+N\beta_j.$$
Also, the set $\left\{  F_{\beta_{\tau(1)}}^{(n_{\tau(1)})} \dots F_{\beta_{\tau(\ell)}}^{(n_{\tau(\ell)})} : 0\le n_{\tau(j)}< \widetilde{N}_{\beta_{\tau(j)}},\,1\leq j\leq\ell \right\}$ is an $\Acf$-basis of $\Ub_{\bqf,\Acf}^{-}$ for any permutation $\tau$ of $\{1,\ldots,\ell\}$.
\end{corollary}
\pf
By \cite[Theorem 4.8]{HY-shapov} there exist $a_{n_{j+1},\ldots,n_{k-1}}\in\Cset(\nu)$ such that
\begin{align*}
F_{\beta_k}^{M}F_{\beta_j}^{N} &= \bqf_{M\beta_k \, N\beta_j}F_{\beta_j}^{N}F_{\beta_k}^{M}
+\sum a_{n_{j+1},\ldots,n_{k-1}} \, F_{\beta_{j+1}}^{m_{j+1}}\cdots F_{\beta_{k-1}}^{m_{k-1}},
\end{align*}
where the sum is over those $(n_{j+1},\ldots,n_{k-1})\in\N_0^{k-j-1}$ as above. From here we get
\begin{align*}
F_{\beta_k}^{(M)}F_{\beta_j}^{(N)} &= \bqf_{M\beta_k \, N\beta_j}F_{\beta_j}^{(N)}F_{\beta_k}^{(M)}
+\sum \mathtt{a}_{n_{j+1},\ldots,n_{k-1}} \, F_{\beta_{j+1}}^{(m_{j+1})}\cdots F_{\beta_{k-1}}^{(m_{k-1})},
\end{align*}
where $\mathtt{a}_{n_{j+1},\ldots,n_{k-1}}= a_{n_{j+1},\ldots,n_{k-1}} \frac{\prod_{i=j+1}^{k-1} (m_i)_{\bqf_{\beta_i\beta_i}}^!}{(M)_{\bqf_{\beta_k\beta_k}}^!(N)_{\bqf_{\beta_j\beta_j}}^!}
\in\Cset(\nu)$. Thus, we have expressed $F_{\beta_k}^{(M)}F_{\beta_j}^{(N)}$ as a $\Cset(\nu)$-linear combination of elements of the basis of $\Ub_{\bqf,\Acf}^{-}$ given in \eqref{eq:PBW-basis-Ures+-Ures-}. On the other hand, Proposition \ref{prop:dual-A-bases}  says that there exists a unique $\Acf$-linear combination of these elements expressing $F_{\beta_k}^{(M)}F_{\beta_j}^{(N)}$. Comparing the two expressions, we deduce that $\mathtt{a}_{n_{j+1},\ldots,n_{k-1}}\in\Acf$.

For the second statement, we combine the first one together with the ideas in the proof of \cite[Theorem 4.9]{HY-shapov}.
\epf

\begin{proposition}\label{prop:U-res-+-Hopf-subalg}
$\Ub_{\bqf,\Acf}^{\pm}$ is a braided Hopf $\Acf$-subalgebra of $U_{\bqf}^{\pm}$.
\end{proposition}
\pf
We already know that they are $\Acf$-subalgebras of $U_{\bqf}^{\pm}$. Thus, it remains to prove that they are $\Acf$-subcoalgebras. We  prove this for $\Ub_{\bqf,\Acf}^{-}$ as the proof for $\Ub_{\bqf,\Acf}^+$ is analogous.

Let $F\in \Ub_{\bqf,\Acf}^{-}$ and $\underline{\Delta}$ be the comultiplication of $\Ub_{\bqf,\Acf}^{-}$. As \eqref{eq:PBW-basis-U+-U-} is a basis of $U_{\bqf}^{-}$, there exist scalars $c_{\mund\nund}\in\Cset(\nu)$ such that
\begin{align*}
\underline{\Delta}(F) &=\sum_{\mund,\nund\in\Lt} c_{\mund\nund} \, F^{\mund} \otimes F^{\nund} =\sum_{\mund,\nund\in\Lt} d_{\mund\nund} \, F^{(\mund)} \otimes F^{(\nund)}.
\end{align*}
Here, $d_{\mund\nund}=c_{\mund\nund}\left(\prod_{j=1}^{\ell} (m_j)_{\bqf_{\beta_j\beta_j}}^!(n_j)_{\bqf_{\beta_j\beta_j}}^! \right) \in\Cset(\nu)$ and all but finitely many of these scalars are zero. As $\{F^{(\nund)} | \nund\in \Lt\}$ is a basis of $\Ub_{\bqf,\Acf}^{-}$, we have to prove that $d_{\mund\nund}\in\Acf$ for all $\mund,\nund\in\Lt$. Since $E^{\nund},E^{\mund}\in U_{\bqf,\Acf}^{+}$ and $U_{\bqf,\Acf}^{+}$ is an $\Acf$-subalgebra, we have that $E^{\nund}E^{\mund}\in U_{\bqf,\Acf}^{+}$ and hence $\langle F,E^{\mund}E^{\nund}\rangle\in\Acf$. Let us compute this value:
\begin{align*}
\langle F,E^{\mund}E^{\nund}\rangle=\langle\underline{\Delta}^{\operatorname{cop}}(F),E^{\mund}E^{\nund}\rangle=\bqf_{\alpha\beta}\langle F^{(1)}, E^{\mund}\rangle \langle F^{(2)},E^{\nund} \rangle
=\bqf_{\alpha\beta}d_{\mund\nund},
\end{align*}
the first equality holds as $\langle-,-\rangle$ is a braided skew-pairing and the second one holds for certain $\alpha,\beta\in\Z^\theta$ by definition of $\underline{\Delta}^{\operatorname{cop}}$. Therefore  $d_{\mund\nund}\in\Acf$ as $\bqf_{\alpha\beta}\in\Acf$.
\epf

As a consequence of this result we get the following two statements.

\begin{corollary}
$\Ub_{\bqf, \Acf}^{\geqslant}$, $\Ub_{\bqf, \Acf}^{\leqslant}$ and  $\Ub_{\bqf,\Acf}$ are Hopf $\Acf$-subalgebras of $U_{\bqf}$. \qed
\end{corollary}

\begin{corollary} 
$\partial_i^{\bqf}(E)\in \Ub_{\bqf,\Acf}^{+}$ for all $E\in \Ub_{\bqf,\Acf}^{+}$ and $i\in\I$. \qed
\end{corollary}

\subsection{The quantum brackets}\label{subsec:quantum-brackets}

One of our aims is to show that $\Ub_{\bqf,\Acf}$ admits a triangular decomposition. While the positive and negative parts are $\Ub_\bqf^\pm$, we need to describe properly the $0$-part, which is larger than $U_{\bqf,\Acf}^0$. We devote this subsection to studying this component. We begin by introducing variations of the bracket symbols defined by Lusztig \cite{L-JAMS,L-book} to express the commutation relations of quantum divided powers.

\begin{definition}
The \emph{quantum bracket} attached to $\beta\in\varDelta^\bqf_+$, $0\le t<\widetilde{N}_{\beta}$ and $n\in\Z$ is
\begin{align}\label{eq:defn-bracket-Kbeta}
\bracket{\beta}{n}{t}_{\bqf}= \frac{\prod_{j=0}^{t-1}\left(\bqf_{\beta}^{j-n}K_{\beta}-L_{\beta}^{-1}\right)}{(t)_{\bqf_{\beta}}^{!}} \quad\mbox{and}\quad
\bracket{\beta}{n}{0}_{\bqf}=1.
\end{align}
When $n=0$, we also write $\bra{\beta}{t}_{\bqf}\coloneqq \bracket{\beta}{0}{t}_{\bqf}$. 

We denote by $\Ub_{\bqf,\Acf}^0$ the $\Acf$-subalgebra generated by $U_{\bqf,\Acf}^0$ and all the quantum brackets $\bracket{\beta}{n}{t}_{\bqf}$ . We will omit the subscript ${\bqf}$ when it is clear from the context.
\end{definition}

\begin{remark}\label{rem:crucial}
The quantum brackets with $\beta\in\fO^{\bqf}_+$ are the most interesting ones. Indeed, $\beta\in\fO^{\bqf}_+$ if and only if $\widetilde{N}_{\beta}=\ord\bqf_{\beta}=\infty$ \cite[Remark 4.1]{AAY}. Instead, if $\beta\notin\fO^{\bqf}_+$, then $N_{\beta}=\widetilde{N}_{\beta}\in\{2,3\}$, and $(k)_{\bqf_{\beta}}^{-1}\in\Acf$ for all $0\leq k< N_{\beta}=\widetilde{N}_{\beta}$, so $\bracket{\beta}{n}{t}_{\bqf}\in U_{\bqf,\Acf}^0$.
\end{remark}

It is not immediate from the definition that the quantum brackets belong to $\Ub_{\bqf,\Acf}$; nevertheless this is the case, see Corollary \ref{le:U-beta-0-generators}. They arise naturally when considering divided powers, as we see in the next result. Let us note before that they behave well under the action of the Lusztig isomorphisms. That is, if $\beta=s_{i_1}^{\bqf}\cdots s_{i_{k-1}}(\gamma)\in\varDelta^{\bqf}_+$ for some $\gamma\in\varDelta^{\bpf}_+$, then $K_{\beta}=T_{i_1}^{\bqf}\cdots T_{i_{k-1}} (K_{\gamma})$, $L_{\beta}=T_{i_1}^{\bqf}\cdots T_{i_{k-1}} (L_{\gamma})$, and $\bqf_{\beta}=\bpf_{\gamma}$, and therefore
\begin{align}\label{eq:T de bracket}
\bracket{\beta}{n}{t}_{\bqf}=T_{i_1}^{\bqf}\cdots T_{i_{k-1}} \bracket{\gamma}{n}{t}_{\bpf}.
\end{align}

\begin{lemma}\label{lemmma:Ealpha-M-Falpha-N}
For each $\beta\in\varDelta^{\bqf}_+$ and $M,N\in\N$, we have that
\begin{align*}
E_{\beta}^{(M)}F_{\beta}^{(N)} &=\sum_{t=0}^{\min\{M,N\}} F_{\beta}^{(N-t)} \bracket{\beta}{M+N-t-1}{t}E_{\beta}^{(M-t)}.
\end{align*}
\end{lemma}
\pf
For any $i\in\I$, this formula for the simple root $\beta=\alpha_i$ is exactly \cite[Lemma 3.7]{HY-shapov}. Now, we fix $\beta=\beta_k=s_{i_1}^{\bqf}\cdots s_{i_{k-1}}(\alpha_{i_k})$. By definition of $E_\beta$ and $F_\beta$, and \eqref{eq:T de bracket} we get
\begin{align*}
E_{\beta}^{(M)} F_{\beta}^{(N)} = & T_{i_1}^{\bqf}\cdots T_{i_{k-1}} \left( E_{i_k}^{(M)}F_{i_k}^{(N)} \right)
\\
= & T_{i_1}^{\bqf}\cdots T_{i_{k-1}} \left( \sum_{t=0}^{\min\{M,N\}} F_{i_k}^{(N-t)} \bracket{\alpha_{i_k}}{M+N-t-1}{t}E_{i_k}^{(M-t)} \right)
\\
= & \sum_{t=0}^{\min\{M,N\}}  T_{i_1}^{\bqf}\cdots T_{i_{k-1}} \left(F_{i_k}^{(N-t)}\right)  T_{i_1}^{\bqf}\cdots T_{i_{k-1}} \bracket{\alpha_{i_k}}{M+N-t-1}{t}
\\
&\qquad\qquad\qquad\qquad\qquad\qquad
\qquad\qquad\qquad \qquad \times T_{i_1}^{\bqf}\cdots T_{i_{k-1}} \left(E_{i_k}^{(M-t)} \right)
\\
= & \sum_{t=0}^{\min\{M,N\}} F_{\beta}^{(N-t)} \bracket{\beta}{M+N-t-1}{t}E_{\beta}^{(M-t)},
\end{align*}
which proves the statement.
\epf

The quantum brackets also satisfy relations analogous to \cite[(g3)--(g5)]{L-JAMS}. Namely, for each $\beta\in\fO^{\bqf}_+$ we have that
\begin{align}
\label{eq:g3}
\bra{\beta}{t}\bracket{\beta}{-t}{s}&=\binom{t+s}{t}_{\bqf_{\beta}}\bra{\beta}{t+s}\quad (s,t\geq0),\\
\noalign{\smallskip}
\label{eq:g4}
\bracket{\beta}{n+1}{t}&=\bracket{\beta}{n}{t}+\bqf_\beta^{-n-1}(1-\bqf_\beta)K_\beta\bracket{\beta}{n}{t-1}\quad(n\in\Z,\,t\geq1),\\
\noalign{\smallskip}
\label{eq:beta01}\bra{\beta}{1}&=K_\beta-L_\beta^{-1}.
\end{align}
Moreover, they obey relations analogous to \cite[(g8)--(g10)]{L-JAMS} as we see below. In particular, we deduce that all quantum brackets can be generated from those with $n=0$.

\begin{lemma}\label{lem:bracket-gen-by-case-n=0}
Let $\beta\in\fO^{\bqf}_+$, $n\in\Z$ and $M,N\in\N_0$. If $M\geq1$ and $N\geq0$, then
\begin{align}\label{eq:g8}
\binom{M+N}{M}_{\bqf_{\beta}}\bra{\beta}{M+N}=\sum_{j=0}^N (\bqf_{\beta}-1)^{j} \binom{M+j-1}{j}_{\bqf_{\beta}} K_{\beta}^{j} \bra{\beta}{M}\bra{\beta}{N-j}.
\end{align}
If $n<0$, then
\begin{align}\label{eq:bracket-n<0}
\bracket{\beta}{n}{N}=\sum_{j=0}^N (\bqf_{\beta}-1)^{N-j} \binom{-n+N-j-1}{N-j}_{\bqf_{\beta}} K_{\beta}^{N-j} \bra{\beta}{j}.
\end{align}
If $n\ge 0$, then
\begin{align}\label{eq:bracket-nge0}
\bracket{\beta}{n}{N}=\sum_{j=0}^N \bqf_{\beta}^{\binom{N-j}{2}-n(N-j)}(1-\bqf_{\beta})^{N-j} \binom{n}{N-j}_{\bqf_{\beta}} K_{\beta}^{N-j} \bra{\beta}{j}.
\end{align}
\end{lemma}

\begin{proof}
 The equality \eqref{eq:g8} follows from \eqref{eq:g3} and \eqref{eq:bracket-n<0}. In order to prove the other two equalities, we first claim that for all $\beta\in\fO^{\bqf}_+$, $N\in\N$, $n\in\Z$, there exist $f_{\beta,n}^N(j)\in\Cset(\nu)$ such that
\begin{align}\label{eq:bracket-n-general-formula}
\bracket{\beta}{n}{N}=\sum_{j=0}^N f_{\beta,n}^N(j) \, K_{\beta}^{N-j} \bra{\beta}{j},
\end{align}
where $f_{\beta,n}^N(N)=1$, and the following recursive formula holds:
\begin{align}\label{eq:bracket-n-recursive-formula-coeff}
\begin{aligned}
f_{\beta,n}^{N+1}(0) &= \frac{f_{\beta,n}^{N}(0) (\bqf_{\beta}^{N-n}-1)}{(N+1)_{\bqf_{\beta}}}, \\
f_{\beta,n}^{N+1}(j) &= \frac{(\bqf_{\beta}^{N-n}-\bqf_{\beta}^j) f_{\beta,n}^{N}(j) + (j)_{\bqf_{\beta}} f_{\beta,n}^{N}(j-1)}{(N+1)_{\bqf_{\beta}}}, &
&1\le j \le N.
\end{aligned}
\end{align}
The proof is by induction on $N$. Indeed, for $N=1$ we have
\begin{align*}
\bracket{\beta}{n}{1}=\bqf_{\beta}^{-n} K_{\beta}-L_{\beta}^{-1}=(\bqf_{\beta}^{-n}-1) K_{\beta}+\bra{\beta}{1}.
\end{align*}
Assume now that \eqref{eq:bracket-n-general-formula} holds for $N$. Then, using the inductive hypothesis,
\begin{align*}
\bracket{\beta}{n}{N+1} &= \tfrac{\bqf_{\beta}^{N-n}K_{\beta}-L^{-1}_{\beta}}{(N+1)_{\bqf_{\beta}}} \bracket{\beta}{n}{N}
= \sum_{j=0}^N \tfrac{\bqf_{\beta}^{N-n}K_{\beta}-\bqf_{\beta}^{j}K_{\beta}+\bqf_{\beta}^{j}K_{\beta}-L^{-1}_{\beta}}{(N+1)_{\bqf_{\beta}}} f_{\beta,n}^N(j) \, K_{\beta}^{N-j} \bra{\beta}{j}
\\
&= \sum_{j=0}^N \tfrac{\bqf_{\beta}^{N-n}-\bqf_{\beta}^{j}}{(N+1)_{\bqf_{\beta}}} f_{\beta,n}^N(j) \, K_{\beta}^{N+1-j} \bra{\beta}{j}
+\sum_{j=0}^N \tfrac{(j+1)_{\bqf_{\beta}}}{(N+1)_{\bqf_{\beta}}} f_{\beta,n}^N(j) \, K_{\beta}^{N-j} \bra{\beta}{j+1}
\\
&= \tfrac{(\bqf_{\beta}^{N-n}-1)f_{\beta,n}^N(0)}{(N+1)_{\bqf_{\beta}}}  K_{\beta}^{N+1}
+\sum_{j=1}^N \tfrac{(\bqf_{\beta}^{N-n}-\bqf_{\beta}^{j}) f_{\beta,n}^N(j)+(j)_{\bqf_{\beta}} f_{\beta,n}^N(j-1)}{(N+1)_{\bqf_{\beta}}}  \, K_{\beta}^{N+1-j} \bra{\beta}{j}
\\
& \qquad + f_{\beta,n}^N(N) \, \bra{\beta}{N+1},
\end{align*}
and the inductive step follows.

Thus it remains to compute the terms $f_{\beta,n}^N(j)$ using \eqref{eq:bracket-n-recursive-formula-coeff}. First we consider $j=0$:
\begin{align*}
f_{\beta,n}^{N}(0) &= \frac{(\bqf_{\beta}^{N-1-n}-1) (\bqf_{\beta}^{N-2-n}-1)\cdots (\bqf_{\beta}^{-n}-1)}{(N)_{\bqf_{\beta}}^{!}}
=\frac{\prod_{\ell=0}^{N-1}\bqf_{\beta}^{\ell-n}-1}{(N)_{\bqf_{\beta}}^!}
\\
&=\begin{cases}
(\bqf_{\beta}-1)^N \binom{-n+N-1}{N}_{\bqf_{\beta}}, & n<0,
\\
\bqf_{\beta}^{\binom{N}{2}-nN}(1-\bqf_{\beta})^N \binom{n}{N}_{\bqf_{\beta}}, & n\ge 0.
\end{cases}
\end{align*}
Next we claim that
\begin{align*}
f_{\beta,n}^{N}(j) &= f_{\beta,n}^{N-j}(0)=\frac{\prod_{\ell=0}^{N-j-1}\bqf_{\beta}^{\ell-n}-1}{(N-j)_{\bqf_{\beta}}^!} & &\text{ for all }j\in\I_N.
\end{align*}
We argue again by induction: the case $N=1$ is straightforward. Assume that this holds for $N$. Using the recursive formula \eqref{eq:bracket-n-recursive-formula-coeff} we have
\begin{align*}
f_{\beta,n}^{N+1}(j) &= \frac{\bqf_{\beta}^j(\bqf_{\beta}^{N-n-j}-1) \frac{\prod_{\ell=0}^{N-j-1}\bqf_{\beta}^{\ell-n}-1}{(N-j)_{\bqf_{\beta}}^!} + (j)_{\bqf_{\beta}} \frac{\prod_{\ell=0}^{N-j}\bqf_{\beta}^{\ell-n}-1}{(N-j+1)_{\bqf_{\beta}}^!}}{(N+1)_{\bqf_{\beta}}}
\\
&= \frac{\prod_{\ell=0}^{N-j}\bqf_{\beta}^{\ell-n}-1}{(N+1-j)_{\bqf_{\beta}}^!}
\frac{\bqf_{\beta}^j(N-j+1)_{\bqf_{\beta}}  + (j)_{\bqf_{\beta}}}{(N+1)_{\bqf_{\beta}}}
= \frac{\prod_{\ell=0}^{N-j}\bqf_{\beta}^{\ell-n}-1}{(N+1-j)_{\bqf_{\beta}}^!}=f_{\beta,n}^{N+1-j}(0),
\end{align*}
and the statement follows.
\end{proof}

The following is an immediate consequence of the above lemma.

\begin{corollary}\label{cor:U-0-generators} As an $\Acf$-algebra, $\Ub_{\bqf,\Acf}^0$ is generated by $K_i$, $L_i$, $i\in\I$, and the quantum brackets $\bra{\beta}{t}$, $\beta\in\fO^{\bqf}_+$, $t\in\N_0$.\qed
\end{corollary}

We next see that in fact all quantum brackets belong to $\Ub_{\bqf,\Acf}$.
Given $\beta\in\fO^{\bqf}_+$, we denote by $\widetilde\Ub_\beta$ the $\Acf$-subalgebra generated by $K_\beta$, $L_\beta^{-1}$, $E_\beta^{(n)}$ and $F_\beta^{(n)}$, $n\in\N_0$. 

\begin{corollary}\label{le:U-beta-0-generators}
Let $\beta\in\fO^{\bqf}_+$. Then $\bracket{\beta}{n}{t}\in\widetilde\Ub_\beta$ for all $n\in\Z$ and $t\in\N_0$. Therefore $\Ub_{\bqf,\Acf}^0$ is contained in $\Ub_{\bqf,\Acf}$.
\end{corollary}

\begin{proof}
We proceed by induction on $t$. Clearly, $\bracket{\beta}{n}{0}, \bracket{\beta}{n}{1}\in\widetilde\Ub_\beta$ for all $n\in\Z$. Now, by Lemma \ref{lemmma:Ealpha-M-Falpha-N} with $M=N=t+1$ and the inductive hypothesis, we have that $\bracket{\beta}{t}{t+1}\in\widetilde\Ub_\beta$. By applying repeatedly \eqref{eq:g4} and the inductive hypothesis, we deduce that $\bra{\beta}{t+1}\in\widetilde\Ub_\beta$. Therefore $\bracket{\beta}{n}{t+1}\in\widetilde\Ub_\beta$ for all $n\in\Z$ by Lemma \ref{lem:bracket-gen-by-case-n=0} and the inductive hypothesis. This finishes the proof.
\end{proof}

For $\beta\in\fO^{\bqf}_+$, let us denote $\widetilde\Ub_\beta^0$ the $\Acf$-subalgebra generated by $K_\beta$ and $\bra{\beta}{t}$, $t\in\N_{0}$. Thus, Lemma \ref{lem:bracket-gen-by-case-n=0} implies that $\bracket{\beta}{n}{t}$ belongs to $\widetilde\Ub_\beta^0$ for all $n\in\Zset$ and $t\in\N_{0}$. We also note that $L_\beta^{-1}=-K_\beta-\bra{\beta}{1}$ belongs to $\widetilde\Ub_\beta^0$.

\begin{lemma}\label{le:tilde U 0 beta}
Let $\beta\in\fO^{\bqf}_+$ and $t\in\N_{0}$. Then
\begin{align}\label{eq:coproduct-bracket-int-form}
\Delta\bra{\beta}{t}=\sum_{0\leq s\leq t}L_\beta^{-s}\bra{\beta}{t-s}\otimes K_\beta^{t-s}\bra{\beta}{s}
\end{align}
and hence $\widetilde\Ub_\beta^0$ is an $\Acf$-subbialgebra of $U_{\bqf,\Acf}$. Moreover,
\begin{align}\label{eq:A gen of Vbeta}
\left\{K_{\beta}^{s}\,\bra{\beta}{t}\right\}_{s,t\in\N_0}
\end{align}
is an $\Acf$-basis of $\widetilde\Ub_\beta^0$.
\end{lemma}

\begin{proof}
The formula for the comultiplication of $\bra{\beta}{t}$ follows by induction on $t$. This implies that $\widetilde\Ub_\beta^0$ is an $\Acf$-subbialgebra.

Working in $\widetilde\Ub_\beta^0\otimes_\Acf\ku(\nu)$, there exist $a_{i,j}\in\ku(\nu)$ such that
\begin{align*}
K_{\beta}^{s}\,\bra{\beta}{t}=a_{s,t}K_\beta^sL_\beta^{-t}+\sum_{\substack{ s<i\leq s+t\\0\leq j<t}}a_{i,j}K_\beta^iL_\beta^{-j}
\end{align*}
Then the set $\left\{K_{\beta}^{s}\,\bra{\beta}{t}\right\}_{s,t\in\N_0}$ is linearly independent over $\ku(\nu)$ and hence over $\Acf$. In order to show that they also generate $\widetilde\Ub_\beta^0$, it is enough to see that $\bra{\beta}{t}\bra{\beta}{u}$ belongs to the $\Acf$-submodule generated by them. We prove this claim by induction on $u$. First, we deduce from \eqref{eq:defn-bracket-Kbeta} that
\begin{align*}
\bra{\beta}{t}\bra{\beta}{1}&=
(1-\bqf_{\beta}^t)K_\beta\bra{\beta}{t}+(t+1)_{\bqf_{\beta}}\bra{\beta}{t+1}
\end{align*}
for all $t\in\N$. Next, \eqref{eq:g8} implies that
\begin{align*}
\bra{\beta}{t}\bra{\beta}{u}=
\binom{t+u}{t}_{\bqf_{\beta}}&\bra{\beta}{t+u}\\
-\sum_{j=0}^{u-1}& (\bqf_{\beta}-1)^{u-j} \binom{t+u-j-1}{u-j}_{\bqf_{\beta}} K_{\beta}^{u-j} \bra{\beta}{t}\bra{\beta}{j}
\end{align*}
for all $t\in\N$, and our claim follows by the inductive hypothesis.
\end{proof}

\begin{remark}\label{rem:ZPi-proper-subgroup}
	
One may wonder whether a basis of $\Ub_{\bqf,\Acf}^0$ can be described proceeding as in Lemma \ref{le:tilde U 0 beta}, using that $\varPi^\bqf$ is a $\Z$-basis of the root lattice of $\g_{\bq}$. Here there are two problems to be solved.
The first one is that we do not know if for $\wgamma\notin\varPi^{\bqf}$, $\bra{\gamma}{n}$ belongs to the $\Acf$-subalgebra $\widetilde\Ub_{\bqf,\Acf}^0$ of $\Ub_{\bqf,\Acf}^0$ generated by $K_\beta$ and $\bra{\beta}{t}$, $\wbeta\in\varPi^\bqf$, $t\in\N_0$.

The second problem comes from the fact that in some cases $\Z\varPi^\bqf$ is a proper $\Z$-submodule of $\Z^\theta$. More precisely, see Lemma \ref{le:U-0-eval-A-repr},
\begin{itemize}
\item If $\bqf$ is of type $\supera{m}{n}$, then $\Z^\theta/\Z\varPi^\bqf \simeq \Z$;
\item If $\bqf$ is of type $\superd{m}{n}$, $\mathbf{D}(2,1;\alpha)$, $\mathbf{F}(4)$ or $\texttt{br}(2)$, then $\Z^\theta/\Z\varPi^\bqf \simeq \Z/2\Z$;
\item If $\bqf$ is of type $\texttt{wk}(4)$, then $\Z^\theta/\Z\varPi^\bqf \simeq \Z/3\Z$.
\end{itemize}
In these cases, $\langle K_{\beta}^{\pm1}, L_{\beta}^{\pm1} \mid \wbeta\in\varPi^\bqf \rangle\subsetneq\langle K_{i}^{\pm1}, L_{i}^{\pm1} \mid i\in\I \rangle$ and hence $\widetilde\Ub_{\bqf,\Acf}^0 \subsetneq \Ub_{\bqf,\Acf}^0$.
\end{remark}

We finish this subsection with the description of a suitable family of $\Acf$-representations of $\Ub_{\bqf,\Acf}^0$ of rank one which will be useful in the following section. 

\begin{lemma}\label{le:U-0-A-repr}
Suppose $\chi:U_{\bqf}^0\to\Cset(\nu)$ is a $\Cset(\nu)$-algebra map such that $\chi(K_i),\chi(L_i)\in\Acf$, $i\in\I$, and $\chi(K_{\beta})=\bqf_{\beta}^{m_{\beta}}$ and $\chi(L_{\beta})=\bqf_{\beta}^{n_{\beta}}$, $\wbeta\in\varPi^{\bqf}$, for some $m_{\beta},n_{\beta}\in\N_0$. Then $\chi$ restricts to an $\Acf$-algebra map $\chi:\Ub_{\bqf,\Acf}^0\to\Acf$.
\end{lemma}

\pf

By Remark \ref{rem:Zq central},
for each $\beta\in\fO^{\bqf}_+$ there exists $m_{\beta},n_{\beta}\in\N_0$ such that $\chi(K_{\beta})=\bqf_{\beta}^{m_{\beta}}$, $\chi(L_{\beta})=\bqf_{\beta}^{n_{\beta}}$. Thus, by direct computation, we have that
\begin{align*}
\chi \left(\bra{\beta}{t}\right) &
=\frac{\prod_{j=0}^{t-1}\bqf_{\beta}^{m_{\beta}+j} -\bqf_{\beta}^{-n_{\beta}}}{(t)_{\bqf_{\beta}}^!}
=\bqf_{\beta}^{-n_{\beta}t}(\bqf_{\beta}-1)^t \binom{m_{\beta}+n_{\beta}+t-1}{t}_{\bqf_{\beta}}\in\Acf,
\end{align*}
for all $\beta\in\fO^{\bqf}_+$ and $t\in\N_0$. Since $\Ub_{\bqf,\Acf}^0$ is generated by the quantum brackets $\bra{\beta}{t}$ and the elements $K_i$ and $L_i$, the statement follows. 
\epf

\subsection{Triangular decomposition}\label{subsec:div-powers-triang-decomp}

This is the main result of the section.

\begin{theorem}\label{thm:triang-decomp-divided-powers}
The multiplication gives an isomorphism of $\Zset^{\I}$-graded $\Acf$-modules
\begin{align*}
\Ub_{\bqf,\Acf}^{+} \otimes_{\Acf} \Ub_{\bqf,\Acf}^{0} \otimes_{\Acf} \Ub_{\bqf,\Acf}^{-} \overset{\sim}{\longrightarrow} \Ub_{\bqf,\Acf}
\end{align*}
and $\Ub_{\bqf,\Acf}^{\pm}\Ub_{\bqf,\Acf}^0=\Ub_{\bqf,\Acf}^0\Ub_{\bqf,\Acf}^{\pm}$. Moreover, for each $i\in\I$ the Lusztig isomorphism restricts to an $\Acf$-algebra isomorphism
\begin{align*}
T_i^{\bqf}:\Ub_{\rho_i (\bqf),\Acf} \to \Ub_{\bqf,\Acf}.
\end{align*}
\end{theorem}

\pf

We begin by proving that $\Ub_{\bqf,\Acf}^{\pm}\Ub_{\bqf,\Acf}^0=\Ub_{\bqf,\Acf}^0\Ub_{\bqf,\Acf}^{\pm}$.
Let $i\in\I$ and $\gamma\in\varDelta_+^\bqf$. The defining relations of $U_\bqf$ imply that
\begin{align}\label{eq:KiLi-commutes-Fbeta-Ebeta}
K_iF_{\gamma}^{(N)} &= \bqf_{\alpha_i\,\gamma}^{-N} F_{\gamma}^{(N)}K_i,& 
K_iE_{\gamma}^{(N)} &= \bqf_{\alpha_i\,\gamma}^{N} E_{\gamma}^{(N)}K_i,&
\\
L_iF_{\gamma}^{(N)} &= \bqf_{\gamma \, \alpha_i}^{-N} F_{\gamma}^{(N)}L_i,& 
L_iE_{\gamma}^{(N)} &= \bqf_{\gamma \, \alpha_i}^{N} E_{\gamma}^{(N)}L_i.&
\end{align}
Let $\beta\in\fO^{\bqf}_+$ and recall $a_{\beta\gamma}$ from \eqref{eq:cartan-root}. Then
\begin{align}\label{eq:bracket-commutes-Fbeta-Ebeta}
\begin{split}
\bracket{\beta}{n}{t}F_{\gamma}^{(N)} &= \bqf_{\gamma\beta}^{tN}  F_{\gamma}^{(N)}\bracket{\beta}{n+N a_{\beta\gamma}}{t}\quad\mbox{and}\\
\noalign{\smallskip}
\bracket{\beta}{n}{t}E_{\gamma}^{(N)} &= \bqf_{\gamma\beta}^{-tN} E_{\gamma}^{(N)}\bracket{\beta}{n-N a_{\beta\gamma}}{t}.
\end{split}
\end{align}
by direct computation. 

Regarding the multiplication map, \eqref{eq:triangular decomposition of Ubqf} implies that it is injective because $\Ub_{\bqf,\Acf}^{\star} \subset U_{\bqf}^{\star}$ for $\star \in\{+,-,0\}$. Thus, we have to show that it is surjective. Set $W_{\bqf}:=\Ub_{\bqf,\Acf}^{-}\Ub_{\bqf,\Acf}^{0}\Ub_{\bqf,\Acf}^{+}$. As $1\in W_{\bqf}$, it is enough to check that $W_{\bqf}$ is a left ideal. It is immediate that $F_{\beta}^{(M)} W_{\bqf} \subseteq W_{\bqf}$ for all $\beta\in\varDelta_+^{\bqf}$ and all $M\in\N_0$ as $\Ub_{\bqf,\Acf}^{-}$ is an $\Acf$-subalgebra. Also, $X W_{\bqf} \subseteq W_{\bqf}$ for $X\in\{K_i^{\pm1},L_i^{\pm1},\bracket{\alpha}{\beta}{t}\}$ by \eqref{eq:bracket-commutes-Fbeta-Ebeta}. In order to see that $\Ub_{\bqf,\Acf}^{+}W_\bqf\subseteq W_\bqf$, we first claim that $T_i^{\bqf}(W_{\rho_i(\bqf)})=W_{\bqf}$ for all $i\in\I$.

Indeed, let us fix $i\in\I$ and a reduced expression 
$w_0^{\bqf}=s_{i'_1}^{\bqf} \cdots s_{i'_{\ell}}$ with $i'_1=i$ and set $\bpf:=\rho_i(\bqf)$. 
As $s_{i'_2}^{\bpf} \cdots s_{i'_{\ell}}$ is reduced, it can be extended to a reduced expression
$w_0^{\bpf}=s_{i'_2}^{\bpf} \cdots s_{i'_{\ell}}s_{i'_{\ell+1}}$ for some $i'_{\ell+1}\in\I$. For $1\leq j\leq \ell$, let
\begin{align*}
\beta'_j=s_{i_1}^{\bqf} \cdots s_{i_{j-1}}(\alpha_{i_j})\in\varDelta_+^{\bqf}\quad\mbox{and}\quad
\gamma_j=s_{i_2}^{\bpf} \cdots s_{i_{j}}(\alpha_{i_{j+1}})\in\varDelta_+^{\bpf}.
\end{align*}
Then $\varDelta_+^{\bqf}=\{\beta'_j \mid 1\leq j\leq\ell\}$ and  $\varDelta_+^{\bpf}=\{\gamma_j \mid 1\leq j\leq\ell\}$.
Since $\beta_1'=\alpha_i$ and $s_i^{\bq}(\gamma_j)=\beta'_{j+1}$ for all $1\leq j\leq\ell-1$, we have that $\{\gamma_j \mid 1\leq j\leq\ell-1\} = \varDelta_+^{\bpf}\setminus\{\alpha_i\}$ and $\gamma_{\ell}=\alpha_i$ by \cite[Proposition 1]{H-Weyl gpd}.
We set accordingly
\begin{align*}
E_{\beta_j'} &= T_{i'_1}^{\bqf} \cdots T_{i'_{j-1}}(E_{i'_j})\in \Ub_{\bqf,\Acf}^{+}, & 
F_{\beta_j'} &= T_{i'_1}^{\bqf} \cdots T_{i'_{j-1}}(F_{i'_j})\in \Ub_{\bqf,\Acf}^{-},
\\
\underline{E}_{\gamma_j} &:= T_{i'_2}^{\bpf} \cdots T_{i'_{j}}(E_{i'_{j+1}})\in \Ub_{\bpf,\Acf}^{+}, &
\underline{F}_{\gamma_j} &:= T_{i'_2}^{\bpf} \cdots T_{i'_{j}}(F_{i'_{j+1}})\in \Ub_{\bpf,\Acf}^{-}.
\end{align*}
Notice that $T_i^{\bqf}(\underline{E}_{\gamma_j})=E_{\beta'_{j+1}}$ and $T_i^{\bqf}(\underline{F}_{\gamma_j})=F_{\beta'_{j+1}}$ for all $1\leq j\leq\ell-1$.
Also,
\begin{align*}
\underline{E}_{\gamma_{\ell}}=f\, \underline{E}_{i}\quad\mbox{and}\quad \underline{F}_{\gamma_{\ell}}=g\, \underline{F}_{i}
\end{align*}
for some $f,g\in\Acf$ as $\gamma_{\ell}=\alpha_i$. Then, by Proposition \ref{prop:dual-A-bases} and Corollary \ref{coro:convex-order-divided-powers}, $W_{\bqf}$ and $W_{\bpf}$ are spanned by the sets
\begin{align*}
&\left\{  F_{\beta'_{1}}^{(n_{1})} \cdots F_{\beta'_{\ell}}^{(n_{\ell})} \,  b  \, E_{\beta'_{\ell}}^{(m_{\ell})} \cdots E_{\beta'_1}^{(m_1)} \mid 0\le m_j,n_j< \widetilde{N}_{\beta_j}, b\in \Ub_{\bqf,\Acf}^{0} \right\}\quad\mbox{and}
\\
&\left\{  \underline{F}_{\gamma_{1}}^{(u_{1})} \cdots \underline{F}_{\gamma_{\ell}}^{(u_{\ell})} \,  \underline{b} \, \underline{E}_{\gamma_{\ell}}^{(v_{\ell})} \cdots \underline{E}_{\gamma_1}^{(v_1)} \mid 0\le u_j,v_j< \widetilde{N}_{\gamma_j}, \underline{b} \in \Ub_{\bpf,\Acf}^{0} \right\},
\end{align*}
respectively. Moreover, although the PBW generators $\{F_{\beta'_{1}}, \dots, F_{\beta'_{\ell}}\}$ and $\{E_{\beta'_{1}}, \dots, F_{\beta'_{\ell}}\}$ are not necessarily equal to $\{F_{\beta_{1}}, \dots, F_{\beta_{\ell}}\}$ and $\{E_{\beta_{1}}, \dots, E_{\beta_{\ell}}\}$, they give rise to the same $\Acf$-subalgebras $\Ub^-_{\bqf,\Acf}$ and $\Ub^+_{\bqf,\Acf}$ thanks to \cite[Theorem 6.20]{H-Lusztig}.

We now fix $0\le u_j,v_j< \widetilde{N}_{\gamma_j}$ and $\underline{b} \in \Ub_{\bpf,\Acf}^{0}$. Then there exists $\underline{b}'\in \Ub_{\bpf,\Acf}^{0}$ such that
$
\underline{b} \underline{E}_i^{(v_\ell)}= \underline{E}_i^{(v_\ell)}\underline{b}'
$
by \eqref{eq:bracket-commutes-Fbeta-Ebeta} and $b':=T_i^{\bqf}(\underline{b}')\in \Ub_{\bqf,\Acf}^{0}$ by \eqref{eq:T de bracket}. Thus, using the computations above, Lemma \ref{lemmma:Ealpha-M-Falpha-N} and \eqref{eq:bracket-commutes-Fbeta-Ebeta} we have that
\begin{align*}
T_i^{\bqf} & \big(\underline{F}_{\gamma_{1}}^{(u_{1})} \cdots \underline{F}_{\gamma_{\ell}}^{(u_{\ell})} \,  \underline{b} \, \underline{E}_{\gamma_{\ell}}^{(v_{\ell})} \cdots \underline{E}_{\gamma_1}^{(v_1)} \big)
\\
&= f^{v_{\ell}}g^{u_{\ell}} \,
T_i^{\bqf} \big( \underline{F}_{\gamma_{1}}^{(u_{1})} \cdots \underline{F}_{\gamma_{\ell-1}}^{(u_{\ell-1})} \underline{F}_{i}^{(u_{\ell})}
\underline{E}_{i}^{(v_{\ell})} \,  \underline{b}' \, \underline{E}_{\gamma_{\ell-1}}^{(v_{\ell-1})} \cdots \underline{E}_{\gamma_1}^{(v_1)} \big)
\\
&= q_{ii}^{-\binom{u_{\ell}}{2}-\binom{v_{\ell}}{2}}f^{v_{\ell}}g^{u_{\ell}} \,
F_{\beta'_2}^{(u_1)} \cdots F_{\beta'_{\ell}}^{(u_{\ell-1})}
K_i^{-u_{\ell}} E_{i}^{(u_{\ell})}
F_{i}^{(v_{\ell})}L_i^{v_{\ell}}
\, b' \, E_{\beta'_{\ell}}^{(v_{\ell-1})} \cdots E_{\beta'_2}^{(v_1)}
\\
&\in 
\sum_{t\ge 0} \Acf \, F_{\beta'_2}^{(u_1)} \cdots F_{\beta'_{\ell}}^{(u_{\ell-1})}
K_i^{-u_{\ell}}
F_{i}^{(v_{\ell}-t)} \bracket{\alpha_i}{u_{\ell}+v_{\ell}-t-1}{t}E_{i}^{(u_{\ell}-t)}
\\
& \qquad\qquad\qquad\qquad\qquad\times L_i^{v_{\ell}} b' \, E_{\beta'_{\ell}}^{(v_{\ell-1})} \cdots E_{\beta'_2}^{(v_1)}
\\
& \subseteq \sum_{t\ge 0} F_{\beta'_2}^{(u_1)} \cdots F_{\beta'_{\ell}}^{(u_{\ell-1})}F_{i}^{(v_{\ell}-t)}
\Ub_{\bqf,\Acf}^{0} 
E_{i}^{(u_{\ell}-t)}E_{\beta'_{\ell}}^{(v_{\ell-1})} \cdots E_{\beta'_2}^{(v_1)}.
\end{align*}
As $\Ub_{\bqf,\Acf}^{\pm}$ are $\Acf$-subalgebras, we have that
\begin{align*}
F_{\beta'_2}^{(u_1)} \cdots F_{\beta'_{\ell}}^{(u_{\ell-1})}F_{i}^{(v_{\ell}-t)}  \in \Ub_{\bqf,\Acf}^{-}\quad\mbox{and}\quad
E_{i}^{(u_{\ell}-t)}E_{\beta'_{\ell}}^{(v_{\ell-1})} \cdots E_{\beta'_2}^{(v_1)}  \in \Ub_{\bqf,\Acf}^{+}
\end{align*}
for all $t\ge 0$ and hence $T_i^{\bqf} \big(\underline{F}_{\gamma_{1}}^{(u_{1})} \cdots \underline{F}_{\gamma_{\ell}}^{(u_{\ell})} \,  \underline{b} \, \underline{E}_{\gamma_{\ell}}^{(v_{\ell})} \cdots \underline{E}_{\gamma_1}^{(v_1)} \big) \in \Ub_{\bqf,\Acf}^{-}\Ub_{\bqf,\Acf}^{0}\Ub_{\bqf,\Acf}^{+}=W_{\bqf}$. Hence, $T_i^{\bqf}(W_{\bpf})\subseteq W_{\bqf}$.
An analogous argument shows that $T_i^{\bpf}(W_{\bqf})\subseteq W_{\bpf}$, and then $T_i^{\bqf}(W_{\bpf})=W_{\bqf}$ as claimed.

We are ready to check that $\Ub_{\bqf,\Acf}^{+}W_\bqf\subseteq W_\bqf$. For $\beta\in\varDelta_+^\bqf$ and $M\in\N_0$, we have that 
\begin{align*}
E_{\beta}^{(M)}&=T_{i_1}^{\bqf}\cdots T_{i_{j-1}}(E_{i_j}^{(M)}) = \bqf_{i_ji_j}^{-\binom{M}{2}} T_{i_1}^{\bqf}\cdots T_{i_{j}}(K_{i_j}^{-M}F_{i_j}^{(M)})
\end{align*}
by \eqref{eq:lusztig-isom-defn}. Therefore, by our previous claim and the stability of $W_{\bpf}$ under left multiplication by elements of $\Ub_{\bqf,\Acf}^{-}$ and $\Ub_{\bqf,\Acf}^{0}$ already proved, we deduce that
\begin{align*}
E_{\beta}^{(M)} W_{\bqf} &= \bqf_{i_ji_j}^{-\binom{M}{2}} T_{i_1}^{\bqf}\cdots T_{i_{j}}(K_{i_j}^{-M}F_{i_j}^{(M)}) T_{i_1}^{\bqf}\cdots T_{i_{j}}(W_{\bpf})
\\
& \subseteq  T_{i_1}^{\bqf}\cdots T_{i_{j}}(K_i^{-M}F_{i_j}^{(M)} W_{\bpf}) \subseteq T_{i_1}^{\bqf}\cdots T_{i_{j}}(W_{\bpf}) =W_{\bqf}
\end{align*}
for $\bpf=\rho_{i_j}\cdots\rho_{i_1}(\bqf)$, as desired. 
Hence $W_{\bqf}$ is a left ideal and the multiplication map is also surjective.

Finally, since $T_i^{\bqf}(W_{\rho_i( \bqf)})=W_{\bqf}=\Ub_{\bqf,\Acf}$, the last claim of the statement holds.
\epf

\section{Quantum group of divided powers}\label{sec:specialization}

Along the lines of Lusztig's ideas \cite{L-JAMS,L-book}, in this section we study what we call the \emph{quantum group of divided powers} of $\bq$, a suitable Hopf algebra quotient of the specialization of $\Ub_{\bqf,\Acf}$ at $\bq$. Explicitly, it is the Hopf algebra
\begin{align}\label{eq:U bq}
	\Ub_{\bq}= \Ub_{\bqf,\Acf}\otimes_{\Acf} \Bbbk/\Ib_{\bq},
\end{align}
where $\Ib_{\bq}$ is the central Hopf ideal of $\Ub_{\bqf,\Acf}\otimes_{\Acf} \Bbbk$ generated by $K_\beta^{N_\beta}-1$, $\wbeta\in\varPi^\bq$, recall Remark~\ref{rem:Zq central}.

By abuse of notation, we simply write $x$ for any element $x\otimes 1$ in the specialization, or even the class in the quotient algebra $\Ub_{\bq}$. 
We denote by $\Ub_\bq^\star$ the subalgebra of $\Ub_\bq$ obtained as image of $\Ub_{\bqf,\Acf}^\star$ for $\star\in\{\pm,0,\pm0,\geqslant,\leqslant\}$. We also consider the $\Bbbk$-subalgebra $\uf_\bq$ of $\Ub_{\bq}$ generated by $E_i$, $F_i$, $K_i^{\pm1}$ and $L_i^{\pm1}$, $i\in\I$, and its $\Bbbk$-subalgebras $\uf^\star_\bq$ defined in the obvious way, for $\star\in\{\pm,0,\pm0,\geqslant,\leqslant\}$.

\subsection{First properties}\label{subsec:specialization-first-prop}
We start by observing that the specialization at $\bq$ respects the triangular decomposition. That is, Theorem \ref{thm:triang-decomp-divided-powers} immediately implies the following. 

\begin{proposition}\label{prop:triangular decomposition for Ub bq}
The multiplication gives an isomorphism $\Ub_{\bq}^{+} \otimes \Ub_{\bq}^0 \otimes \Ub_{\bq}^{-} \overset{\sim}{\longrightarrow} \Ub_{\bq}$. \qed
\end{proposition}

We see next that the 0-part of $\Ub_{\bq}$ can be generated by finitely many quantum brackets, namely those associated with Cartan roots and their heights. 

\begin{lemma}\label{lem:powers-Kbeta-Lbeta}
The subalgebra $\Ub_{\bq}^0$ is generated as an algebra by $K_i$, $L_i$ and $\bra{\beta}{N_{\beta}}$,  $i\in\I$ and $\beta \in \fO^{\bq}_+$. Moreover, we have that
\begin{align}\label{eq:Uq-specialization-relation-grouplike}
K_{\beta}^{N_{\beta}}=L_{\beta}^{N_{\beta}}=1, \quad \bra{\beta}{s}\in \uf_{\bq}^0, \quad \bra{\beta}{mN_{\beta}} = \frac{1}{m!} \prod_{j=0}^{m-1}\left(\bra{\beta}{N_{\beta}}+j(\bq_{\beta}-1)^{N_{\beta}}\right)
\end{align}
for all $\beta\in\fO_+^\bq$, $0\le s<N_\beta$ and $m\in\N$.
\end{lemma}

\begin{proof}
By definition, $\bra{\beta}{0}=1$, and for $s>0$,
$(s)_{\bqf_{\beta}}^!\bra{\beta}{s}=   \prod_{j=0}^{s-1} \left(\bqf_{\beta}^{j}K_{\beta} -L_{\beta}^{-1}\right)$. 
Then, after specialization, we get $0=\prod_{j=0}^{N_\beta-1}\left(\bq_{\beta}^{j}K_{\beta} -L_{\beta}^{-1}\right)$ when $s=N_\beta$ since $\bq_\beta$ is a root of unity of order $N_{\beta}$, and hence $K_{\beta}^{N_\beta}L_{\beta}^{N_\beta}=1$. As  $K_\beta^{N_\beta}=1$, we deduce that $L_{\beta}^{N_\beta}=1$, and the first two relations of \eqref{eq:Uq-specialization-relation-grouplike} hold. If $0<s<N_\beta$, then $(s)_{\bq_{\beta}}^!\neq0$, and hence
$\bra{\beta}{s}=\frac{1}{(s)_{\bq_{\beta}}^!} \prod_{j=0}^{s-1}(\bq_{\beta}^{j}K_{\beta} -L_{\beta}^{-1})\in \uf_{\bq}^0$.

Let $U$ be the subalgebra generated by $\{K_i^{\pm1}, L_i^{\pm 1} | i\in\I\}\cup \left\{ \bra{\beta}{N_{\beta}} | \beta \in \fO^{\bq}_+ \right\}$. We show next that $\Ub_\bq^0=U$.
We have to prove that $\bra{\beta}{n}\in U$ for all $\beta \in \fO^{\bq}_+$ and all $n\in\N_0$. We fix $\beta$ and argue by induction on $n$. If $n<N_\beta$, then 
$\bra{\beta}{n}\in U$ by the above, and it also holds when $n=N_{\beta}$ by definition of $U$.
Let $n>N_{\beta}$ and assume that $\bra{\beta}{n'}\in U$ for all $n'<n$.
 We split the proof of the inductive step into two cases. First, if $N_{\beta}$ does not divide $n$, then, by specializing \eqref{eq:g8} with $M=n-1$ and $N=1$ we get
\begin{align*}
(n)_{\bq_{\beta}}\bra{\beta}{n}=\bra{\beta}{n-1}\bra{\beta}{1}
+(\bq_{\beta}-1) (n-1)_{\bq_{\beta}} K_{\beta} \bra{\beta}{n-1}.
\end{align*}
The right hand side belongs to $U$ by inductive hypothesis. As $(n)_{\bq_{\beta}}\ne0$ we conclude that $\bra{\beta}{n}\in U$. The second case is when $n=N_{\beta} m$ for some $m>1$. Then, by \eqref{eq:g8} with $M=(m-1)N_{\beta}$ and $N=N_{\beta}$, we have that
\begin{align*}
\binom{mN_{\beta}}{N_{\beta}}_{\bqf_{\beta}}  \bra{\beta}{mN_{\beta}}=\sum_{j=0}^{N_{\beta}} (\bq_{\beta}-1)^{j} \binom{(m-1)N_{\beta}+j-1}{j}_{\bq_{\beta}} K_{\beta}^{j} \bra{\beta}{(m-1)N_{\beta}}\bra{\beta}{N_{\beta}-j}.
\end{align*}
As $\binom{mN_{\beta}}{N_{\beta}}_{\bq_{\beta}}=m\neq0$, $\binom{mN_\beta-1}{N_\beta}_{\bq_\beta}=m-1$, and $\binom{(m-1)N_{\beta}+j-1}{j}_{\bq_{\beta}}=0$ for $j\neq0, N_{\beta}$, 
\begin{align*}
m\bra{\beta}{mN_{\beta}}= \bra{\beta}{(m-1)N_{\beta}}\left(\bra{\beta}{N_{\beta}}+(m-1)(\bq_{\beta}-1)^{N_{\beta}}  \right).
\end{align*}
Again, the right hand side belongs to $U$, and therefore $\bra{\beta}{n}\in U$. This completes the induction. Moreover, if we assume that $\bra{\beta}{(m-1)N_{\beta}}= \frac{1}{(m-1)!} \prod_{t=0}^{m-2}\left(\bra{\beta}{N_{\beta}}+t(\bq_{\beta}-1)^{N_{\beta}}\right)$, the above identity implies that
\begin{align*}
\bra{\beta}{mN_{\beta}}= \frac{1}{m} \bra{\beta}{(m-1)N_{\beta}}\bra{\beta}{N_{\beta}}=\frac{1}{m!} \prod_{j=0}^{m-1}\left(\bra{\beta}{N_{\beta}}+j(\bq_{\beta}-1)^{N_{\beta}}\right),
\end{align*}
that is, the last relation of \eqref{eq:Uq-specialization-relation-grouplike} holds.
\end{proof}

Concerning the divided powers, the specialization yields the immediate identities
\begin{align}\label{eq:Uq-specialization-relation-PRV=0}
& E_{\beta}^{N_{\beta}}=F_{\beta}^{N_{\beta}}=0
\end{align}
for all $\beta \in \fO^{\bq}_+$, since $E_{\beta}^{N_\beta}=(N_\beta)_{\bqf_{\beta}}^! E_{\beta}^{(N_\beta)}$ and $F_{\beta}^{N_\beta}=(N_\beta)_{\bqf_{\beta}}^! F_{\beta}^{(N_\beta)}$ by definition. Likewise, the positive and negative parts can be generated by the divided powers of the indecomposable Cartan roots at their heights. The key fact is that there are $\Bbbk$-algebra anti-isomorphisms
\begin{align}\label{eq:isom-specialization-positive,negative}
\phi^-:\Ub_\bq^{-}\longrightarrow \luq\quad\mbox{and}\quad
\phi^+:(\Ub_\bq^{+})^{\underline{\opp}}\longrightarrow \luqm
\end{align}
which identify the divided powers with the appropriate elements in the dual bases of the PBW bases of $\wtoba_{\bq}$ and $\wtoba_{\bq^{(-1)}}$, recall \eqref{eq:fdual-definition}, \eqref{eq:PBW-basis-Lusztig-alg} and Proposition \ref{prop:dual-A-bases}. Indeed, \cite[Proposition 5.9]{AAY} constructs such isomorphisms for $U_{\bqf,\mathsf{A}}^{\pm}\otimes_\mathsf{A}\Bbbk$ where $\mathsf{A}=\Acf\otimes_\Z\Cset$ and hence the same construction works here.
\begin{lemma}\label{lem:Uq-specialization-relations-bases}
The subalgebra $\Ub_{\bq}^+$ is generated as an algebra by $E_i$, $i\in\I$, and $E_\beta^{(N_\beta)}$, $\wbeta \in \varPi^{\bq}$, and the set
\begin{align}\label{eq:Uqpos-special-basis}
\left\{ E_{\beta_{1}}^{(m_{1})} \cdots E_{\beta_\ell}^{(m_\ell)} : 0\le m_j< \widetilde{N}_{\beta_j},\,1\leq j\leq\ell \right\}
\end{align}
is a basis. 
Analogously, $\Ub_{\bq}^-$ is generated as an algebra by $F_i$, $i\in\I$, and $F_\beta^{(N_\beta)}$, $\wbeta \in \varPi^{\bq}$, and the set
\begin{align}\label{eq:Uqneg-special-basis}
\left\{ F_{\beta_1}^{(m_1)} \cdots F_{\beta_\ell}^{(m_\ell)} : 0\le m_j< \widetilde{N}_{\beta_j},\,1\leq j\leq\ell \right\}.
\end{align}
is a basis.
\end{lemma}

\pf
The fact that the set is a basis of $\Ub_{\bq}^\pm$ follows from the isomorphism \eqref{eq:isom-specialization-positive,negative}. This isomorphism together with \cite[Remark 4.10]{AAR-MRL} and \cite[Lemma 3.1 (2)]{AC} imply that $E_i$ and $E_\beta^{(N_\beta)}$ (resp. $F_i$ and $F_\beta^{(N_\beta)}$), $i\in\I$ and $\wbeta \in \varPi^{\bq}$, generate $\Ub_\bq^+$ (resp. $\Ub_\bq^-$). 
\epf

We see below that $\uf_\bq$ is a small quantum group in the sense of \cite[Definition 4.1]{vay-small}. By construction, $\uf_\bq$ is pointed and its group of group-like elements is a quotient of $\Z^{2\theta}$, let $N_{\bq}\subseteq \Z^{2\theta}$ be the kernel. We set
\begin{align}
\Gamma_\bq=G(\uf_\bq)\simeq\Z^{2\theta}/N_\bq
\end{align}
We notice that $K_{\beta}^{N_{\beta}}, L_{\beta}^{N_{\beta}}\in N_\bq$ for all $\beta \in \fO_+^\bq$ by Lemma \ref{lem:powers-Kbeta-Lbeta}.

\begin{proposition}
The subalgebra $\uf^{\pm}_\bq$ is isomorphic to the Nichols algebra $\toba_{\bq^{(\pm1)}}$ and $\uf_{\bq}$ is isomorphic to the quotient of the Drinfeld double $D_\bq$ by $g-1$, $g\in N_{\bq}$. In particular,
\begin{align*}
\left\{ E_{\beta_1}^{m_1} \cdots E_{\beta_\ell}^{m_\ell} \mid 0\le m_j< N_{\beta_j},\,1\leq j\leq\ell \right\}\quad\mbox{and}\\
\left\{ F_{\beta_1}^{m_1} \cdots F_{\beta_\ell}^{m_\ell} \mid 0\le m_j< N_{\beta_j},\,1\leq j\leq\ell \right\}
\end{align*}
are bases of $\uf_\bq^+$ and $\uf_\bq^-$, respectively.
\end{proposition}

\begin{proof}
Clearly, we have an epimorphism of $U_\bq$ onto $\uf_\bq$ and hence $\uf_\bq^{\pm}$ is a pre-Nichols algebra of $\bq^{(\pm1)}$. By Lemma \ref{lem:Uq-specialization-relations-bases}, the sets in the statement are bases of $\uf_\bq^{+}$ and $\uf_\bq^{-}$, which have the same cardinality as the PBW bases of the Nichols algebras $\toba_\bq$ and $\toba_{\bq^{(-1)}}$. Therefore $\uf^{\pm}_\bq\simeq\toba_{\bq^{(\pm1)}}$ and $\uf_\bq\simeq D_\bq/(g-1\mid g\in N_\bq)$.
\end{proof}

\subsection{The quantum Frobenius morphism}\label{subsec:extension}

Our next objective is to prove that $\Ub_\bq$ fits into an extension of Hopf algebras between the small quantum group and the enveloping algebra of the semisimple Lie algebra $\g_{\bq}$, recall Subsection \ref{subsec:cartan-roots}. We start with some technical results needed to show that $\uf_{\bq}$ is a normal Hopf subalgebra. 

To simplify the notation, we set 
\begin{align}
\hbt_{\wbeta}= \bra{\beta}{N_{\beta}},\,  \beta \in \fO^{\bq}_+.
\end{align}
The images of these elements will generate the Cartan subalgebra of $\g_{\bq}$. Below, $\Ss$ stands for the antipode of $\Ub_{\bq}$, as usual.

\begin{lemma}\label{lem:adjoint-action-bracket}
For each $\beta \in \fO^{\bq}_+$, 
\begin{align}\label{eq:coproduct-bracket-Nbeta}
\Delta \left(\hbt_{\wbeta} \right) -\hbt_{\wbeta} \otimes 1 - 1\otimes \hbt_{\wbeta} \in \uf_{\bq}^{0}\otimes \uf_{\bq}^{0}.
\end{align}
Therefore, $\Ss\left(\hbt_{\wbeta}\right) +\hbt_{\wbeta}\in \uf_{\bq}^{0}$ and then
\begin{align}\label{eq:adj-action-bracket-Nbeta}
\left(\ad \hbt_{\wbeta} \right)(x) -\left[\hbt_{\wbeta},x\right] \in \uf_{\bq}
\end{align}
for all $x\in\uf_{\bq}$.
\end{lemma}

\pf
By \eqref{eq:coproduct-bracket-int-form} and Lemma \ref{lem:powers-Kbeta-Lbeta},
\begin{align*}
\Delta \left(\hbt_{\wbeta} \right) -\hbt_{\wbeta} \otimes 1 - 1\otimes \hbt_{\wbeta} &=\Delta\bra{\beta}{N_{\beta}}-\bra{\beta}{N_{\beta}}\otimes 1 - 1\otimes \bra{\beta}{N_{\beta}}
\\
&=\sum_{s=1}^{N_{\beta}-1} L_\beta^{-s}\bra{\beta}{N_{\beta}-s}\otimes K_\beta^{N_{\beta}-s}\bra{\beta}{s}\in \uf_{\bq}^{0}\otimes \uf_{\bq}^{0}.
\end{align*}

As $\epsilon\left(\bra{\beta}{s}\right)=0$ for all $s>0$, the previous equality implies that
\begin{align*}
0 &= m(\Ss\otimes\id)\Delta \left(\hbt_{\wbeta} \right)= \Ss\left(\hbt_{\wbeta}\right)+\hbt_{\wbeta}
+\sum_{s=1}^{N_{\beta}-1} \Ss\left(\bra{\beta}{N_{\beta}-s}\right) L_\beta^{s}K_\beta^{N_{\beta}-s}\bra{\beta}{s},
\end{align*}
and the second statement follows, again by Lemma \ref{lem:powers-Kbeta-Lbeta}. For the last statement,
\begin{align*}
\left(\ad \hbt_{\wbeta} \right)(x) &=
\hbt_{\wbeta} x + x \Ss\left(\hbt_{\wbeta}\right)
+\sum_{s=1}^{N_{\beta}-1} L_\beta^{-s}\bra{\beta}{N_{\beta}-s}x K_\beta^{s}\Ss\left(\bra{\beta}{s}\right),
\end{align*}
so \eqref{eq:adj-action-bracket-Nbeta} follows from Lemma \ref{lem:powers-Kbeta-Lbeta}, the previous characterization of $\Ss\left(\hbt_{\wbeta}\right)$ and the fact that $\uf_{\bq}$ is a Hopf subalgebra.
\epf

The previous result allows us to characterize the coradical of $\Ub_{\bq}$ as follows. 

\begin{lemma}\label{lem:Uq-pointed}
The Hopf algebra $\Ub_{\bq}$ is pointed and $G(\Ub_{\bq})=\Gamma_{\bq}$.
\end{lemma}
\begin{proof}
Let $\Sigma_n$, $n\in\N_0$, be the direct sum of all the $\Z$-components of $\Ub_{\bq}$ of degree $k$ with $|k|\le n$. This is a coalgebra filtration and hence the coradical of $\Ub_{\bq}$ is contained in $\Sigma_0=\Ub_{\bq}^0$. As $K_i$ and $L_i$ are group-like elements and $\Ub_{\bq}^0$ is commutative, the quotient $H$ of $\Ub_{\bq}^0$ by the ideal generated by $K_i-1$, $L_i-1$, $1\le i\le \theta$, is a Hopf algebra generated by the images of $\hbt_{\wbeta}$, $\beta \in \fO^{\bq}_+$, which are primitive in $H$, by Lemma \ref{lem:powers-Kbeta-Lbeta} and \eqref{eq:coproduct-bracket-Nbeta}. Hence $H$ is connected, which implies that $\Ub_{\bq}^0$ is pointed with coradical generated by $K_i$, $L_i$, $1\le i\le \theta$. From here, $G(\Ub_{\bq})=\Gamma_{\bq}$.
\end{proof}

\begin{lemma}\label{lem:adjoint-action-EbetaNbeta}
For each $\beta\in\fO^{\bq}_+$,
\begin{align}\label{eq:coproduct-EbetaNbeta}
\Delta(E_{\beta}^{(N_{\beta})}) -E_{\beta}^{(N_{\beta})}\otimes 1 - 1\otimes E_{\beta}^{(N_{\beta})} \in \uf_{\bq}^{\ge}\otimes \uf_{\bq}^{+},
\\ \label{eq:coproduct-FbetaNbeta}
\Delta(F_{\beta}^{(N_{\beta})}) -F_{\beta}^{(N_{\beta})}\otimes 1 - 1\otimes F_{\beta}^{(N_{\beta})} \in \uf_{\bq}^{-}\otimes \uf_{\bq}^{\le}.
\end{align}
Therefore, $\Ss(E_{\beta}^{(N_{\beta})})+E_{\beta}^{(N_{\beta})}\in \uf_{\bq}^{\ge}$, $\Ss(F_{\beta}^{(N_{\beta})})+F_{\beta}^{(N_{\beta})}\in \uf_{\bq}^{\le}$ and then
\begin{align*}
\left(\ad E_{\beta}^{(N_{\beta})}\right)(x) -\left[E_{\beta}^{(N_{\beta})},x\right] \in \uf_{\bq}\quad\mbox{and}\quad
\left(\ad F_{\beta}^{(N_{\beta})}\right)(x) -\left[F_{\beta}^{(N_{\beta})},x\right]  \in \uf_{\bq}
\end{align*}
for all $x\in\uf_{\bq}$.
\end{lemma}
\pf
We start by proving that
\begin{align}\label{eq:coproduct-braided-EbetaNbeta}
\uDe(E_{\beta}^{(N_{\beta})}) -E_{\beta}^{(N_{\beta})}\otimes 1 - 1\otimes E_{\beta}^{(N_{\beta})} &\in \uf_{\bq}^{+}\otimes \uf_{\bq}^{+} 
\end{align}
for all $\beta\in\fO^{\bq}_+$.
Evoking the skew-pairing as in the proof of Proposition \ref{prop:U-res-+-Hopf-subalg}, we have to prove that if the product of two elements of the PBW basis of $\wtoba_\bq$, say
\begin{align*}
\mathbf{b}=F_{\beta_1}^{m_1} \cdots F_{\beta_\ell}^{m_{\ell}} \ne 1\quad\mbox{and}\quad
\mathbf{c}=F_{\beta_1}^{n_1} \cdots F_{\beta_\ell}^{n_{\ell}} \ne 1
\end{align*}
with $0\le n_k<\widetilde{N}_k$, has $E_{\beta}^{N_{\beta}}$ with non-zero coefficient (when we write it as a linear combination of the elements of the PBW), then $m_k,n_k<N_{\beta_k}$ for all $\beta_k\in\fO^{\bq}_+$. As $\wtoba_{\bq}$ is $\N_0^{\theta}$-graded, we may assume that $N_{\beta}\beta=\sum_j (m_j+n_j)\beta_j$.
Suppose on the contrary that there exists $\beta_k\in\fO^{\bq}_+\setminus\{\beta\}$ such that $\max \{m_k,n_k\}\ge N_{\beta_k}$. If $m_k\ge N_{\beta_k}$, then set
\begin{align*}
\mathbf{b}'=F_{\beta_1}^{m_1} \cdots  F_{\beta_{k-1}}^{m_{k-1}}F_{\beta_k}^{m_k-N_{\beta_k}}F_{\beta_{k+1}}^{m_{k+1}} \cdots F_{\beta_\ell}^{m_{\ell}},
\end{align*}
which is another term of the PBW basis, and $\mathbf{b}=F_{\beta_k}^{N_{\beta_k}}\mathbf{b}'$ as $F_{\beta_k}^{N_{\beta_k}}$ is central by Remark \ref{rem:Zq central}. Let $a_{r_1,\cdots,r_{\ell}}\in\Bbbk$ be such that
\begin{align*}
\mathbf{b}'\mathbf{c}=\sum_{r_1,\cdots,r_{\ell}} a_{r_1,\cdots,r_{\ell}} \, F_{\beta_1}^{r_1} \cdots F_{\beta_\ell}^{r_{\ell}},
\end{align*}
where the sum is over those $r_1,\cdots,r_{\ell}$ with $\sum_j r_j\beta_j=N_{\beta}\beta-N_{\beta_k}\beta_k$. Then
\begin{align*}
\mathbf{b}\mathbf{c} =F_{\beta_k}^{N_{\beta_k}}\mathbf{b}'\mathbf{c}=\sum_{r_1,\cdots,r_{\ell}} a_{r_1,\cdots,r_{\ell}} \, 
F_{\beta_1}^{r_1} \dots F_{\beta_{k-1}}^{r_{k-1}}F_{\beta_k}^{r_k+N_{\beta_k}}F_{\beta_{k+1}}^{r_{k+1}} \cdots F_{\beta_\ell}^{r_{\ell}},
\end{align*}
and hence $F_{\beta}^{N_{\beta}}$ does not appear in the expression of $\mathbf{b}\mathbf{c}$ in terms of the PBW basis.
Otherwise, $n_k\ge N_{\beta_k}$ and the proof is analogous. In any case we get a contradiction, and hence \eqref{eq:coproduct-braided-EbetaNbeta} holds. 

Now, recall that the coproduct for the bosonization is given by
\begin{align*}
\Delta(E_{\beta}^{(N_{\beta})})= \left(E_{\beta}^{(N_{\beta})}\right)^{(1)} \left(\left(E_{\beta}^{(N_{\beta})}\right)^{(2)}\right)_{(-1)} \otimes \left(\left(E_{\beta}^{(N_{\beta})}\right)^{(2)}\right)_{(0)}
\end{align*}
where $x_{(-1)}\otimes x_{(0)} \in \Ub_{\bq}^{0}\otimes \Ub_{\bq}^{+}$ denotes the coaction. Since $x_{(-1)}\otimes x_{(0)} \in \uf_{\bq}^{0}\otimes \uf_{\bq}^{+}$ for all $x\in\uf_{\bq}^+$, and 
\begin{align*}
\left(E_{\beta}^{(N_{\beta})}\right)_{(-1)} \otimes \left( E_{\beta}^{(N_{\beta})}\right)_{(0)} &= K_{\beta}^{N_{\beta}} \otimes E_{\beta}^{(N_{\beta})} = 1 \otimes E_{\beta}^{(N_{\beta})}
\end{align*}
\eqref{eq:coproduct-EbetaNbeta} follows from \eqref{eq:coproduct-braided-EbetaNbeta}.

The statements about $\Ss$ and the adjoint action of $E_{\beta}^{(N_{\beta})}$ follow as in Lemma \ref{lem:adjoint-action-bracket}. The three statements for $F_{\beta}^{(N_{\beta})}$ follow by the same steps, taking into account that the bosonization decomposition of $\Ub_{\bq}^{\le 0}$ is on the right.
\epf

\begin{proposition}\label{prop:uq-normal-in-Uq}
The Hopf subalgebra $\uf_{\bq}$ is normal.
\end{proposition}
\pf
We have to check that $\ad y(\uf_{\bq})\subseteq \uf_{\bq}$ for all $y\in \Ub_{\bq}$. As $\ad$ is an action, it suffices to check it on a set of generators of $\Ub_{\bq}$, and as $\uf_{\bq}$ is a Hopf subalgebra, we have to check it for those generators not in $\uf_{\bq}$. All in all, Lemmas \ref{lem:powers-Kbeta-Lbeta} and \ref{lem:Uq-specialization-relations-bases} together with previous analysis reduces the question to check the following statement:
\begin{align}\label{eq:uq-normal-generators}
\ad y(\uf_{\bq})&\subseteq \uf_{\bq} & \text{for all }&y\in \left\{ \hbt_{\wbeta}, E_{\beta}^{(N_{\beta})}, F_{\beta}^{(N_{\beta})} \vert \beta \in \fO^{\bq}_+ \right\}.
\end{align}
We will prove this case by case. For $y=\hbt_{\wbeta}$ we have that $\ad \hbt_{\wbeta}(\uf_{\bq})\subseteq \uf_{\bq}$ if and only if $\left[\hbt_{\wbeta},\uf_{\bq}\right]\subseteq \uf_{\bq}$ by Lemma \ref{lem:adjoint-action-bracket}. As $\left[\hbt_{\wbeta},-\right]$ is a derivation, it is enough to prove that $\left[\hbt_{\wbeta},x\right]\in\uf_{\bq}$ for $x\in \{E_i,K_i,F_i\vert i\in\I\}$. Clearly, $\left[\hbt_{\wbeta},K_i\right]=0$ for all $i\in\I$. Also, using \eqref{eq:bracket-commutes-Fbeta-Ebeta}, \eqref{eq:condition-cartan-roots}, \eqref{eq:bracket-n<0} and \eqref{eq:bracket-nge0}, we get
\begin{align*}
\hbt_{\wbeta} F_i&= \bra{\beta}{N_{\beta}}F_{i}= \bq_{\alpha_i\beta}^{N_{\beta}}  F_{i}\bracket{\beta}{a_{\beta\alpha_i}}{N_{\beta}}
= F_{i}\bracket{\beta}{a_{\beta\alpha_i}}{N_{\beta}}
\\ 
&= F_{i}\hbt_{\wbeta}+\sum_{j=0}^{N_{\beta}-1} (\bq_{\beta}-1)^{N_{\beta}-j} \binom{-a_{\beta\alpha_i}+N_{\beta}-j-1}{N_{\beta}-j}_{\bq_{\beta}} F_iK_{\beta}^{-j} \bra{\beta}{j}\quad\mbox{and}
\\
\hbt_{\wbeta}E_i &= \bra{\beta}{N_{\beta}}E_{i} = \bq_{\alpha_i\beta}^{-N_{\beta}} E_{i}\bracket{\beta}{-a_{\beta\alpha_i}}{N_{\beta}}
=E_{i}\bracket{\beta}{-a_{\beta\alpha_i}}{N_{\beta}}
\\
&= E_{i}\hbt_{\wbeta}+\sum_{j=0}^{N_{\beta}-1} \bq_{\beta}^{\binom{N_{\beta}-j}{2}-a_{\beta\alpha_i}j}(1-\bq_{\beta})^{N_{\beta}-j} \binom{-a_{\beta\alpha_i}}{N_{\beta}-j}_{\bq_{\beta}} E_iK_{\beta}^{-j} \bra{\beta}{j}.
\end{align*}
Therefore $\left[\hbt_{\wbeta},F_i\right], \left[\hbt_{\wbeta},E_i\right]\in \uf_{\bq}$  by Lemma \ref{lem:powers-Kbeta-Lbeta}.

Similarly, if $y=E_{\beta}^{(N_{\beta})}$, it is enough to check that $\left[E_{\beta}^{(N_{\beta})},\uf_{\bq}\right]\subseteq \uf_{\bq}$ by Lemma \ref{lem:adjoint-action-EbetaNbeta}, which in turn reduces to prove that $\left[E_{\beta}^{(N_{\beta})},x\right]\in\uf_{\bq}$ for all $x\in \{E_i,K_i,F_i\vert i\in\I\}$. The case $x=E_i$ follows by \cite[Proposition 3.2]{AAR-Belgium}. For $x=K_i$,
\begin{align*}
\left[E_{\beta}^{(N_{\beta})},K_i\right] &= (1-\bq_{\alpha_i,\beta}^{N_{\beta}}) E_{\beta}^{(N_{\beta})}K_i
\overset{\eqref{eq:condition-cartan-roots}}{=} 0.
\end{align*}
For $x=F_i$, $\left[E_{\beta}^{(N_{\beta})},F_i\right]=(\id\otimes\langle F_i,-\rangle)\uDe(E_{\beta}^{(N_{\beta})})$: This identity holds in $U_\bqf=D_\bqf$ \cite[Lemma 4.20]{H-Lusztig}, 
then keeps in $\Ub_{\bqf,\Acf}$, and then to the evaluation. Then $\left[E_{\beta}^{(N_{\beta})},F_i\right]\in \uf_{\bq}^+$ by \eqref{eq:coproduct-braided-EbetaNbeta}. 

The proof if $y=F_{\beta}^{(N_{\beta})}$ is analogous to that for the previous case.
\epf

Now, since $\uf_{\bq}$ is a normal Hopf subalgebra, we can consider the quotient Hopf algebra
$\Ub_{\bq}/\langle (\uf_{\bq})^+ \rangle$. Let $\pi:\Ub_{\bq}\to \Ub_{\bq}/\langle (\uf_{\bq})^+ \rangle$ be the canonical projection and set
\begin{align}
\hb_{\wbeta} &=\tfrac{1}{(\bq_\beta-1)^{N_\beta}} \pi \left(\hbt_{\wbeta}\right), & \eb_{\wbeta} &= \tfrac{1}{(\bq_\beta-1)^{N_\beta}}\pi \left(E_{\beta}^{(N_{\beta})}\right), &
\fb_{\wbeta} &= \pi \left( F_{\beta}^{(N_{\beta})}\right), &  \beta & \in \fO^{\bq}_+.
\end{align}
These (normalized) elements generate $\Ub_{\bq}/\langle (\uf_{\bq})^+ \rangle$ as an algebra, are primitive by Lemmas \ref{lem:adjoint-action-bracket} and \ref{lem:adjoint-action-EbetaNbeta}, and will give rise to Chevalley generators.

\begin{proposition}\label{prop:enveloping de gq}
The quotient Hopf algebra
$\Ub_{\bq}/\langle (\uf_{\bq})^+ \rangle$ is isomorphic to the enveloping algebra $U(\g_{\bq})$ of the semisimple Lie algebra $\g_{\bq}$, with the elements $\eb_{\wbeta}$, $\fb_{\wbeta}$, $\hb_{\wbeta}$, $\beta\in \fO^{\bq}_+$, providing a realization of $\g_{\bq}$.
\end{proposition}

\pf
Let us denote by $H$ the quotient $\Ub_{\bq}/\langle (\uf_{\bq})^+ \rangle$. We observe that, as the ideal $\langle (\uf_{\bq})^+ \rangle$ is generated by the homogeneous elements
$E_i$, $F_i$, $K_i-1$, $L_i-1$, $H$ is a $\Z^{\I}$-graded Hopf algebra, with $\Z^{\I}$-degree induced by $\pi$. Let us denote by $\g$ the subspace of $\cP(H)$, the space of primitive elements of $H$, spanned by the elements $\eb_{\wbeta}$, $\fb_{\wbeta}$, $\hb_{\wbeta}$, $\beta\in \fO^{\bq}_+$. Since $H$ is generated by $\g$ as an algebra, \cite[Lemma 5.5.1]{Mont-book} applies with $A_0=\Bbbk 1$ and $A_1=\g+A_0$: $A_n=\left(A_1\right)^n$ gives place to a coalgebra filtration and then $H$ is connected, i.e. its coradical is $H_0=\Bbbk 1$. Therefore, $H=U(\cP(H))$, see e.g. \cite[Theorem 5.6.5]{Mont-book}. Once again, as $H$ is generated by $\g$, we have $\g=\cP(H)$, and hence $H=U(\g)$.

Thus, it remains to prove that $\g\simeq \g_{\bq}$ as Lie algebras. Notice that $\g$ is $\Z^{\I}$-graded, with
\begin{align*}
\deg \hb_{\wbeta} &=0, & \deg \eb_{\wbeta}&= \wbeta, & \deg \fb_{\wbeta}&= -\wbeta, & \beta &\in \fO^{\bq}_+,
\end{align*}
and $\g$ has a triangular decomposition $\g=\n_+\oplus\h\oplus\n_-$, where $\n_+$, respectively, $\h$, $\n_-$, is the subspace spanned by 
$(\eb_{\wbeta})_{\beta \in \fO^{\bq}_+}$, respectively $(\hb_{\wbeta})_{\beta \in \fO^{\bq}_+}$, $(\fb_{\wbeta})_{\beta \in \fO^{\bq}_+}$.
We will check that each defining relation of $\g_{\bq}$ holds in $\g$. First, $[\hb_{\wbeta},\hb_{\wgamma}]=0$ for all $\beta,\gamma \in \fO^{\bq}_+$ since $\h$ is the image of the commutative algebra $\Ub_{\bq}^0$.

For elements in the positive part of $\g_{\bq}$, there is nothing to check as $\n_+\simeq \n_+^{\bq}$ 
by \cite{AAR3}. In particular, $\n_+$ is generated by $(\eb_{\wbeta})_{\wbeta \in \varPi^{\bq}}$ as a Lie algebra, and these elements satisfy Serre relations: if $A_{\bq}=(a_{\wbeta\wgamma})_{\wbeta,\wgamma\in \varPi^{\bq}}$ is the Cartan matrix of $\g_{\bq}$, then $(\ad \eb_{\wbeta})^{a_{\wbeta\wgamma}}\eb_{\wgamma}=0$ for all $\wbeta,\wgamma\in \varPi^{\bq}$. Dually, $\n_- \simeq \n_-^{\bq}$ is generated by $(\fb_{\wbeta})_{\wbeta \in \varPi^{\bq}}$ as Lie algebras, and $(\ad \fb_{\wbeta})^{a_{\wbeta\wgamma}}\fb_{\wgamma}=0$ for all $\wbeta,\wgamma\in \varPi^{\bq}$.

Let $\wbeta\ne\wgamma \in \varPi^{\bq}$. Then $\wbeta-\wgamma\notin\wfO^{\bq}_+\cup(-\wfO^{\bq}_+)\cup\{0\}$, so $\g_{\wbeta-\wgamma}=0$. As $[\eb_{\wbeta},\fb_{\wgamma}]\in\g_{\wbeta-\wgamma}$, we have that $[\eb_{\wbeta},\fb_{\wgamma}]=0$. Also, by Lemma \ref{lemmma:Ealpha-M-Falpha-N},
\begin{align*}
\eb_{\wbeta}\fb_{\wbeta} &= \tfrac{1}{(\bq_\beta-1)^{N_\beta}} \pi \left(E_{\beta}^{(N_{\beta})}F_{\beta}^{(N_{\beta})}\right) = \sum_{t=0}^{N_{\beta}} \tfrac{1}{(\bq_\beta-1)^{N_\beta}}\pi \left(F_{\beta}^{(N_{\beta}-t)} \bracket{\beta}{2N_{\beta}-t-1}{t}E_{\beta}^{(N_{\beta}-t)}\right)
\\
&= \tfrac{1}{(\bq_\beta-1)^{N_\beta}}\pi \left(F_{\beta}^{(N_{\beta})}E_{\beta}^{(N_{\beta})}\right)+\tfrac{1}{(\bq_\beta-1)^{N_\beta}}\pi\left(\bracket{\beta}{N_{\beta}-1}{N_{\beta}}\right)
=\fb_{\wbeta}\eb_{\wbeta}+\hb_{\wbeta},
\end{align*}
since $E_{\beta}^{(j)}\in(\uf_{\bq})^+$ for all $0<j<N_{\beta}$. Thus $[\eb_{\wbeta},\fb_{\wbeta}]=\hb_{\wbeta}$.

Finally, using \eqref{eq:bracket-commutes-Fbeta-Ebeta}, \eqref{eq:condition-cartan-roots} and \eqref{eq:bracket-nge0}, respectively \eqref{eq:bracket-n<0}, we get
\begin{align*}
\hb_{\wbeta}\fb_{\wgamma} &= \tfrac{1}{(\bq_\beta-1)^{N_\beta}}\pi\left(\bra{\beta}{N_{\beta}}F_{\gamma}^{(N_{\gamma})} \right) 
=\tfrac{1}{(\bq_\beta-1)^{N_\beta}}\pi \left(\bq_{\gamma\beta}^{N_{\beta}N_{\gamma}}  F_{\gamma}^{(N_{\gamma})}\bracket{\beta}{N_{\gamma}a_{\beta\gamma}}{N_{\beta}}\right)
\\
&= \fb_{\wgamma}\pi \left( \sum_{j=0}^{N_\beta}(\bq_\beta-1)^{-j}\binom{-N_\gamma a_{\beta\gamma} +N_\beta-j-1}{N_\beta-j}_{\bq_\beta}K_\beta^{N_\beta-j}\bra{\beta}{j} \right)
\\
&= \fb_{\wgamma}\pi \left(\binom{-N_\beta a_{\wbeta\wgamma}+N_\beta-1}{N_\beta}_{\bq_\beta}K_\beta^{N_\beta}\right)
+\tfrac{1}{(\bq_\beta-1)^{N_\beta}}\fb_{\wgamma}\pi \left( \bq_{\beta}^{-N_{\gamma}a_{\beta\gamma}N_{\beta}}\bra{\beta}{N_{\beta}} \right)
\\
&= -a_{\wbeta\wgamma}\fb_{\wgamma}+\fb_{\wgamma}\hb_{\wbeta},\quad\mbox{and}
\\
\hb_{\wbeta}\fb_{\wbeta} &= \tfrac{1}{(\bq_\beta-1)^{N_\beta}}\pi\left(\bra{\beta}{N_{\beta}}F_{\beta}^{(N_{\beta})} \right) 
=\tfrac{1}{(\bq_\beta-1)^{N_\beta}}\pi \left(\bq_{\beta}^{N_{\beta}^2}  F_{\beta}^{(N_{\beta})}\bracket{\beta}{2N_{\beta}}{N_{\beta}}\right)
\\
&= -2\fb_{\wbeta}+\fb_{\wbeta}\hb_{\wbeta}.
\end{align*}
Thus, $[\hb_{\wbeta},\fb_{\wgamma}]=-a_{\wbeta\wgamma}\fb_{\wgamma}$ for all $\wbeta, \wgamma\in\varPi^{\bq}$. Now, $[\hb_{\wbeta},\eb_{\wgamma}]\in\g_{\wgamma}=\Bbbk\eb_{\wgamma}$ since $\g$ is $\Z^{\I}$-graded, and using the relations already proved we have
\begin{align*}
0 &= [\hb_{\wbeta},\hb_{\wgamma}] = [\hb_{\wbeta},[\eb_{\wgamma},\fb_{\wgamma}]]
= [[\hb_{\wbeta},\eb_{\wgamma}],\fb_{\wgamma}] -a_{\wbeta\wgamma}\hb_{\wgamma}.
\end{align*}
Hence, $[\hb_{\wbeta},\eb_{\wgamma}]=a_{\wbeta\wgamma}\eb_{\wgamma}$.

In conclusion, we have a morphism $\Phi:\g_{\bq}\twoheadrightarrow\g$ of $\Z^{\I}$-graded Lie algebras, which keeps the triangular decomposition, and the restrictions to the positive and the negative parts are isomorphisms. Moreover, this is surjective because $\h$ is spanned by $(\hb_{\wbeta})_{\wbeta \in \varPi^{\bq}}$. Indeed, let $\wbeta\in\wfO^{\bq}_+\setminus\varPi^{\bq}$. Then there exist $\walpha \neq \wgamma\in \wfO^{\bq}_+$ such that $\wbeta=\walpha+\wgamma$, and by \cite[Lemma 3.11]{AAR3}, there exists $c\in\Bbbk^{\times}$ such that $c[\eb_{\walpha},\eb_{\wgamma}]=\eb_{\wbeta}$. Also, $d[\fb_{\walpha},\fb_{\wgamma}]=\fb_{\wbeta}$ for some $d\in\Bbbk^{\times}$. Using Jacobi identity and the relations above,
\begin{align*}
\hb_{\wbeta} &= [\eb_{\wbeta},\fb_{\wbeta}] = d\, \left[ \eb_{\wbeta}, [\fb_{\walpha},\fb_{\wgamma}] \right]
= d\, \left[ \left[\eb_{\wbeta}, \fb_{\walpha}\right],\fb_{\wgamma} \right]+
d\, \left[ \fb_{\walpha}, \left[ \eb_{\wbeta}, \fb_{\wgamma}\right] \right]
\\
&= cd\, \left[ \left[[\eb_{\walpha},\eb_{\wgamma}], \fb_{\walpha}\right],\fb_{\wgamma} \right]+
cd\, \left[ \fb_{\walpha}, \left[ [\eb_{\walpha},\eb_{\wgamma}], \fb_{\wgamma}\right] \right]
\\
&= cd\, \left[ \left[\hb_{\walpha}, \eb_{\wgamma} \right],\fb_{\wgamma} \right]+
cd\, \left[ \fb_{\walpha}, \left[ \eb_{\walpha},\hb_{\wgamma}\right] \right]
= cd \left(a_{\walpha\wgamma}\hb_{\wgamma}+a_{\wgamma\walpha} \hb_{\walpha}\right).
\end{align*}
Recursively we can prove that $\h$ is spanned by $(\hb_{\wbeta})_{\wbeta \in \varPi^{\bq}}$. 

Let us see that $\Phi$ is injective. Let $h\in\ker \Phi$, then $h\in \h_{\bq}$ and hence $\Phi(h)$ acts by zero on every $\eb_{\wbeta}$ as $\Phi$ is a Lie algebra map, but as $A_{\bq}$ is invertible, we have that $h=0$.

Therefore $\Phi$ is an isomorphism of Lie algebras and the statement follows.
\epf

The above results immediately yield the distinctive structural property of the quantum group of divided powers of $\bq$.

\begin{theorem}\label{thm:Uq-extension-Hopf}

There exists an extension of Hopf algebras
\begin{equation*}
\begin{tikzcd}[column sep=large, arrow style=tikz, >=To]
	\uf_{\bq} \arrow[r, hook] & \Ub_{\bq} \arrow[r, two heads, "\operatorname{Fr}"] & U(\g_{\bq}).
\end{tikzcd}
\end{equation*}
\qed
\end{theorem}

In particular, the restriction of this extension to the zero components gives the algebra isomorphism
\begin{align}\label{eq:U0q = Gammaq ot Uhq}
\Ub_{\bq}^0\simeq\Gamma_\bq\,\otimes\,U(\h_\bq)
\end{align}
Also, as a direct consequence, we obtain natural bases for $\Ub_{\bq}^0$ and $\Ub_{\bq}$.

\begin{corollary}\label{thm:Uq-special}
The set $\left\{ g \prod_{\wbeta\in\varPi^{\bq}} \bra{\beta}{N_{\beta}}^{c_{\beta}} : g\in\Gamma_{\bq}, c_{\beta} \in\N_{0} \right\}$ is a basis of $\Ub_{\bq}^0$, and
\begin{align*}
\left\{  F_{\beta_{1}}^{(n_{1})} \dots F_{\beta_{\ell}}^{(n_{\ell})} \,  g \prod_{\wbeta\in\varPi^{\bq}} \bra{\beta}{N_{\beta}}^{c_{\beta}}  \, E_{\beta_{1}}^{(m_1)} \dots E_{\beta_{\ell}}^{(m_{\ell})} : 0\le m_j,n_j< \widetilde{N}_{\beta_j}, g\in\Gamma_{\bq}, c_{\beta} \in\N_{0} \right\}.
\end{align*}
is a basis of $\Ub_{\bq}$.
\end{corollary}

\pf
The first set is a basis of $\Ub_{\bq}^0$ since $\uf_{\bq}^0=\Bbbk\Gamma_{\bq}$, and $\{\prod_{\wbeta\in\varPi^{\bq}} \bra{\beta}{N_{\beta}}^{c_{\beta}} \mid  c_{\beta} \in\N_{0} \}$ is the preimage of the basis $\{\prod_{\wbeta\in\varPi^{\bq}} \hb_{\wbeta}^{c_{\beta}} \mid  c_{\beta} \in\N_{0} \}$ of $U(\h_\bq)$. 
The second set is a basis of $\Ub_\bq$ by Proposition \ref{prop:triangular decomposition for Ub bq} and Lemma \ref{lem:Uq-specialization-relations-bases}.
\epf

\subsection{Lusztig isomorphisms}\label{subsec:specialization-Lusztig-iso}

Next we show the existence of Lusztig isomorphisms for quantum groups of divided powers and then derive some consequences of this fact.

\begin{theorem}\label{thm:Uq-special-Lusztig-isom}
For each $i\in \I$, the Lusztig isomorphism $T_{i}^{\bqf}:\Ub_{\rho_i(\bqf), \Acf}\to \Ub_{\bqf, \Acf}$ induces an isomorphism $T_{i}^{\bq}:\Ub_{\rho_i(\bq)}\to \Ub_{\bq}$.
\end{theorem}
\pf
The isomorphism $T_i^{\bqf}:\Ub_{\rho_i(\bqf),\Acf} \to \Ub_{\bqf,\Acf}$ of $\Acf$-algebras in Theorem \ref{thm:triang-decomp-divided-powers} induces an isomorphism of $\Bbbk$-algebras 
$T_i^{\bq}:\Ub_{\rho_i(\bqf),\Acf}\otimes_{\Acf}\Bbbk \to \Ub_{\bqf,\Acf}\otimes_{\Acf}\Bbbk$. As $T_i^{\bq}(K_\beta)=K_{s_i^{\bq}(\beta)}$ and $N_{\beta}=N_{s_i^{\bq}(\beta)}$ for all $\wbeta\in\varPi^\bq$, and
\begin{align}
\label{eq:si-Pi-q}
s_i^{\bq}(\varPi^{\rho_i(\bq)}) &=\begin{cases}
\left(\varPi^{\bq}\setminus\{\alpha_i\}\right)\cup\{-\alpha_i\} & \text{if }\alpha_i\text{ is Cartan},
\\
\varPi^{\rho_i(\bq)}, & \text{if }\alpha_i\text{ is not Cartan}.
\end{cases}
\end{align}
we have that $T_i^{\bq}\left(\Ib_{\rho_i(\bq)}\right)=\Ib_{\bq}$, so it induces the desired isomorphism $T_{i}^{\bq}:\Ub_{\rho_i(\bq)}\to \Ub_{\bq}$.
\epf

\begin{remark}\label{rem:Uq-special-Lusztig-isom-Nichols}
The restriction of $T_{i}^{\bq}$ to $\uf_{\rho_i (\bq)}$ is the Lusztig isomorphism $T_{i}^{\bq}:\uf_{\rho_i(\bq)}\to \uf_{\bq}$ attached to the Drinfeld doubles of Nichols algebras of $\bq$ and $\rho_i(\bq)$. Also, it restricts to a group isomorphism $T_{i}^{\bq}:\Gamma_{\rho_i (\bq)}\to \Gamma_{\bq}$.
\end{remark}

We use the Lusztig isomorphisms to describe explicitly the group $\Gamma_{\bq}$ under a technical condition related to Remark \ref{rem:ZPi-proper-subgroup}.

\begin{lemma}\label{le:U-0-eval-A-repr}
Let $N=\ord \xi$ and assume that $2$ does not divide $N$ if $\bqf$ is of type $\superd{m}{n}$, $\mathbf{D}(2,1;\alpha)$, $\mathbf{F}(4)$ or $\texttt{br}(2)$, and $3$ does not divide $N$ if $\bqf$ is of type $\texttt{wk}(4)$. Then, given $0\le m_{\beta},n_{\beta}<N_{\beta}$, $\wbeta\in\varPi^\bq$, there exists an algebra map $\chi:\Ub^0_{\bq}\to\Bbbk$ such that 
\begin{align}\label{eq:U-0-eval-character}
\chi(K_{\beta})=\bq_{\beta}^{m_{\beta}}\quad\mbox{and}\quad \chi(L_{\beta})=\bq_{\beta}^{n_{\beta}}.
\end{align}
Therefore, a presentation for the group $\Gamma_{\bq}$ is
\begin{align*}
\Gamma_{\bq} =\left(
K_{i}, L_{i},\, i\in\I\mid K_{\beta}^{N_{\beta}}=1= L_{\beta}^{N_{\beta}},\, \wbeta \in \varPi^{\bq}\right).
\end{align*}

\end{lemma}

\pf
For the first claim, Lusztig isomorphisms allow us to reduce the proof of the existence of such an algebra map for just one braiding on each Weyl equivalence class. We pick the diagrams in \cite[Appendix A]{AAY}. 
We prove first that there exists a $\Cset(\nu)$-algebra map $\chi:U_{\bqf}^0\to\Cset(\nu)$ such that $\chi(K_i),\chi(L_i)\in\Acf$, for all $i\in\I$,
\begin{align*}
\chi(K_{\beta})=\bqf_{\beta}^{\widetilde{m}_{\beta}}\quad\mbox{and}\quad \chi(L_{\beta})=\bqf_{\beta}^{\widetilde{n}_{\beta}}
\end{align*}
for all $\wbeta\in\varPi^\bq$, where $\widetilde{m}_{\beta}$ and  $\widetilde{n}_{\beta}$ are integers satisfying that $N$ divides $\widetilde{m}_{\beta}-m_{\beta}$ and $\widetilde{n}_{\beta}-n_{\beta}$.
As $U_{\bqf}^0=\Cset(\nu)[K_i^{\pm1},L_i^{\pm1}]$, this amounts to choosing suitable values for $c_i\coloneqq\chi(K_i)$ and $c_i'\coloneqq \chi(L_i)$ in $\Cset(\nu)^{\times}$.

If $\bqf$ is of Cartan type, then $\varPi^\bq=\{\alpha_i, \, 1\le i\le \theta\}$ and  it is enough to choose $\widetilde{m}_{\beta}=m_{\beta}$ and $\widetilde{n}_{\beta}=n_{\beta}$. If $\bqf$ is of type $\supera{k}{\theta-k}$, then $\varPi^\bq=\{\alpha_i, \, i\ne k\}$. Hence, we choose $\widetilde{m}_{\beta}=m_{\beta}$ and $\widetilde{n}_{\beta}=n_{\beta}$, and we can choose any value for $c_k$ and $c_k'$.

In all the other cases, $\varPi^\bq=\{\alpha_i, \, i\ne k\}\cup\{\gamma\}$ for suitable $k\in\{1,\cdots,\theta\}$ and $\gamma\in\N_0^{\theta}$ with $k$ in the support of $\gamma$. Then, for $i\ne k$ we set $c_i=\bqf_{ii}^{m_i}$, $c_i'=\bqf_{ii}^{n_i}$ except possibly in type $\superda{\alpha}$ as we will explain later, and it remains to choose $c_k$ and $c_k'$. We do it case by case for $c_k$, the proof for $c_k'$ is analogous.
\begin{description}[leftmargin=*]
\item[$\superb{k}{\theta-k}$] Here, $\gamma=\alpha_k+\cdots+\alpha_{\theta}$, so $\chi(K_{\gamma})=\bqf_{\gamma}^{\widetilde{m}_{\gamma}}$ if and only if
\begin{align*}
c_k &=\bqf_{\gamma}^{\widetilde{m}_{\gamma}}\bqf_{k+1,k+1}^{-m_{k+1}} \cdots \bqf_{\theta\theta}^{-m_{\theta}}
=(-1)^{\widetilde m_\gamma}\nu^{\widetilde{m}_{\gamma}-2m_{k+1}- \cdots -2m_{\theta-1}-m_{\theta}}.
\end{align*}
Thus we set $\widetilde{m}_{\gamma}=m_{\gamma}$ and the corresponding value of $c_k$.

\item[$\superd{k}{\theta-k}$] $\gamma=\alpha_{k-1}+2\alpha_k+\cdots+2\alpha_{\theta-1}+\alpha_{\theta}$, so $\chi(K_{\gamma})=\bqf_{\gamma}^{\widetilde{m}_{\gamma}}$ if and only if
\begin{align*}
c_k^2 &=\bqf_{\gamma}^{\widetilde{m}_{\gamma}}\bqf_{k-1,k-1}^{-m_{k-1}}\bqf_{k+1,k+1}^{-2m_{k+1}} \cdots \bqf_{\theta-1,\theta-1}^{-2m_{\theta-1}} \bqf_{\theta\theta}^{-m_{\theta}}
=\nu^{-\widetilde{m}_{\gamma}+m_{k-1}-2m_{k+1}-\cdots-2m_{\theta-1}-2m_{\theta}}.
\end{align*}
Thus we can choose $\widetilde{m}_{\gamma}\in \{m_{\gamma}, m_{\gamma}+N\}$ such that the exponent of $\nu$ is even since $N$ is odd by assumption, and then 
$c_k=\nu^{\frac{1}{2}(-\widetilde{m}_{\gamma}+m_{k-1}-2m_{k+1}-\cdots-2m_{\theta-1}-2m_{\theta})}$.

\item[$\superda{\alpha}$] $\gamma=\alpha_1+2\alpha_2+\alpha_3$, $k=2$, so $\chi(K_{\gamma})=\bqf_{\gamma}^{\widetilde{m}_{\gamma}}$ if and only if
\begin{align*}
c_2^2 &=\bqf_{\gamma}^{\widetilde{m}_{\gamma}}\bqf_{11}^{-\widetilde m_1}\bqf_{33}^{-\widetilde m_3} =\nu^{-\widetilde{m}_{\gamma}(d_1+d_3)-\widetilde{m}_{1}d_1-\widetilde{m}_{3}d_3}.
\end{align*}
If $m_{1}d_1+m_{3}d_3$ is even, then we choose $\widetilde{m}_i=m_i$ for $i=1,3$, $\widetilde{m}_{\gamma}\in \{m_{\gamma}, m_{\gamma}+N\}$ even. 
If $m_{1}d_1+m_{3}d_3$ is odd, then either $m_1$ and $d_1$ are odd, or $m_3$ and $d_3$ are so: assuming $m_1$, $d_1$ odd, we choose $\widetilde{m}_1=m_1+N$, $\widetilde{m}_3=m_3$, so $\widetilde{m}_{1}d_1+\widetilde{m}_{3}d_3$ becomes even and choose again $\widetilde{m}_{\gamma}\in \{m_{\gamma}, m_{\gamma}+N\}$ even. In any case,
$c_2=\nu^{-\frac{1}{2}(\widetilde{m}_{\gamma}(d_1+d_3)+\widetilde{m}_{1}d_1+\widetilde{m}_{3}d_3)}$.

\item[$\superf$] $\gamma=\alpha_1+2\alpha_2+3\alpha_3+2\alpha_4$, $k=4$, so $\chi(K_{\gamma})=\bqf_{\gamma}^{\widetilde{m}_{\gamma}}$ if and only if
\begin{align*}
c_4^2 &=\bqf_{\gamma}^{\widetilde{m}_{\gamma}}\bqf_{11}^{-m_1}\bqf_{22}^{-2m_2} \bqf_{33}^{-3m_3}
=\nu^{-3\widetilde{m}_{\gamma}-2m_1-4m_2-3m_3}.
\end{align*}
Thus we can choose $\widetilde{m}_{\gamma}\in \{m_{\gamma}, m_{\gamma}+N\}$ such that the exponent of $\nu$ is even by our assumption on $N$, and set 
$c_4=\nu^{-\frac{1}{2}(3\widetilde{m}_{\gamma}+2m_1+4m_2+3m_3)}$.

\item[$\superg$] $\gamma=\alpha_1+2\alpha_2+\alpha_3$, $k=1$, so $\chi(K_{\gamma})=\bqf_{\gamma}^{\widetilde{m}_{\gamma}}$ if and only if
\begin{align*}
c_1 &=\bqf_{\gamma}^{\widetilde{m}_{\gamma}}\bqf_{22}^{-2m_2}\bqf_{33}^{-m_3}
=(-1)^{m_{\gamma}}\nu^{-\widetilde{m}_{\gamma}-2m_2-3m_3}.
\end{align*}
Thus we set $\widetilde m_\gamma=m_\gamma$, $c_1=(-1)^{m_{\gamma}}\nu^{-m_{\gamma}-2m_2-3m_3}$.

\item[$\mathtt{wk}(4)$] $\gamma=\alpha_1+2\alpha_2+3\alpha_3+\alpha_4$, $k=3$, so $\chi(K_{\gamma})=\bqf_{\gamma}^{\widetilde{m}_{\gamma}}$ if and only if
\begin{align*}
c_3^3 &=\bqf_{\gamma}^{\widetilde{m}_{\gamma}}\bqf_{11}^{-m_1}\bqf_{22}^{-2m_2} \bqf_{44}^{-m_4}
=(-1)^{\widetilde{m}_{\gamma}+m_4}\nu^{-\widetilde{m}_{\gamma}+m_4-m_1-2m_2}.
\end{align*}
Thus we choose $\widetilde{m}_{\gamma}\in \{m_{\gamma}, m_{\gamma}+N, m_{\gamma}+2N\}$ such that $3$ divides the exponent of $\nu$, and set 
$c_3=(-1)^{\widetilde{m}_{\gamma}+m_4}\nu^{-\frac{1}{3}(\widetilde{m}_{\gamma}+m_1+2m_2-m_4)}$.

\item[$\mathtt{br}(2)$] Here, $\gamma=2\alpha_1+\alpha_2$, $k=1$, so $\chi(K_{\gamma})=\bqf_{\gamma}^{\widetilde{m}_{\gamma}}$ if and only if
$c_1^2 =\bqf_{\gamma}^{\widetilde{m}_{\gamma}}\bqf_{22}^{-m_2}=\zeta^{\widetilde{m}_{\gamma}}\nu^{-\widetilde{m}_{\gamma}-m_2}$.
Thus we choose $\widetilde{m}_{\gamma}\in \{m_{\gamma}, m_{\gamma}+N\}$ such that $\widetilde{m}_{\gamma}+m_2$ is even and set 
$c_1=\zeta^{2\widetilde{m}_{\gamma}}\nu^{-\frac{1}{2}(\widetilde{m}_{\gamma}+m_2)}$.

\end{description}

Hence the hypothesis of Lemma \ref{le:U-0-A-repr} holds, so $\chi$ restricts to an $\Acf$-algebra map $\chi:\Ub_{\bqf,\Acf}^0\to\Acf$. Therefore we can extend scalars to get an algebra map $\chi:\Ub_{\bq}^0\to\Bbbk$ such that $\chi(K_{\beta})=\bq_{\beta}^{\widetilde{m}_{\beta}}=\bq_{\beta}^{m_{\beta}}$ and $\chi(L_{\beta})=\bq_{\beta}^{\widetilde{n}_{\beta}}=\bq_{\beta}^{n_{\beta}}$  for all $\wbeta\in\Pi^\bq$, as desired. We highlight that, as a consequence, we get $M\coloneqq \Big(\prod_{\wbeta \in \varPi^{\bq}} N_{\beta}\Big)^2$ different characters of $\Gamma_{\bq}$.

In order to prove the isomorphism $\Gamma_{\bq}=\Z^{2\theta}/N_\bq \simeq \Lambda_{\bq}\coloneqq \Z^{2\theta}/(K_{\beta}^{N_{\beta}}, L_{\beta}^{N_{\beta}}| \wbeta \in \varPi^{\bq})$, we recall that $(K_{\beta}^{N_{\beta}}, L_{\beta}^{N_{\beta}}| \wbeta \in \varPi^{\bq})\subseteq N_\bq$ by Lemma \ref{lem:powers-Kbeta-Lbeta}, so there exists an epimorphism $\Lambda_{\bq}\twoheadrightarrow \Gamma_{\bq}$. Let us check, case by case according to Remark \ref{rem:ZPi-proper-subgroup}, that it is indeed an isomorphism by showing that $\Gamma_\bq$ admits as many different characters as $\Lambda_\bq$. In fact,  $\Lambda_{\bq}$ is finite of order $M$ if $\Z^{2\theta}=(K_{\beta}, L_{\beta}| \wbeta \in \varPi^{\bq})$. If $\bq$ is of type $\supera{k}{\theta-k}$, then $\Lambda_{\bq}$ is a product of a finite group of order $M$ with a free abelian group generated by $K_k$ and $L_k$, and we can choose arbitrarily the values on $K_k$ and $L_k$ for each of the $M$ characters constructed above. If $\Z^{2\theta}/(K_{\beta}, L_{\beta}| \wbeta \in \varPi^{\bq})\simeq\Z/2\Z\times\Z/2\Z$, then the order of $\Lambda_{\bq}$ is $4M$. In this case, let $\chi_{K_k},\chi_{L_k}:\Ub_{\bq}\rightarrow\Bbbk$ be the algebra maps such that $\chi_{K_k}(K_k)=\chi_{L_k}(L_k)=-1$ and that are equal to $1$ on the remaining generators. Since $N$ is odd, we obtain $4M$ different characters for $\Gamma_\bq$. Finally, if $\bq$ is of type $\mathtt{wk}(4)$, then $\Z^{2\theta}/(K_{\beta}, L_{\beta}| \wbeta \in \varPi^{\bq})\simeq\Z/3\Z\times\Z/3\Z$ and the order of $\Lambda_{\bq}$ is $9M$. In this case, let $\zeta$ be a third root of unity and $\Acf'=\Acf[\zeta]$. Since $\Ub_{\bq}= \Ub_{\bqf,\Acf}\otimes_{\Acf'}\Acf'\otimes_{\Acf} \Bbbk/\Ib_{\bq}$, we can get algebra maps $\chi_{K_k},\chi_{L_k}:\Ub_{\bq}\rightarrow\Bbbk$ such that $\chi_{K_k}(K_k)=\chi_{L_k}(L_k)=\zeta$ and that are equal to $1$ on the remaining generators. Hence, we obtain $9M$ different characters for $\Gamma_\bq$. This completes the proof.
\epf

We finish this subsection by analyzing the effect of Lusztig isomorphisms on the enveloping algebras. We notice that 
$\g_{\rho_i (\bq)} \simeq \g_{\bq}$, and $s_i^{\bq}$ moves $\varPi^{\rho_i(\bq)}$ according to \eqref{eq:si-Pi-q}.

\begin{theorem}\label{thm:Uq-special-Lusztig-isom-Lie-alg} 
For each $i\in\I$, the isomorphism $T_i^{\bq}$ induces an isomorphism of Hopf algebras $\overline{T}_i^{\bq}: U(\g_{\rho_i(\bq)}) \to U(\g_{\bq})$.
\end{theorem}
\pf
The isomorphism $T_{i}^{\bq}:\Ub_{\rho_i(\bq)}\to \Ub_{\bq}$ from Theorem \ref{thm:Uq-special-Lusztig-isom} restricts to the Lusztig isomorphism 
$T_{i}^{\bq}:\uf_{\rho_i(\bq)}\to \uf_{\bq}$, as noticed in Remark \ref{rem:Uq-special-Lusztig-isom-Nichols}. Looking at the definition, 
$T_i^{\bq}(X)\in(\uf_{\bq})^+$ for all $X\in\{K_i-1,L_i-1,E_i,F_i\}$, so $T_i^{\bq}\left( (\uf_{\rho_i(\bq)})^+\right) \subseteq(\uf_{\bq})^+$. Hence $T_i^{\bq}$ induces an algebra map $\overline{T}_i^{\bq}: U(\g_{\rho_i(\bq)}) \to U(\g_{\bq})$, and analogously we have an algebra map $\overline{T}_i^{\rho_i(\bq)}: U(\g_{\bq}) \to U(\g_{\rho_i(\bq)})$. As $T_i^{\rho_i(\bq)}$ is the inverse of $T_i^{\bq}$ on $\Ub_{\bq}$, the same happens for the induced maps on $U(\g_{\bq})$. It remains to prove that this is a Hopf algebra map; to do so, we will prove that $\overline{T}_i^{\bq}(\g_{\rho_i(\bq)})=\g_{\bq}$.

As in the proof of Theorem \ref{thm:triang-decomp-divided-powers}, we fix $i\in\I$, a reduced expression 
$w_0^{\bqf}=\sigma_{i_1}^{\bqf} \cdots \sigma_{i_{\ell}}$ such that $i_1=i$ and $\bpf\coloneqq\rho_i(\bqf)$. Then $\bp\coloneqq\rho_i(\bq)$ is obtained by evaluation of $\bpf$. Recall that $\sigma_{i_2}^{\bpf} \cdots \sigma_{i_{\ell}}$ can be extended to a reduced expression
$w_0^{\bpf}=\sigma_{i_2}^{\bpf} \cdots \sigma_{i_{\ell}}\sigma_{i_{\ell+1}}$ for some $i_{\ell+1}\in\I$. As there we set
\begin{align*}
\beta_j&:=s_{i_1}^{\bqf} \cdots s_{i_{j-1}}(\alpha_{i_j})\in\varDelta_+^{\bqf}, & 
\gamma_j&:=s_{i_2}^{\bpf} \cdots s_{i_{j}}(\alpha_{i_{j+1}})\in\varDelta_+^{\bpf}, & &i\in\I_{\ell},
\end{align*}
so $\varDelta_+^{\bqf}=\{\beta_j | j\in\I_{\ell}\}$, $\varDelta_+^{\bpf}=\{\gamma_j | j\in\I_{\ell}\}$, $s_i^{\bq}(\gamma_j)=\beta_{j+1}$ for all $i\in\I_{\ell-1}$ and $\gamma_{\ell}=\alpha_i$, and
\begin{align*}
E_{\beta_j} &:= T_{i_1}^{\bqf} \cdots T_{i_{j-1}}(E_{i_j})\in \Ub_{\bqf,\Acf}^{+}, & 
F_{\beta_j} &:= T_{i_1}^{\bqf} \cdots T_{i_{j-1}}(F_{i_j})\in \Ub_{\bqf,\Acf}^{-},
\\
\underline{E}_{\gamma_j} &:= T_{i_2}^{\bpf} \cdots T_{i_{j}}(E_{i_{j+1}})\in \Ub_{\bpf,\Acf}^{+}, &
\underline{F}_{\gamma_j} &:= T_{i_2}^{\bpf} \cdots T_{i_{j}}(F_{i_{j+1}})\in \Ub_{\bpf,\Acf}^{-}.
\end{align*}
By definition, $T_i^{\bqf}(\underline{E}_{\gamma_j})=E_{\beta_{j+1}}$ and $T_i^{\bqf}(\underline{F}_{\gamma_j})=F_{\beta_{j+1}}$ for all $j\in\I_{\ell-1}$.
As $\gamma_{\ell}=\alpha_i$,
\begin{align*}
\underline{E}_{\gamma_{\ell}}&=f\, \underline{E}_{i}, & \underline{F}_{\gamma_{\ell}}&=g\, \underline{F}_{i}, & \text{ for some }f,g\in\Acf.
\end{align*}
As $K_{\gamma_\ell}-L_{\gamma_{\ell}}^{-1}=[\underline{E}_{\gamma_{\ell}},\underline{F}_{\gamma_{\ell}}]=fg \, [\underline{E}_i,\underline{F}_i]=K_i-L_i^{-1}$, we have that $fg=1$.

Assume first that $i$ is not a Cartan vertex, that is, $\alpha_i\notin \fO^{\bq}_+$. For each $\beta_j\in\fO^{\bq}_+$ we have that $\gamma_{j-1}\in\fO^{\bp}_+$, $\bqf_{\beta_j}=\bpf_{\gamma_{j-1}}$, and then $E_{\beta_j}^{(N_{\beta_j})}=T_i^{\bqf}\left(\underline{E}_{\gamma_{j-1}}^{(N_{\gamma_{j-1}})}\right)$. After evaluation we get the same equality for $T_i^{\bq}$, so
\begin{align*}
\eb_{\wbeta_j} &= \pi \left(E_{\beta_j}^{(N_{\beta_j})}\right)=\pi \circ T_i^{\bq}\left(\underline{E}_{\gamma_{j-1}}^{(N_{\gamma_{j-1}})}\right)
=\overline{T}_i^{\bq} \circ \pi \left(\underline{E}_{\gamma_{j-1}}^{(N_{\gamma_{j-1}})}\right) = \overline{T}_i^{\bq} (\eb_{\wgamma_{j-1}}).
\end{align*}
Analogously, $\fb_{\wbeta_j}= \overline{T}_i^{\bq} (\fb_{\wgamma_{j-1}})$, so
\begin{align*}
\hb_{\wbeta_j} &= [\eb_{\wbeta_j}, \fb_{\wbeta_j}] = \left[ \overline{T}_i^{\bq} (\eb_{\wgamma_{j-1}}), \overline{T}_i^{\bq} (\fb_{\wgamma_{j-1}}) \right]
= \overline{T}_i^{\bq} \left(  \left[\eb_{\wgamma_{j-1}}, \fb_{\wgamma_{j-1}}\right] \right)
= \overline{T}_i^{\bq} \left( \hb_{\wgamma_{j-1}}\right),
\end{align*}
and the statement follows, as $\eb_{\wgamma_{j-1}}$, $\hb_{\wgamma_{j-1}}$, $\fb_{\wgamma_{j-1}}$, $\gamma_{j-1}\in\fO^{\bp}_+$, generate $\g_{\rho_i \bq}$.

Now assume that $i$ is a Cartan vertex, that is, $\alpha_i\in \fO^{\bq}_+$. For each $\beta_j\in\fO^{\bq}_+$, $\beta_j\ne \alpha_i$, the proof follows as in the previous case, so consider the case $\beta_1=\alpha_i$. Here,
\begin{align*}
T_i^{\bqf}\left(E_i^{(N_i)}\right) &= \tfrac{1}{(N_i1)_{\bqf}^!}(F_iL_i)^{N_i}
 = \bqf_{ii}^{\binom{N_i+1}{2}} \tfrac{1}{(N_i)_{\bqf}^!}L_i^{N_i}F_i^{N_i}
=\bqf_{ii}^{\binom{N_i+1}{2}} L_i^{N_i}F_i^{(N_i)},
\end{align*}
so after evaluation, $T_i^{\bq}\left(E_i^{(N_i)}\right)=\bq_{ii}^{\binom{N_i+1}{2}} L_i^{N_i}F_i^{(N_i)}=\pm F_i^{(N_i)}$. Hence,

\begin{align*}
\overline{T}_i^{\bq} (\eb_{\wgamma_{\ell}}) &= f^{N_i} T_i^{\bq}\circ \pi \left(E_i^{(N_i)}\right) 
= f^{N_i} \pi \circ T_i^{\bq} \left(E_i^{(N_i)}\right) = \pm f^{N_i} \pi \left( F_i^{(N_i)} \right) = \pm f^{N_i} \fb_{\beta_1}.
\end{align*}
Also, $\overline{T}_i^{\bq} (\fb_{\wgamma_{\ell}}) =\pm g^{N_i} \eb_{\beta_1}$, and combining these two computations,
$\overline{T}_i^{\bq} (\hb_{\wgamma_{\ell}}) =- \hb_{\beta_1}$.

Hence, in both cases we get that $\overline{T}_i^{\bq}(\g_{\rho_i (\bq)})=\g_{\bq}$ and the statement follows.
\epf

\begin{remark}
The map induced by $T_i^{\bq}$ between $U(\g_{\rho_i(\bq)})$ and $U(\g_{\bq})$ can be thought as the identity if $i$ is not Cartan, as we may identify $\varPi^{\rho_i(\bq)}$ with $\varPi^{\bq}$ using $s_i^{\bq}$, and as the action of the $i$-th reflection of the Weyl group, up to scalars, if $i$ is Cartan.
\end{remark}

\subsection{An example}\label{subsec:example-B11}
We finish the paper with the presentation by generators and relations of the Hopf algebra $\Ub_{\bq}$ and the associated extension of Hopf algebras when the braiding  is of super type $\superb{1}{1}$. We summarize all the necessary data of this case in Figure \ref{fig:root-data}. First,  we have to consider the algebra $U_{\bqf}$ which is generated by $E_i$, $K_i^{\pm 1}$, $L_i^{\pm 1}$, $F_i$, $i=1,2$, subject to the commutating relations,
\begin{equation}\label{eq:def-rels-superB-integral}
\begin{aligned}
K_iE_j = \bqf_{ij} E_jK_i,&\qquad
L_iE_j = \bqf_{ji} E_jL_i, \\[2pt]
K_iF_j = \bqf_{ij}^{-1} F_jK_i,&\qquad
L_iF_j = \bqf_{ji}^{-1} F_jL_i, \\[2pt]
E_iF_j-F_jE_i &= \delta_{ij}(K_i-L_i^{-1}),
\end{aligned}
\end{equation}
the defining relations of the Nichols algebras $\toba_\bqf$ and $\toba_{\bqf^{(-1)}}$,
\begin{align}\label{eq:rel Nichols bqf}
E_1^2 &=(\ad_c E_2)^3 E_1=0 &F_1^2&=(\ad_c F_2)^3 F_1=0,
\end{align}
and the relations stating that the generators \(K_i^{\pm1}\) and \(L_i^{\pm1}\) generate a free abelian group.

\begin{figure}[ht]
\centering
\begin{tikzpicture}[
>=stealth,
every node/.style={font=\small},
panel/.style={draw=accent!55, rounded corners=4pt, fill=accent!6,
inner xsep=10pt, inner ysep=9pt, align=center,
text width=0.88\linewidth},
title/.style={font=\scriptsize\scshape, text=accent,
fill=white, inner sep=2pt, anchor=west}
]

\node[panel] (p1) {
$\displaystyle
\bqf=\begin{pmatrix} -1 & \nu^{-1}\\ \nu^{-1} & \nu \end{pmatrix}
\quad\xrightarrow{\;\;\nu\,\mapsto\,\xi\;\;}\quad
\bq=\begin{pmatrix} -1 & \xi^{-1}\\ \xi^{-1} & \xi \end{pmatrix}$
};
\node[title] at ([xshift=12pt]p1.north west) {Braiding matrix};

\node[panel, below=5mm of p1] (p2) {
\renewcommand{\arraystretch}{1.6}%
$\begin{array}{r|cccc}
\varDelta^\bq_+ & \alpha_1
& \textcolor{accent}{\alpha_1+\alpha_2}
& \alpha_1+2\alpha_2
& \textcolor{accent}{\alpha_2}\\
\hline
\bq_{\beta,\beta} & -1 & -\xi^{-1} & -1 & \xi\\
N_\beta^\bq & 2
& \textcolor{accent}{%
\displaystyle M=\begin{cases}
\scriptstyle 2N, & \scriptstyle N\text{ odd},\\
\scriptstyle N, & \scriptstyle N=4k,\\
\scriptstyle N/2, & \scriptstyle N=4k+2
\end{cases}}
& 2 & \textcolor{accent}{N}
\end{array}$
};
\node[title] at ([xshift=12pt]p2.north west) {Positive roots};

\node[panel, below=5mm of p2] (p3) {
$\displaystyle
\wfO^{\bq}_+=\varPi^{\bq}=
\bigl\{\textcolor{accent}{M(\alpha_1+\alpha_2)},\,
\textcolor{accent}{N\alpha_2}\bigr\}
\qquad\text{of type } A_1\times A_1$
};
\node[title] at ([xshift=12pt]p3.north west) {Cartan roots};

\end{tikzpicture}
\caption{Braiding of super type $\superb{1}{1}$, for an indeterminate $\nu$ and $\xi$ a primitive root of unity of order $N\geq 3$; cf.~\cite[Table~1, Row~4]{AAR-Belgium}.}
\label{fig:root-data}
\end{figure}

We now construct the non-simple PBW generators of $U_{\bqf}$,
\begin{align*}
E_{12} &= E_1E_2-\nu^{-1}E_2E_1, & 
E_{122} &= E_{12}E_2-E_2E_{12}, 
\\
F_{12} &=\frac{F_1F_2-\nu F_2F_1}{\nu^{-2}-1}, & 
F_{122} &= \frac{F_{12}F_2-F_2F_{12}}{(\nu-\nu^{-1})^{2}},
\end{align*}
and the quantum brackets and quantum divided powers associated to the Cartan roots,
\begin{align}\label{eq:B(1|1)-generators-Cartan-roots}
&E_2^{(N)}, &&E_{12}^{(M)},&& \hbt_2\coloneqq \bra{2}{N},&& \hbt_{12} \coloneqq \bra{12}{M},&& F_2^{(N)}, && F_{12}^{(M)}.
\end{align}
Here, $\Acf=\mathbb Z\bigl[\nu^{\pm1},(\nu-1)^{-1},(\nu^2-1)^{-1}\bigr]$. We consider the $\Acf$-subalgebra $\Ub_{\bqf,\Acf}$ generated by these quantum elements and the generators of $U_{\bqf}$, and following \eqref{eq:U bq},
\begin{align*}
\Ub_{\bq}= \Ub_{\bqf,\Acf}\otimes_{\Acf} \Bbbk/ \langle K_2^N-1, K_1^MK_2^M-1 \rangle.
\end{align*}
Thus $\Ub_{\bq}$ is generated by the elements $E_i$, $K_i$, $L_i$, $F_i$, with $i\in\{1,2\}$, and those in \eqref{eq:B(1|1)-generators-Cartan-roots}, subject to the relations in $U_{\bqf}$ and the relations \eqref{eq:rel B11 0 part}--\eqref{eq:last relation} which we describe below.

The relations among the elements in the $0$-part $\Ub_{\bq}^0$ follow from Lemma \ref{lem:powers-Kbeta-Lbeta}. Namely,
\begin{align}\label{eq:rel B11 0 part}
\begin{aligned}
K_2^N &= 1,  &  L_2^N&=1, & \bra{2}{mN}&=\frac{1}{m!}\prod_{k=0}^{m-1}\big(\hbt_2+k(\xi-1)^N\big)
\\
K_{1}^MK_2^M &= 1, & L_{1}^ML_2^M&=1, & \bra{12}{mM}&=\frac{1}{m!}\prod_{k=0}^{m-1}\big(\hbt_{12}+k(-\xi^{-1}-1)^M\big).
\end{aligned}
\end{align}
By Lemma \ref{le:U-0-eval-A-repr} and \eqref{eq:U0q = Gammaq ot Uhq},
\begin{align*}
\Ub_{\bq}^0 \simeq \Bbbk( \Z/N\Z \times \Z/M\Z)^2 \otimes \Bbbk[x_1,x_{2}].
\end{align*}

The relations between the generators of $\Ub_{\bq}^0$ and those of $\Ub_{\bq}^{\pm}$ come from \eqref{eq:KiLi-commutes-Fbeta-Ebeta}--\eqref{eq:bracket-commutes-Fbeta-Ebeta}, together with those in \eqref{eq:def-rels-superB-integral}. Explicitly,
\begin{align}
K_iE_{\gamma}^{(N_\gamma^\bq)} &= E_{\gamma}^{(N_\gamma^\bq)}K_i,&
K_iF_{\gamma}^{(N_\gamma^\bq)} &= F_{\gamma}^{(N_\gamma^\bq)}K_i,
\\
L_iE_{\gamma}^{(N_\gamma^\bq)} &= E_{\gamma}^{(N_\gamma^\bq)}L_i,&
L_iF_{\gamma}^{(N_\gamma^\bq)} &= F_{\gamma}^{(N_\gamma^\bq)}L_i,
\\
\hbt_2E_1&=E_1\bracket{2}{2}{N},& 
\hbt_2F_1&=F_1\bracket{2}{-2}{N},
\\
\hbt_2E_2&=E_2\bracket{2}{-2}{N},&
\hbt_{2}E_{2}^{(N)} &=E_{2}^{(N)}\hbt_{2}+2(\xi-1)^NE_{2}^{(N)},
\\
\hbt_2F_2&=F_2\bracket{2}{2}{N},&
\hbt_{2}F_{2}^{(N)} &=F_{2}^{(N)}\hbt_{2}-2(\xi-1)^NF_{2}^{(N)},
\\
\hbt_{2}E_{12}^{(M)} &=E_{12}^{(M)}\hbt_{2}, &
\hbt_{2}F_{12}^{(M)} &=F_{12}^{(M)}\hbt_{2},
\\
\hbt_{12}E_1 &= E_1\bracket{12}{-2}{M}, &
\hbt_{12}F_1 &= F_1\bracket{12}{2}{M},
\\
\hbt_{12}E_2 &= E_2\hbt_{12},&
\hbt_{12}F_2 &= F_2\hbt_{12},&
\\
\hbt_{12}E_{2}^{(N)} &=E_{2}^{(N)}\hbt_{12},&
\hbt_{12}E_{12}^{(M)} &=E_{12}^{(M)}\hbt_{12}+2(-1-\xi^{-1})^ME_{12}^{(M)},
\\
\hbt_{12}F_{2}^{(N)} &=F_{2}^{(N)}\hbt_{12}, &
\hbt_{12}F_{12}^{(M)} &=F_{12}^{(M)}\hbt_{12}-2(-1-\xi^{-1})^MF_{12}^{(M)},
\end{align}
for $\wgamma\in\varPi^{\bq}$; the quantum brackets $\bracket{2}{\pm2}{N}$ and $\bracket{12}{\pm2}{M}$ can be expressed as linear combinations of powers of $\hbt_{2}$, $\hbt_{12}$, $K_2$ and $K_{12}=K_1K_2$ using \eqref{eq:bracket-n<0} and \eqref{eq:bracket-nge0}.

The relations between generators of $\Ub_{\bq}^+$, respectively $\Ub_{\bq}^-$, besides those in  \eqref{eq:def-rels-superB-integral}, are the following. The defining relations of the Nichols algebras $\toba_\bq$ and $\toba_{\bq^{(-1)}}$, besides those in \eqref{eq:rel Nichols bqf}, are
\begin{align}
E_2^N &= E_{12}^M=0, & F_2^N &=F_{12}^M=0.
\end{align}
Corollary \ref{coro:convex-order-divided-powers} gives a formula for the bracket of $F_2^{(N)}$ and $F_{12}^{(M)}$ before evaluation, and analogous formulas hold for $E_2^{(N)}$ and $E_{12}^{(M)}$. By degree constraints we get
\begin{align}
E_2^{(N)}E_{12}^{(M)}&=E_{12}^{(M)}E_2^{(N)}, & F_2^{(N)}F_{12}^{(M)}&=F_{12}^{(M)}F_2^{(N)}.
\end{align}
We can also check the following formulas:
\begin{align}
E_1E_{12}^{(M)} &= E_{12}^{(M)} E_1, & E_1E_2^{(N)} &= \xi^{-2}E_2^{(N-2)}E_{122} + \xi^{-1}E_2^{(N-1)}E_{12} + E_2^{(N)}E_1
\\
E_2E_2^{(N)}&=E_2^{(N)}E_2, & E_{12}^{(M)}E_2 &= E_{12}^{(M-1)}E_{122} + E_2E_{12}^{(M)},
\\
F_1F_{12}^{(M)} &= F_{12}^{(M)} F_1, & F_1F_2^{(N)} &= \xi^{2}F_2^{(N-2)}F_{122} + \xi F_2^{(N-1)}F_{12} + F_2^{(N)}F_1
\\
F_2F_2^{(N)}&=F_2^{(N)}F_2, & F_{12}^{(M)}F_2 &= F_{12}^{(M-1)}F_{122} + F_2F_{12}^{(M)}.
\end{align}

Finally we give relations between generators of $\Ub_{\bq}^+$ and $\Ub_{\bq}^-$, besides the last one in  \eqref{eq:def-rels-superB-integral}. From Lemma \ref{lemmma:Ealpha-M-Falpha-N} we get, after evaluation,
\begin{align}
E_2^{(N)}F_2 &= F_2E_2^{(N)}+(\xi K_2-L_2^{-1})E_2^{(N-1)},
\\
E_2F_2^{(N)} &= F_2^{(N)}E_2+F_2^{(N-1)}(\xi K_2-L_2^{-1}),
\\
E_2^{(N)}F_2^{(N)} &=F_2^{(N)}E_2^{(N)}+\sum_{t=1}^{N} F_2^{(N-t)} \bracket{2}{-t-1}{t}E_2^{(N-t)},
\\
E_{12}^{(M)}F_{12}^{(M)} &=F_{12}^{(M)}E_{12}^{(M)}+\sum_{t=1}^{M} F_{12}^{(M-t)} \bracket{12}{-t-1}{t}E_{12}^{(M-t)}.
\end{align}

By induction we get relations involving $E_i$, $F_i$ and the  divided powers $E_2^{(n)}$, $F_{12}^{(m)}$, $E_{12}^{(m)}$ and $F_2^{(n)}$, before evaluation. After evaluation we have:
\begin{align}
E_1F_2^{(N)} &= F_2^{(N)}E_1, \qquad E_2^{(N)}F_1 = F_1E_2^{(N)},
\\
E_1F_{12}^{(M)} &= (1-\xi^{-2})F_{122}F_{12}^{(M-2)}L_1^{-1}+(1-\xi^{-2})F_{2}F_{12}^{(M-1)}L_1^{-1}+F_{12}^{(M)}E_1,
\\
E_{12}^{(M)}F_1 &= (1-\xi^{-2})K_1^{-1}E_{122}E_{12}^{(M-2)}+(1-\xi^{-2})K_1^{-1}E_{2}E_{12}^{(M-1)}+F_1E_{12}^{(M)},
\\
E_2F_{12}^{(M)} &=(1-\xi^{-2})K_2F_1F_{12}^{(M-1)}+F_{12}^{(M)}E_2,
\\
E_{12}^{(M)}F_{2} &=(1-\xi^{-2})E_{12}^{(M-1)}E_1K_2+F_{2}E_{12}^{(M)},
\\
E_2^{(N)}F_{12}^{(M)} &=(1-\xi^{-2})K_2F_1F_{12}^{(M-1)}E_2^{(N-1)}+F_{12}^{(M)}E_2^{(N)},
\\
\label{eq:last relation}
E_{12}^{(M)}F_{2}^{(N)} &=(1-\xi^{-2})F_{2}^{(N-1)}E_{12}^{(M-1)}E_1K_2+F_{2}^{(N)}E_{12}^{(M)}.
\end{align}
The extension of Hopf algebras of Theorem \ref{thm:Uq-extension-Hopf} is in this case
\begin{equation*}
	\begin{tikzcd}[column sep=large, arrow style=tikz, >=To]
		\uf_{\bq} \arrow[r, hook] & \Ub_{\bq} \arrow[r, two heads] & U(\mathfrak{sl}_2\oplus\mathfrak{sl}_2).
	\end{tikzcd}
\end{equation*}
An explicit PBW basis of $\Ub_{\bq}$ is
\begin{align*}
& \left\{ E_2^{(m_1)} E_{122}^{(m_2)} E_{12}^{(m_3)} E_1^{(m_4)} K_{12}^{p_1}K_2^{p_2}L_{12}^{q_1}L_2^{q_2} \hbt_{12}^{r_1}\hbt_{2}^{r_2}  F_2^{(n_1)} F_{122}^{(n_2)} F_{12}^{(n_3)} F_1^{(n_4)} |  \right. \\
& \left. 0\le m_1,m_3, n_1,n_3, r_1, r_2 <\infty, 0\le m_2,m_4, n_2,n_4 <2, 0 \le p_1,q_1<M, 0\le p_2,q_2<N \right\}.
\end{align*}

\end{document}